\RequirePackage{fix-cm}
\documentclass[smallextended]{svjour3}       % onecolumn (second format)
\smartqed  % flush right qed marks, e.g. at end of proof
\usepackage{graphicx}
\usepackage{fullpage}
\usepackage{amsmath}
\usepackage{amsfonts}
\usepackage{amssymb}
\usepackage{hyperref}
\usepackage[]{algorithm2e}
\usepackage{comment}
\usepackage{tabularx}
\usepackage{todonotes}
\usepackage{subcaption}
\usepackage{float}
\usepackage{tikz}
\usepackage{color}
\usepackage{wrapfig}

\usepackage{mathtools}

\DeclarePairedDelimiter\abs{\lvert}{\rvert}%
\DeclarePairedDelimiter\norm{\lVert}{\rVert}%   

\makeatletter
\let\oldabs\abs
\def\abs{\@ifstar{\oldabs}{\oldabs*}}
\let\oldnorm\norm
\def\norm{\@ifstar{\oldnorm}{\oldnorm*}}
\makeatother

\begin{document}

\title{Neural network based limiter with transfer learning
}


\author{R\'{e}mi Abgrall         \and
        Maria Han Veiga %etc.
}


\institute{Maria Han Veiga \at
              University of Zurich \\
              Tel.: +41 44 63 56191\\
              %Fax: +00000\\
              \email{hmaria@physik.uzh.ch}           %  \\
           \and
           R\'{e}mi Abgrall \at
           University of Zurich\\
           \email{remi.abgrall@math.uzh.ch}
}

\date{Received: date / Accepted: date}

\maketitle

\begin{abstract}
Inspired by \cite{deepray2018,mhvabgrall}, a neural network is trained using simulation data from a Runge Kutta discontinuous Galerkin (RKDG) method and a modal high order limiter \cite{lydia}. With this methodology, we design one and two-dimensional black-box shock detection functions. Furthermore, we describe a strategy to adapt the shock detection function to different numerical schemes without the need of a full training cycle and large dataset.  We evaluate the performance of the neural network on a RKDG scheme for validation. To evaluate the domain adaptation properties of this neural network limiter, our methodology is verified on a residual distribution scheme (RDS), both in one and two-dimensional problems, and on Cartesian and unstructured meshes. Lastly, we report on the quality of the numerical solutions when using a neural based shock detection method, in comparison to more traditional limiters, as well as on the computational impact of using this method in existing codes.
\end{abstract}

\section{Introduction}
\label{introduction}
When dealing with nonlinear conservation laws \eqref{eq:conslaw}, it is well known that discontinuous solutions can emerge, even for smooth initial data. The numerical approximation of the discontinuous solution will develop non-physical oscillations around the discontinuity, which in turn will negatively impact the accuracy of the numerical scheme. There are different methods to control these oscillations, for example, adding a viscous term as in \eqref{eq:viscous} or limiting the solution.

\begin{equation}
\label{eq:conslaw}
\frac{\partial}{\partial t} u + \nabla \cdot f(u) = 0
\end{equation}

\begin{equation}
\label{eq:viscous}
\frac{\partial}{\partial t} u + \nabla \cdot f(u) = \nabla \cdot (\nu(u)\nabla u)
\end{equation}

Neural networks gained new popularity due to the computational tractability of the back-propagation algorithm, used for the learning of weights and biases in a deep neural network. Furthermore, it has been empirically shown to generate robust models for classification in many areas of applications \cite{vision,s2s} and theoretically, to generate \textit{universal} classifiers and \textit{universal} function approximators \cite{rojas,rojas2,ilya}. 

One can write a stabilisation method as a function $f$ that takes some local properties of the solution (denoted by the map $X(u(x))$), and returns a modified solution value $\tilde{u}(x)$, which has some desired properties, e.g. non-oscillatory, maximum principle preserving,

\[f \left( u(x), X\left(u(x)\right)\right) = \tilde{u}(x). \]

For example, in the context of finite volume schemes, where the usual notation applies: $\bar{u}_j$ denotes the cell average of the solution $u(x)$ in the interval $[x_{j-1/2},x_{j+1/2}]$, $u_j$ the pointwise solution value at node $x_j$, we can write the \textbf{Minmod limiter} as:

\[ f \left( u(x), X\left(u(x)\right)\right) = \bar{u}_j + \mbox{minmod}(X(u(x)))\frac{(x-x_j)}{x_{j+1/2}-x_{j-1/2}} = \tilde{u}(x) \quad x\in[x_{j-1/2},x_{j+1/2}]\]
where the map $X(u(x))$ encodes local properties of the solution around $[x_{j-1/2},x_{j+1/2}]$. In this case, 
\[X(u(x)) = (\bar{u}_j - u^+_{j-1/2}, \Delta^+ \bar{u}_j,\Delta^-\bar{u}_j),\]
where $\Delta^+\bar{u}_j := \bar{u}_{j+1}-\bar{u}_j$, $\Delta^- \bar{u}_j:= \bar{u}_j-\bar{u}_{j-1}$. It is well known that this limiter clips the extrema for smooth solutions.

To minimize the clipping of extrema for smooth solutions, one can consider a \textbf{TVD limiter} \cite{tvdlimiter}:
\[ f \left( u(x), X\left(u(x)\right)\right), M) = \bar{u}_j - \mbox{TVD}(X(u(x)),M)\frac{(x-x_j)}{x_{j+1/2}-x_{j-1/2}} = \tilde{u}(x) \quad x\in[x_{j-1/2},x_{j+1/2}], \]
where the map $X(u(x))$ again encodes local properties of the solution around $[x_{j-1/2},x_{j+1/2}]$, and
\[
TVD(X(u(x)),M) =  \left\{
\begin{array}{ll}
      \mbox{minmod}(X(u(x))) & \mbox{else} \\
      \Delta\bar{u}_j & \mbox{if } |\bar{u}_j - u^+_{j-1/2}| \leq M(\Delta x)^2,
\end{array}
\right.
\]
where M is a user defined parameter that gives an estimate of the smoothness of the solution $u(x)$. It is to note that $M$ can take a value in a large range of positive numbers, and that it is usually a global quantity fixed in the beginning of the numerical experiment. Thus, this can be a drawback if the solution has different smoothness properties across the domain (in space and time).

In the context of discontinuous Galerkin method, we can mention, among many, the high-order limiter \cite{biswas1994,lydia}, which does the limiting in a hierarchical manner. \cite{lydia}, in particular, relies on the modal representation of the solution $u(x)$, it is formulated specifically for the modal discontinuous Galerkin method using Legendre polynomials as basis functions, and it does not clip the solution maxima.

Recently, the idea of using artificial neural networks as trouble cell indicators has been explored firstly in \cite{deepray2018}, motivated by the objective to find a universal troubled-cell indicator that can be used for general conservation laws. The authors show that this type of approach is promising, beating traditional, parameter-dependent limiters.

Domain adaptation is concerned with using a model from a particular source data distribution on a different (but related) target data distribution \cite{transfer}. One simple example is the task of \emph{spam filtering}, where a model is used to discriminate between spam and non-spam emails. A model can be trained on the data of a particular user, and adapted to be used on the data of a new user, who might receive significantly different emails. This can be a strategy to train limiters that work in different numerical schemes, while using only (or the majority of the) labeled data from a particular numerical scheme. This can be useful, for example, in the case where one has access to a solver for which labeled data is easy to obtain. In particular, we are interested in training a shock detection algorithm using data from a simple Runge Kutta discontinuous Galerkin method, and adapting this limiter to different numerical schemes.

In this work we are interested in addressing the following questions:
\begin{enumerate}
\item{Is it possible to learn a data-driven stabilisation function which requires minimal user input?}
\item{Can this stabilisation algorithm be used in different numerical schemes, leading to stabilisation methods which are agnostic to the underlying numerical scheme?}
\end{enumerate}

The paper is structured as follows: the methodology of training and using a trained limiter is presented in section \ref{method} (as well as the extension to 2-dimensional problems), in section \ref{sec:transfer} we describe the transfer learning strategy and in section \ref{sec:dataset} the construction of the dataset is described. In sections \ref{sec:1ddetection} and \ref{sec:results2d}, numerical results for 1-dimensional and 2-dimensional problems, respectively, are shown. We conclude the paper with a discussion and outlook, in section \ref{sec:conclusion}.

In the spirit of open and reproducible science, all the datasets, trained models and some solvers\footnote{some of the solvers used are still under development and not publicly available.} are made available in a public repository \cite{repository}. 

\section{Method}
\label{method}
In the following section, we focus on the three main aspects of this method:
\begin{itemize}
\item The set-up of the learning algorithm (section \ref{sec:setup}),
\item The integration of a trained model with an existing CFD code (section \ref{sec:trainedmodel}),
\item A description of the performance measures used to validate this method (section \ref{sec:perfmeasures}).
\end{itemize}

\subsection{Set-up of the learning algorithm}
\label{sec:setup}
In this section we describe the details of the learning algorithm. We use a multilayer Perceptron (MLP) neural network. A general comprehensive introduction of neural networks can be found in \cite{deeplearning}.

We wish to learn a map $f: \bf{X} \to \mathcal{Y}$, where $\bf{X}$ denotes an arbitrary set containing examples that we wish to label with possible outcomes $\mathcal{Y}$. The task at hand is a binary classification, thus $f$ will be our binary classifier and $\mathcal{Y} = \{0,1\}$.

We choose $f$ to be defined by the composition of a sequence of functions $g_1, g_2, ... g_n$, yielding

\[f(x) = g_n( ... g_2(g_1(x))).\]

This is known as a deep neural network. There are many different classifiers which can be used, but it has been shown that deep neural networks perform well on a variety of classification tasks, in particular when the classification plane is nonlinear \cite{deeplearning}.

Each function $g_i(w_i,b_i,h_i(\cdot))$ is parametrized by a matrix $w_i$, called weights, a vector $b_i$ called bias and an activation function $h_i(\cdot)$ which introduces the non-linearity on the neural network.

These parameters are tuned through a \emph{loss function} $\mathcal{L}(x)$, which measures how well the mapping $f$ performs on a given dataset $\mathcal{D}$, using the gradient descent and back-propagation algorithms.

The gradient descent \cite{gradientdescent} is a first-order iterative optimization algorithm for finding the minimum of a function (in this case, the loss function $\mathcal{L}(x)$), relying on the fact that for small enough $\eta$:
\[ a_{n+1} = a_n - \eta \nabla \mathcal{L}(a_n) \]
then $\mathcal{L}(a_{n+1})\leq \mathcal{L}(a_n)$, for a defined and differentiable loss function $\mathcal{L}$. The size of the update on the opposite direction of the steepest gradient is tuned through $\eta$, the \emph{learning rate}. Instead of using a global learning rate, we use the Adam algorithm \cite{adam}, which does not suffer from the problem of choosing the wrong learning rate, by adaptively estimating the size of the update for each parameter (in this case, the weights and biases), using information from the first and second moments of the gradient $\nabla \mathcal{L}(a_n)$ (with respect to parameter $a$). The update follows:

\[ a_{n+1} = a_n - \alpha \frac{\hat{m}_n}{\sqrt{\hat{v}_n} + \varepsilon}, \]
where $\alpha$ denotes the \textit{stepsize}, typically set to $\alpha = 0.001$, $\hat{m}_t$ denotes the \textit{corrected} moving first moment of the gradient, $\hat{v}_t$ the \textit{corrected} moving second moment of the gradient (both with respect to parameter $a$) and $\varepsilon$ a small value to avoid division by zero. For the full algorithm, please refer to appendix \ref{nn-ap:adam}.

Finally, the back-propagation algorithm is used, which gives a rule to update the different $w_i$ and $b_i$ in order to minimize the loss function \cite{backprop}.

Furthermore, two different loss functions are used:

\begin{itemize}
\item{cross-entropy:

\begin{equation}
\label{eq:loss}
\mathcal{L}(\mathcal{D}) = - \frac{1}{N}\sum_i^N y_i\log(\hat{p}_i) + (1 - y_i)\log(1 - \hat{p}_i).
\end{equation}}
\item{weighted cross-entropy:
\begin{equation}
\label{eq:loss}
\mathcal{L}(\mathcal{D}) = - \frac{1}{N}\sum_i^N y_i\log(\hat{p}_i)\omega + (1 - y_i)\log(1 - \hat{p}_i).
\end{equation}}
\end{itemize}

We consider the weighted cross entropy loss function as one particular property of our problem is that it is expected that there will be a class imbalance on the dataset (both during the training phase and prediction phase). In particular, it is more likely to find cells which are in no need for stabilisation than ones which are in need for stabilisation. Furthermore, it is more desirable to overlimit than to miss a cell that needs limiting, as it might lead to unphysical results and potentially crash the code. To bias the learning towards predicting the less represented class (needing stabilisation), it is common practice to use a weighted cost function, which increases the penalty of mislabeling a positive label\cite{imbalance}. The asymmetry in the loss function is added through the coefficient $\omega$.

Finally, we specify the activation functions used: for the hidden layers $1,...,n-1$ a Rectified Linear Unit (ReLU) is used:
\[h(x) = \max(0,x). \]

Although there exist more sophisticated activation functions, typically modifications to ReLU, e.g. leaky ReLU, Parametric ReLU or Randomised leaky ReLU, these require further parameter estimation, add at least one more dependence to the saved model, and the empirical improvement on the performance is not extremely significant \cite{relus}. Furthermore, there are theoretical predictions when using ReLU as an activation function, which we make use in the following section to roughly estimate the size of the neural network.

For the last layer (output layer), a sigmoid function is used:

\[ h(x) = \frac{1}{1 + \exp ^{-x}},\]
in order to attain a value that can be interpreted as a probability.

Lastly, the datasets (described in section \ref{sec:dataset}) are split into disjoint sets of training, validation and test sets.

The training phase is detailed in algorithm \ref{algo:train}. We set the batchsize to be $256$, and the number of epochs to be between $4$ and $10$, the Adam algorithm \cite{adam} as the optimizer and use a stopping criteria to avoid overfitting, that is triggered when the empirical generalisation error (measured through the loss on the validation set) increases for several training cycles \cite{earlystopping}. 

\vspace{10pt}

\begin{algorithm}[H]
 \KwData{ Training set $\mathcal{T}$, test set $\mathcal{Q}$, hyperparameters: batchsize, nepochs, architecture, optimizer }
 \KwResult{ $\mathcal{C}$ (trained classifier) }
 \For{nepochs}{
 	\For{nsubsets}{
  	training\_loss = $\mathcal{C}(\mathcal{T})$\;
  	gradients = compute loss (training\_loss )\;
  	$\mathcal{C}$.update\_weights(gradients)\;
  	test\_loss = $\mathcal{C}(\mathcal{Q})$\;
  	\If{(early\_stopping(test\_loss))}{
		stop training\;  		
  		}
  	}
  }
 \caption{Training phase }
 \label{algo:train}
\end{algorithm}

\subsubsection{Architecture}
\label{sec:arch}

Because a fully connected MLP is used, one can follow the results in \cite{optimalbounds}, which establish lower bounds for the number of required non-zero weights and layers necessary for a MLP (using ReLU activation functions), to approximate classification maps $f:[-1/2,1/2]^d\to \mathcal{Y}$:
\[f = \sum_{k\leq N} c_k \mathcal{X}_{K_k},\]
where $f$ is a piecewise smooth function with $N$ ``pieces". Here each $\mathcal{X}_{K_k}$ denotes an indicator function:
\begin{equation}
\label{eq:gausshat}
\mathcal{X}_{K_k}(x) = \left\{
\begin{array}{ll}
      1, & \mbox{if } x\in K_k  \\
	  0, & \mbox{if } x\not\in K_k \\
\end{array} 
\right.  \quad (x,t)\in[0,1]\times \mathbb{R}^+,
\end{equation}
for a compact set $K_k \subset [-1/2,1/2]^d$ with $\partial K_k \in \mathcal{C}^{\beta}$, the set of $\beta$ times differentiable curves, for all $k\leq N$. Namely, $\beta$ controls the regularity of the boundary of the set $K_k$, $d$ denotes the dimension of the feature map.

Theorem 4.2 \cite{optimalbounds} provides the following lower bound for the minimal number of nonzero weights that a neural network $R_\rho(\Phi ^f_\varepsilon)$ needs to have in order to approximate any function belonging to a function class $\mathcal{C} \in L^2([-1/2,1/2]^d)$, up to a $L^2$ error of at most $\varepsilon$:
\[ \# Weights \ge C\cdot \varepsilon^{-2(d-1)/\beta}/\log _2\left(\frac{1}{\varepsilon}\right), \]
where $C=C(d,\beta,B,C_0)$ is a constant, $B$ is an upper bound on the supremum norm of $R_\rho(\Phi ^f_\varepsilon)$ and $C_0>0$ is related to with the maximum bits needed to encode a weight of the network.

Furthermore, Theorem 3.5 \cite{optimalbounds} provides approximation rates for piecewise constant functions. In particular, the number of required layers depends on the dimension $d$ and for the regularity parameter $\beta$:
\[ \# Layers \ge c' \log_2(2+\beta)\cdot(1+\beta/d). \]

However, the assumptions about the classification function's regularity and desired approximation accuracy (in $L^2$ norm) change the shape and size of the network dramatically, leading to very different lower bounds for the number of layers and number of neurons. In practice, we choose networks which are deeper and have a small width ($\#Weights$), as we did not observe a large performance degradation and it lead to less computationally expensive models. In appendix \ref{ap:arch}, the explicit architectures are detailed. 

\subsubsection{Directional invariance}
\label{sec:directionalinv}
\begin{figure}
\begin{center}
\includegraphics[width=7.0cm]{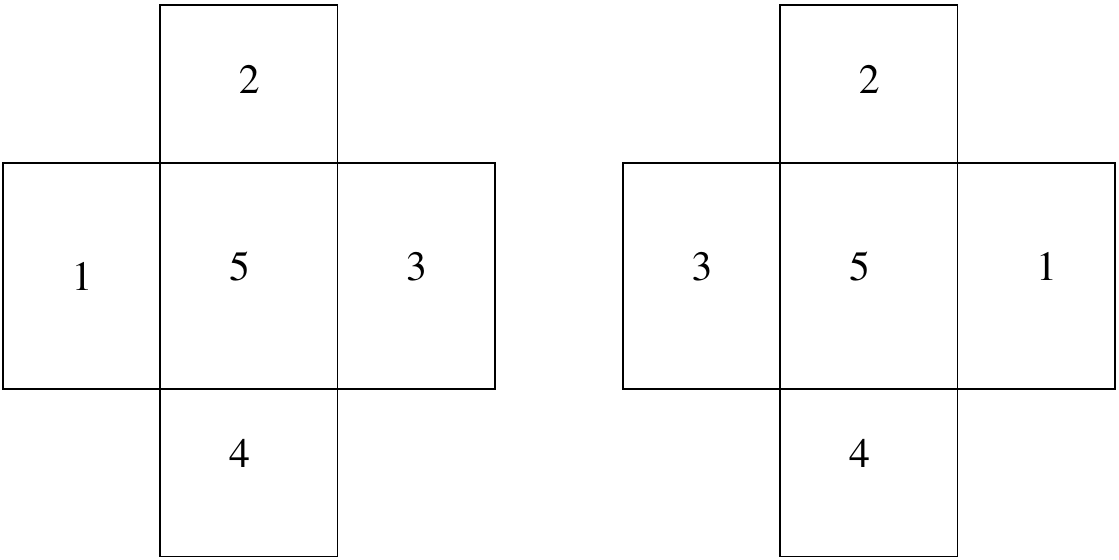}
\caption{Example of desired invariance, with respect to reflection through y-axis.\label{fig:invariance1}}
\end{center}
\end{figure}

To introduce feature invariance (see figure \ref{fig:invariance1}), some models tested use an aggregate way to estimate the label. For a given feature vector $\vec{x}$, several copies of this vector are generated $\vec{x}_1, ... \vec{x}_n$ where permutations between the features are performed and their prediction is evaluated. The final label estimation is given by majority label produced for $\vec{x}, \vec{x}_1, ... \vec{x}_n$. This ensures the response of the model is does not depend on the orientation of the stencil.

\subsection{Integration of the method on a CFD code}
\label{sec:trainedmodel}
We assume that exists a neural network which has been trained offline. Each layer $g_i$ can be fully characterized by the following information:

\[ w_i,b_i,h_i(\cdot) \]
where $w_i$ the weights matrix, $b_i$ a bias vector and $h_i$ denotes the activation function. There are two routines necessary to integrate a trained neural network with an existing code:
\begin{enumerate}
\item{\textbf{Generation of features:} given the local solution $u$, generate the feature quantities $X(u)$.
}
\item{\textbf{Prediction routine:} Given the features $X(u)$, it is necessary to evaluate the function $f$. As noted previously, a layer $g_i$ is fully specified by parameters $w_i, b_i, h_i(\cdot)$. Once the neural network has been trained offline, the weights and biases can be loaded onto the CFD code. What remains to be implemented are the activation functions for the hidden layers and the activation function for the output layer in order to evaluate $f$ at some given input.
}
\end{enumerate}

The generic stabilisation algorithm is detailed in \ref{algo:limiter}.

\begin{algorithm}[H]
 \KwData{solution at cell $i$, $u_i$}
 \KwResult{ stabilised solution $\tilde{u_i}$ }
 X = generateFeatures($u_i$)\;
 label = predict(X)\;
 \eIf{label = 1}{
  $\tilde{u_i}$ = stabilisation($u_i$)\;
  }{
  $\tilde{u_i}$ = $u_i$\;
  }
 \caption{Limiting \label{algo:limiter}}
\end{algorithm}

For systems, each variable is limited independently.

\subsection{Measuring performance of model}
\label{sec:perfmeasures}
In this work we use two sets of performance measures, namely:

\begin{enumerate}
\item Label prediction measures
\item $L^1$-norm of the numerical solution
\end{enumerate}

For the first set of measures, we can use typical metrics used in computational statistics and machine learning communities:
\begin{align*}
\frac{tp + tn}{tp + tn + fp + fn}\quad &\mbox{(accuracy)}\\
\frac{tp}{tp + fn} \quad &\mbox{(recall)}\\
\frac{tp}{tp + fp} \quad &\mbox{(precision)}
\end{align*}
where $tp$ is the number of correctly predicted positive labels, $tn$ the number of correctly predicted negative labels, $fp$ the number of incorrectly predicted positive labels and $fn$ the number of incorrectly predicted negative labels\footnote{Because we are in the regime of a supervised learning problem, we assume true labels exist. The dataset generation is described in section \ref{sec:dataset}.}. It is important to consider recall and precision because the distribution of labels can be imbalanced. If only accuracy is measured, it is possible to attain a high accuracy by solely predicting the majority label. These metrics are computed against the test set.

For the second set of measures, we consider the $L^1$-norm because it is the relevant measure for a CFD code and we can study the effect of this method on the empirical error of the numerical solution.

\section{Dataset}
\label{sec:dataset}

The dataset is an integral part of data-driven analysis. It contains the data for which we want to learn a mapping for. The task is to learn a function which indicates whether a discontinuity is present in the solution or not. The dataset is the set containing $N$ samples $\{X_i, y_i\}_{i=1}^N $, where $X_i$ denotes some local properties of the solution $u_i$ (\textit{features}) and $y_i$ (\textit{labels}) indicates the existence (or not) of a discontinuity.

\begin{figure}
\centering{
\includegraphics[width=5.5cm]{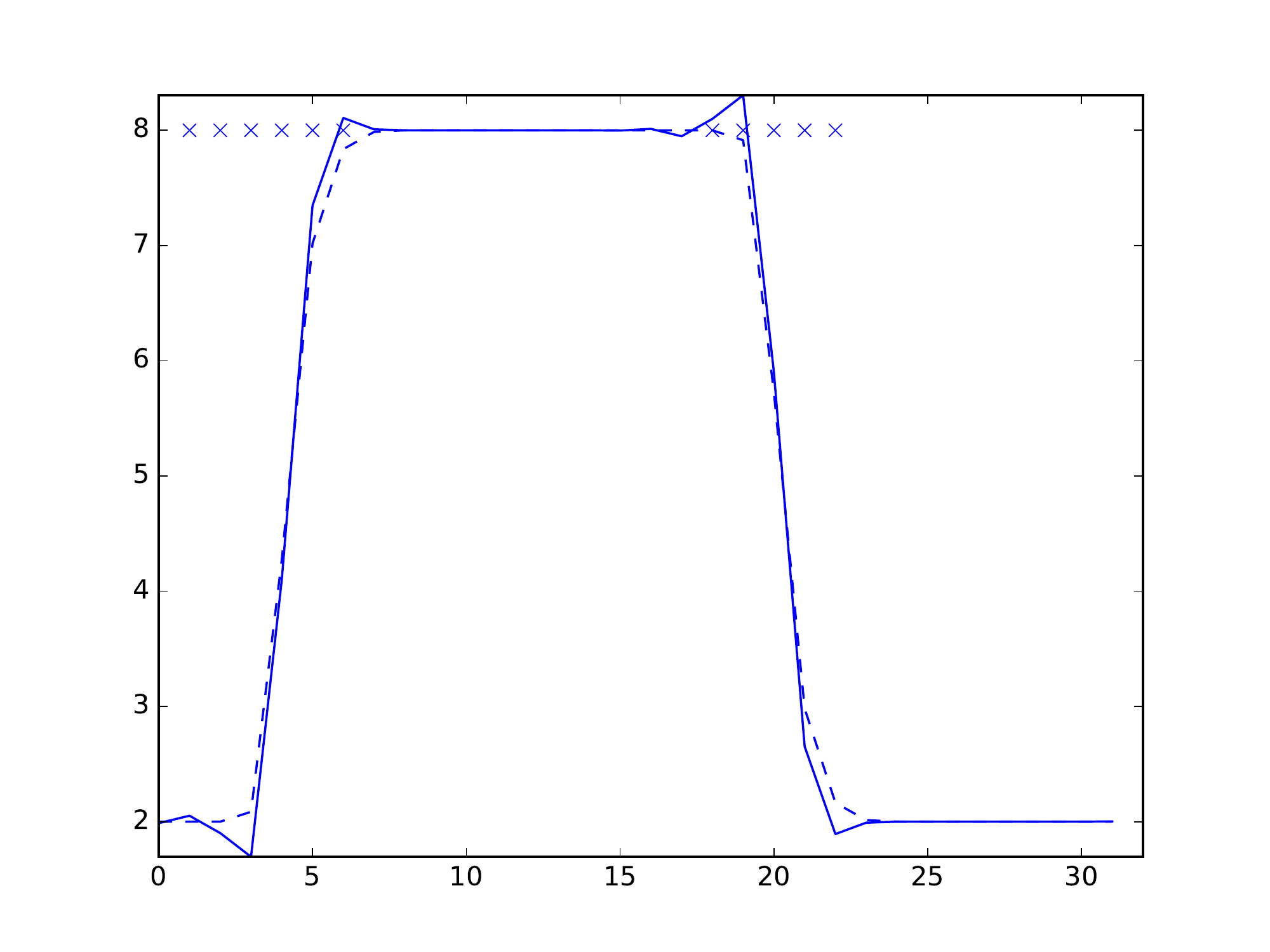}
\caption{Example of a dataset entry, the crosses 'x' denote flagged cells.\label{fig:labelsex1}}}
\end{figure}

For the 1-dimensional case, the dataset is generated by performing many runs on a 1-dimensional discontinuous Galerkin code solving the advection equation for different initial conditions, orders and mesh sizes (see table \ref{table:runs1d}), and the labels are obtained by running a \textit{high-order limiter} for discontinuous Galerkin \cite{lydia} which does hierarchical limiting of the modes (see example in figure \ref{fig:labelsex1}). If the limited solution deviates from the unlimited solution more than $\zeta$, then the cell needs stabilisation, and we get a positive label. We set $\zeta = 0.01$. 

The interesting point of generating the labels in this way is that we can learn a black-box function of a limiter which is specific to a particular framework (in this case, discontinuous Galerkin) and integrate the trained model with other numerical methods. Furthermore, the cell which are flagged by the \textit{high-order limiter} are similar to a well tuned TVD limiter.

\subsection{Features}
The features $X_i$ are the different quantities used to describe the local solution $u_i$. For the sake of generalisation, we choose features which are readily available in different numerical methods, such as: cell mean value, values at interface, divided differences between neighbours (see table \ref{table:features} for the complete description of features). 

Furthermore, to introduce some magnitude invariance, we normalise physical values (such as averages or pointwise values) by applying the following map:
\begin{equation}
\label{eq:normalisation1d}
     u_{normal}(u_*) = \frac{(u_*-u_{min} )}{\mid u_{max} \mid + \mid u_{min} \mid  } - \frac{( u_{max} - u_*)}{\mid u_{max} \mid + \mid u_{min} \mid  }
\end{equation}
where $u_{max} = \max ( \bar{u}_i, \bar{u}_{i+1}, \bar{u}_{i-1})$ and $u_{min} = \min ( \bar{u}_i, \bar{u}_{i+1}, \bar{u}_{i-1})$ (taken only over averages). This is a very important step, as a naive normalisation can lead to a non-informative representation of the feature vector.

Clearly, one could consider features which are specific for a particular method, e.g. the modes of the solution when using a discontinuous Galerkin method, or the smoothness parameters for a WENO method.

\subsection{Extension to 2-dimensional problems}

\begin{figure}
\centering{
\includegraphics[width=3.8cm]{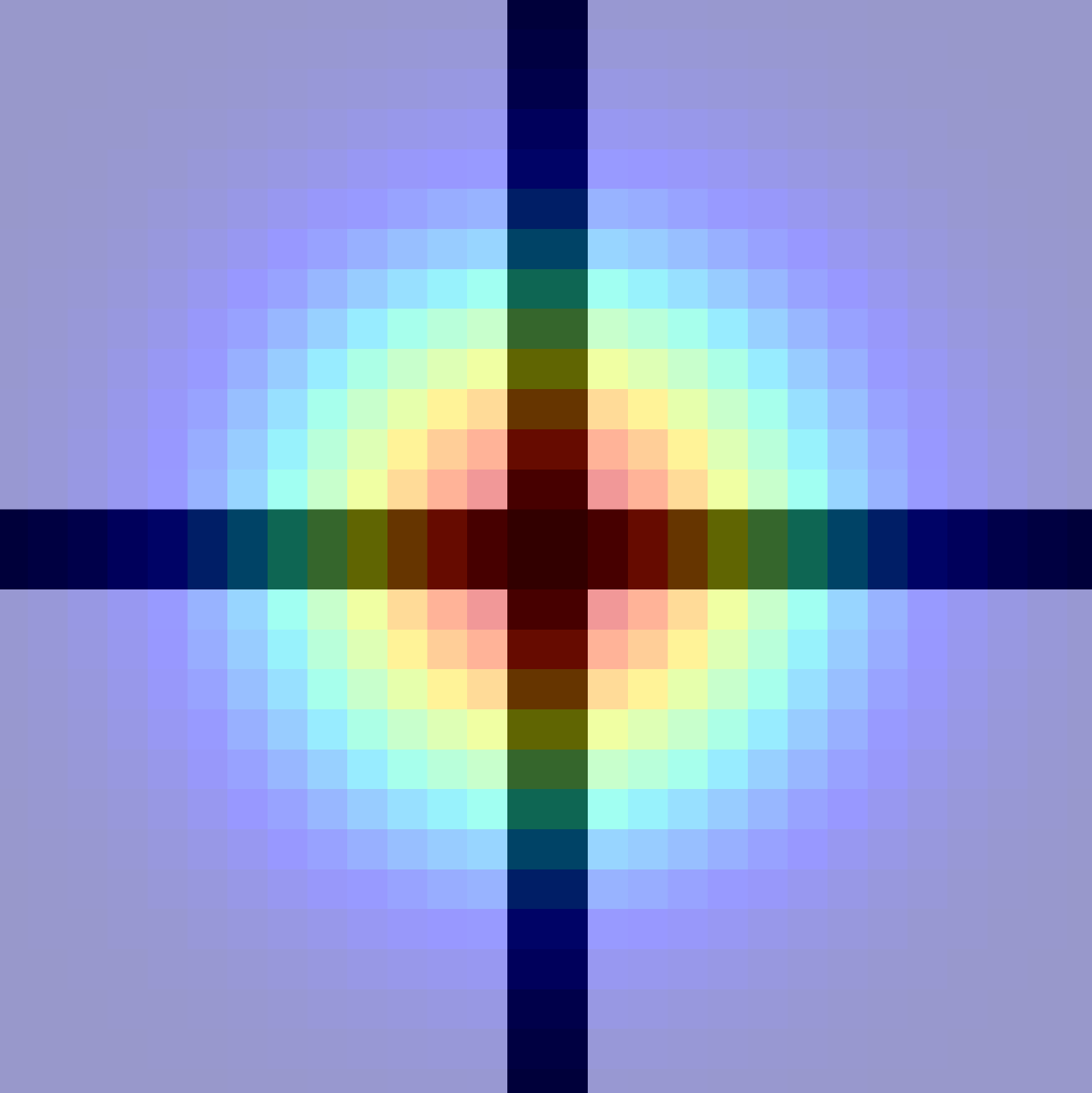}
\caption{Flagging the maximum of a smooth solution.\label{fig:wrongflags}}}
\end{figure}

The extension to 2-dimensional problems is done in a straightforward manner. Using the initial conditions detailed in table \ref{table:runs2d}, we generate a dataset. 

To obtain the labels, the high order limiter\cite{lydia} is used with a threshold $\zeta = 0.25\% $ for maximal difference between the limited and unlimited solution. This is necessary as although the high order limiter might not degrade smooth solutions, it is still triggered (figure \ref{fig:wrongflags}).

Furthermore, the features used are shown on table \ref{table:features2d} and the features are normalised in the same manner as in the 1-dimensional case.
\vspace{10pt}

\begin{table}
\centering{
\resizebox{\columnwidth}{!}{
    \begin{tabular}{ l| l }
   		 \hline
		Initial condition & $\rho = \sin(2\pi x)$,
        $\left\{
          \begin{array}{ll}
                8 & 0.25\leq x \leq 0.75 \\
                2 & \mbox{otherwise} \\
          \end{array} 
          \right.$,
        $\left\{
          \begin{array}{ll}
                \frac{1}{2}\sin(2\pi x) & x \leq 0.3 \\
                0 & \mbox{otherwise} \\
          \end{array} 
          \right. \quad x\in[0,1]$ \\ \hline
        Advection speed & -1, 1 \\ \hline
        Mesh size & 8, 16, 32, 64, 128 \\ \hline
        Order & 2, 3 \\ \hline
    \end{tabular}
    }
	\caption{Runs used to generate the 1-dimensional dataset. \label{table:runs1d}}    
    }
\end{table}

\begin{table}[h]
\centering{
\resizebox{\columnwidth}{!}{
    \begin{tabular}{ l| l l }
      \textbf{ID} & \textbf{Feature Name} & \textbf{Description} \\ \hline
      1 & h & cell width \\ \hline
      2 & $\bar{u}_i$ & average value of solution at cell $i$\\ \hline
      3 & $\bar{u}_{i+1}$ & average value of solution at cell $i+1$ \\ \hline
      4 & $\bar{u}_{i-1}$ & average value of solution at cell $i-1$ \\ \hline
      5 & $u_{i-\frac{1}{2}}^+$ & value of solution at interface $i-1/2$ reconstructed in cell $i$ \\ \hline
      6 & $u_{i+\frac{1}{2}}^-$ & value of solution at interface $i+1/2$ reconstructed in cell $i$\\ \hline
      7 & $u_{i-\frac{1}{2}}^-$ & value of solution at interface $i-1/2$ reconstructed at cell $i-1$ \\ \hline
      8 & $u_{i+\frac{1}{2}}^+$& value of solution at interface $i+1$ reconstructed at cell $i+1$\\ \hline
      9 & $du_{i+1}$ & undivided difference between $\bar{u}_i$ and $\bar{u}_{i+1}$\\ \hline
      10 & $du_{i-1}$ & undivided difference between $\bar{u}_i$ and $\bar{u}_{i-1}$ \\ \hline
      11 & $du_i$ & undivided difference between $\bar{u}_{i+1}$ and $\bar{u}_{i-1}$, divided by $2$ \\ \hline
    \end{tabular}
    }
    \caption{Features table for 1-dimensional problem \label{table:features}}
    }
\end{table}

\begin{table}
\centering{
\resizebox{\columnwidth}{!}{
    \begin{tabular}{ l| l }
   		 \hline
		Initial condition &
        $\left\{
          \begin{array}{ll}
                8 & 0.25\leq \abs{\vec{x}} \leq 0.75 \\
                2 & \mbox{otherwise} \\
          \end{array} 
          \right.$,
        $\exp(-10((x_1-0.5)^2+(x_2-0.5)^2)^2)\quad \vec{x}\in[0,1]^2$ \\ \hline
        Advection speed & $(-1,-1)$, $(1,1)$,  $(1,0)$, $(0,1)$\\ \hline
        Mesh size & 16, 32, 64, 128 \\ \hline
        Order & 2, 3 \\ \hline
    \end{tabular}
    }
    \caption{Runs used to generate the 2-dimensional dataset. \label{table:runs2d}}
    }
\end{table}

\begin{table}
\centering{
\resizebox{\columnwidth}{!}{
    \begin{tabular}{ l| l l }
      \textbf{ID} & \textbf{Feature Name} & \textbf{Description} \\ \hline
      1 & $\Delta x$ & cell x-width \\ \hline
      2 & $\Delta y$ & cell y-width \\ \hline
      3 & $\bar{u}_{i,j}$ & average value of solution at cell $i$\\ \hline
      4 & $\bar{u}_{i+1,j}$ & average value of solution at cell $i+1,j$ \\ \hline
      5 & $\bar{u}_{i-1,j}$ & average value of solution at cell $i-1,j$ \\ \hline
      6 & $\bar{u}_{i,j+1}$ & average value of solution at cell $i,j+1$ \\ \hline
      7 & $\bar{u}_{i,j-1}$ & average value of solution at cell $i,j-1$ \\ \hline
      8 & $u_{i-\frac{1}{2},j}^+$ & value of solution at interface $i-1/2$ reconstructed in cell $i,j$ \\ \hline
      9 & $u_{i+\frac{1}{2},j}^-$ & value of solution at interface $i+1/2$ reconstructed in cell $i,j$\\ \hline
      10 & $u_{i-\frac{1}{2},j}^+$ & value of solution at interface $i-1/2$ reconstructed at cell $i-1,j$ \\ \hline
      11 & $u_{i+\frac{1}{2},j}^-$& value of solution at interface $i+1/2$ reconstructed at cell $i+1,j$\\ \hline
      12 & $u_{i,j-\frac{1}{2}}^+$ & value of solution at interface $j-1/2$ reconstructed in cell $i,j$ \\ \hline
      13 & $u_{i,j+\frac{1}{2}}^-$ & value of solution at interface $i+1/2$ reconstructed in cell $i,j$\\ \hline
      14 & $u_{i,j-\frac{1}{2}}^+$ & value of solution at interface $i-1/2$ reconstructed at cell $i,j-1$ \\ \hline
      15 & $u_{i,j+\frac{1}{2}}^-$& value of solution at interface $i+1/2$ reconstructed at cell $i,j+1$\\ \hline
      16 & $du_{i+1,j}$ & divided difference between $\bar{u}_{i,j}$ and $\bar{u}_{i+1,j}$ \\ \hline
      17 & $du_{i-1,j}$ & divided difference between $\bar{u}_{i,j}$ and $\bar{u}_{i-1,j}$ \\ \hline
      18 & $du_{{i,j}_x}$ & divided difference between $\bar{u}_{i+1,j}$ and $\bar{u}_{i-1,j}$, divided by $2$ \\ \hline
      19 & $du_{i,j+1}$ & divided difference between $\bar{u}_{i,j}$ and $\bar{u}_{i,j+1}$\\ \hline
      20 & $du_{i,j-1}$ & divided difference between $\bar{u}_{i,j}$ and $\bar{u}_{i,j-1}$ \\ \hline
      21 & $du_{{i,j}_y}$ & divided difference between $\bar{u}_{i,j+1}$ and $\bar{u}_{i,j-1}$, divided by $2$ \\ \hline
      22 & $u_{max}$ & maximum value between the averages in the considered patch \\ \hline
      23 & $u_{min}$ & minimum value between the averages in the considered patch \\ \hline
    \end{tabular}}
    \caption{Features table for 2-dimensional problem \label{table:features2d}}
    }
\end{table}

\section{Domain adaptation}
\label{sec:transfer}
With the objective to generate shock detection functions which do not depend on the underlying numerical scheme, we are interested in exploring the idea of using a limiter trained with a particular solver (numerical scheme and mesh type) and observing its performance in a different solver (in particular, for a different numerical scheme and mesh type). Furthermore, we want to explore different strategies which can be used to perform domain adaptation.

The motivations to study this type of problem are: firstly, for some numerical schemes, there are limiters which are designed to be parameter free which rely on a particular feature of the underlying numerical scheme. It would be desirable if limiters designed for a particular numerical scheme could be \textit{generalised} to be used for different numerical schemes. Secondly, there might be a particular numerical solver for which it is easy to generate a lot of labeled data. 

Traditional supervised machine learning operates the assumption that training and testing data are taken from the same input space and the same data distribution. However, this assumption does not always hold. Transfer learning aims to produce an effective model for a target task with limited or no labeled training data, but using knowledge from a different, but related source domain to predict the true label for an unseen target instance.

Formally, the problem of transfer learning can be written as: let $\mathcal{X}_s$ be the \textit{source} instance space. In this space, each instance $x_s \in \mathcal{X}_s$ is represented by a feature vector $\vec{x}_s \in \vec{\mathcal{X}}_s$, where $\vec{\mathcal{X}}_s$ denotes the \textit{source} feature space. Let $\mathcal{X}_t$ be the target instance space and similarly,  $x_t \in \mathcal{X}_t$ is represented by a feature vector $\vec{x}_t \in \vec{\mathcal{X}}_t$, where $\vec{\mathcal{X}}_t$ denotes the \textit{target} feature space. In the case of \textit{heterogeneous transfer learning}, we have that $\vec{\mathcal{X}}_t\neq \vec{\mathcal{X}}_s$. 

Most \textit{heterogeneous transfer learning} solutions involve a transformation of the feature spaces: a symmetric transformation takes both feature spaces $\vec{\mathcal{X}}_t, \vec{\mathcal{X}}_s$ and learns a feature transformation to project each vector onto a common space for adaptation purposes $\vec{\mathcal{X}}_c$ \cite{transfer}, whereas an assymetric transformation transforms the source feature space to align with the target feature space (or vice versa). This approach is appropriate when the source and target have the same class label space and one can transform between $\vec{\mathcal{X}}_t$ and $\vec{\mathcal{X}}_s$.

Concretely, we can consider the problem of transferring a shock detection function trained on a dataset coming from a modal discontinuous Galerkin scheme to be adapted on a dataset generated by a residual distribution scheme (\cite{abg,SWjcp,ricchiuto2010,abgrall2019} for a brief introduction). In particular, the source dataset is generated as in \ref{sec:dataset} and the target dataset is described below, in \ref{subsec:rddata}. Finally, the end goal is to understand whether a limiter designed for one specific class of numerical schemes (in this work, modal Discontinous Galerkin on Cartesian meshes), can be effective on a different numerical scheme (Residual distribution, both for structured and unstructured meshes).

\subsection{Residual distribution dataset}
\label{subsec:rddata}

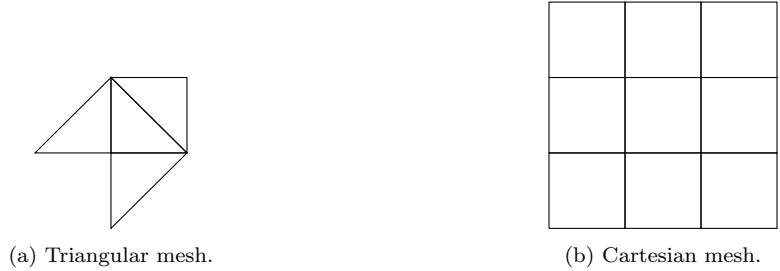
\begin{figure}
     \centering
     \begin{subfigure}[b]{0.45\textwidth}
         \centering
         \begin{tikzpicture}
\draw (0,-1) -- (0,0) -- (1,-1)-- (0,-1);
\draw (1,-1) -- (1,0) -- (0,0)-- (1,-1);
\draw (0,-1) -- (-1,-1) -- (0,0)-- (0,-1);
\draw (0,-1) -- (0,-2) -- (1,-1)-- (0,-1);
\end{tikzpicture}
     \caption{Triangular mesh.}
     \end{subfigure}
     %\hfill
     \centering
     \begin{subfigure}[b]{0.45\textwidth}
     \centering
\begin{tikzpicture}
\draw (-1,0) -- (-1,1) -- (0,1) -- (0,0) -- (-1,0);
\draw (0,0) -- (1,0) -- (1,1) -- (0,1) -- (0,0);
\draw (1,0) -- (2,0) -- (2,1) -- (1,1) -- (1,0); 
\draw (-1,1) -- (-1,2) -- (0,2) -- (0,1) -- (-1,1);
\draw (0,1) -- (1,1) -- (1,2) -- (0,2) -- (0,1);
\draw (1,1) -- (2,1) -- (2,2) -- (1,2) -- (1,1);
\draw (-1,-1) -- (-1,0) -- (0,0) -- (0,-1) -- (-1,-1);
\draw (0,-1) -- (1,-1) -- (1,0) -- (0,0) -- (0,-1);
\draw (1,-1) -- (2,-1) -- (2,0) -- (1,0) -- (1,-1);
\end{tikzpicture}         
      \caption{Cartesian mesh.}
     \end{subfigure}
      \caption{Example meshes.}\label{fig:meshes}
\end{figure}

In order to generate the dataset, small meshes are constructed (see figure \ref{fig:meshes}). We represent the solution in each control volume as linear coefficients of a polynomial basis. Because there is no straightforward way to generate a labeled dataset through simulations, we impose continuous and discontinuous functions, randomly varying different parameters, such as angle of the discontinuity, etc.

Using this method, we can generate a large dataset of examples which is close to the task at hand.

\subsection{1-dimensional case}
We use the trained neural network and we integrate with a 1-dimensional residual distribution code. Note, we must generate \emph{similar enough} features in order to use the model which was trained on a different numerical solver. The 1-dimensional case is not difficult, in particular, because when designing the feature space for the 1-dimensional limiter, we purposely choose quantities which are readily available in most numerical schemes (specifically, quantities based on nodal values, averages and differences).

\begin{comment}
To this end, we use the algorithm detailed in algorithm \ref{algo:rdlimiter1d}.

\begin{algorithm}[H]
 \KwData{solution at cell $i$, $u_i$}
 \KwResult{ stabilised solution $\tilde{u_i}$ }
 X = generateFeatures($u_i$)\;
 label = predict(X)\;
 \eIf{label = 1}{
  $\tilde{u_i}$ = parachute scheme($u_i$)\;
  }{
  $\tilde{u_i}$ = high order scheme($u_i$)\;
  }
 \caption{Limiting \label{algo:rdlimiter1d}}
\end{algorithm}
\end{comment}

\subsection{2-dimensional case}

The two dimensional case is not as simple, as the limiter is trained on Cartesian meshes, and the target problems are defined not only on Cartesian, structured meshes but also on triangular, unstructured meshes.

To this end, we must find a common feature space between these two methods. In particular, we test the two simple strategies:
\begin{itemize}
    \item mapping to feature space of the Cartesian model,
    \item followed by a retraining phase.
\end{itemize}

\subsubsection{Mapping to unstructured mesh}
While the mapping to structured square meshes is quite straightforward, for unstructured meshes we project a triangular element information to the feature space defined in table \ref{table:features2d}. The feature transformation can be found in table \ref{table:remap}.

\subsubsection{Retraining}
A model which was trained with the data generated by the discontinuous Galerkin scheme is loaded and retrained using the dataset as described in Section \ref{subsec:rddata}. %In order to match the feature spaces, we apply the mapping as described in table \ref{table:transform}.

To avoid the phenomena of catastrophic forgetting of neural networks \cite{deeplearning}, which describes the lack of ability to learn different tasks in a sequence, in the retraining phase, a hybrid dataset containing elements from the target and source dataset is used, with a parameter $\lambda$ which determines the ratio to be taken from each dataset.

The retraining algorithm is detailed in algorithm \ref{algo:retraining}.

\begin{algorithm}[H]
 \KwData{trained classifier $\mathcal{C}$, source dataset $\mathcal{S}$, target dataset $\mathcal{T}$, dataset ratio $\lambda$}
 \KwResult{ transferred classifier $\mathcal{C}_t$ }
  $\mathcal{T}_r$ = remapFeatures($\mathcal{T}$)\;
  $\mathcal{D}$ = get\_dataset($\mathcal{S}$,$\mathcal{T}_r$,$\lambda$)\;
  \For{nepochs}{
 	\For{nsubsets}{
  	training\_loss = $\mathcal{C}$.predict($\mathcal{D}_{train}$)\;
  	gradients = compute loss (training\_loss )\;
  	$\mathcal{C}$.update\_weights(gradients)\;
  	test\_loss = $\mathcal{C}$.predict($\mathcal{D}_{test}$)\;
  	\If{(early\_stopping(test\_loss))}{
  		$\mathcal{C}_t = \mathcal{C}$\;
		stop training\;  		
  		}
  	}
  }	 
 \caption{Retraining phase \label{algo:retraining}}
\end{algorithm}

\begin{comment}

\begin{algorithm}[H]
 \KwData{ Training set $\mathcal{T}$, test set $\mathcal{Q}$, validation set $\mathcal{V}$, hyperparameters: batchsize, nepochs, architecture, optimizer }
 \KwResult{ $\mathcal{C}$ (trained classifier) }
 \For{nepochs}{
 	\For{nsubsets}{
  	training\_loss = $\mathcal{C}$.predict($\mathcal{T}$)\;
  	gradients = compute loss (training\_loss )\;
  	$\mathcal{C}$.update\_weights(gradients)\;
  	test\_loss = $\mathcal{C}$.predict($\mathcal{Q}$)\;
  	\If{(early\_stopping(test\_loss))}{
		stop training\;  		
  		}
  	}
  }
 \caption{Training phase }
 \label{algo:train}
\end{algorithm}
\end{comment}

\begin{table}
\begin{tabularx}{\textwidth}{l| l l X}
\textbf{ID} & \textbf{Feature}  & \textbf{Remapped} & \textbf{Description} \\ \hline
      1 & $\Delta x$ & $\sqrt{a}$ & square root of element area \\ \hline
      2 & $\Delta y$ & $\sqrt{a}$ & square root of element area \\ \hline
      3 & $\bar{u}_{i,j}$ & $\bar{u}_{i}$ & average value of solution at triangle $i$\\ \hline
      4 & $\bar{u}_{i+1,j}$ & $\bar{u}_{i,1}$ & average value of solution at neighbouring cell 1 of cell $i$ \\ \hline
      5 & $\bar{u}_{i-1,j}$ & $\bar{u}_{i,2}$ & average value of solution at neighbouring cell 2 of cell $i$ \\ \hline
      6 & $\bar{u}_{i,j+1}$ & $\bar{u}_{i,3}$ & average value of solution at neighbouring cell 3 of cell $i$ \\ \hline
      7 & $\bar{u}_{i,j-1}$ & $\frac{1}{4}(\bar{u}_{i,j} + \sum_k \bar{u}_{k}) $ & average value of solution at patch \\ \hline
      8 & $u_{i-\frac{1}{2},j}^+$ & $u(x_{e_1})$ & value of solution at mid point of edge shared between $i$ and neighbour 1 \\ \hline
      9 & $u_{i+\frac{1}{2},j}^-$ & $u(x_{e_2})$ &value of solution at mid point of edge shared between triangle  $i$ and neighbour 2\\ \hline
      10 & $u_{i-\frac{1}{2},j}^+$ & $u(x_{e_3})$ &value of solution at mid point of edge shared between triangle  $i$ and neighbour 3\\ \hline
      11 & $u_{i+\frac{1}{2},j}^-$& $\bar{u}_{i}$ & average solution $i$ \\ \hline
      12 & $u_{i,j-\frac{1}{2}}^+$ & $u(x_{e_1})$ & value of solution at mid point of edge shared between triangle $i$ and neighbour 1 \\ \hline
      13 & $u_{i,j+\frac{1}{2}}^-$ & $u(x_{e_2})$ &value of solution at mid point of edge shared between triangle  $i$ and neighbour 2\\ \hline
      14 & $u_{i,j-\frac{1}{2}}^+$ & $u(x_{e_3})$ & value of solution at mid point of edge shared between triangle  $i$ and neighbour 3 \\ \hline
      15 & $u_{i,j+\frac{1}{2}}^-$&  $\bar{u}_{i}$ & average solution\\ \hline
      16 & $du_{i+1,j}$ & $\bar{u}_i - \bar{u}_{i,1}$ & undivided difference between $\bar{u}_i$ and $\bar{u}_{i,1}$\\ \hline
      17 & $du_{i-1,j}$ & $\bar{u}_i - \bar{u}_{i,2}$ & undivided difference between $\bar{u}_i$ and $\bar{u}_{i,1}$ \\ \hline
      18 & $du_i$ & $\bar{u}_{i,1} - \bar{u}_{i,2}$  & undivided difference between $\bar{u}_{i,1}$ and $\bar{u}_{i,2}$, divided by $2$ \\ \hline
      19 & $du_{i,j+1}$ & $\bar{u}_i - \bar{u}_{i,3}$ & undivided difference between $\bar{u}_i$ and $\bar{u}_{i,3}$\\ \hline
      20 & $du_{i,j-1}$ & $\bar{u}_i - \bar{u}_i$ & placeholder undivided difference \\ \hline
      21 & $du_j$ & $\bar{u}_{i,1} - \bar{u}_{i,3}$ & undivided difference between $\bar{u}_{i,1}$ and $\bar{u}_{i,3}$, divided by $2$ \\ \hline
\end{tabularx}
\caption{Features transformation \label{table:remap}}
\end{table}
\section{Numerical experiments}
\label{sec:1ddetection}

This section is split in three parts: firstly, we show the performance of several trained neural networks on an unseen validation set by measuring the accuracy, recall and precision. In the second part, we choose a model that performed well and we integrate it with a CFD code. The model runs as a blackbox limiter (denoted as NN) and we compare its performance to the Minmod limiter and hierarchical high order limiter (denoted as HIO) through the $L^1$ error norm. We perform some tests for the linear advection equation and Euler system of equations. The initial conditions are chosen as different from the ones used for the training. Finally, the transferred limiter is tested in the context of a 1-dimensional residual distribution scheme.

\subsection{1-dimensional problems}
\subsubsection{Detection rate}
\label{sub:detection}
We measure the performance of several models on a unseen validation set (table \ref{table:1ddetectionrates}).

\begin{table}

\centering{
\resizebox{\columnwidth}{!}{
\begin{tabular}{ c | c c c c }
      \hline
      \textbf{Model} & \textbf{Description }& \textbf{Accuracy (\%)} & \textbf{Recall (\%)} & \textbf{Precision (\%)} \\ \hline
      Random & randomly guessing & 50.00 & $\theta$\footnote{$\theta$ denotes the estimate of the ratio between positive and negative labels, which could be measured empirically through the training data. Empirically, we have $\hat{\theta}=0.1$} & 50.00 \\ \hline
      
      Model 1 & 2 hidden layers (HL) & 77.96 & 69.17 & 24.13 \\ \hline
      Model 2 & 3 HL & 95.52 & 80.66 & 84.19\\ \hline
      Model 3 & 4 HL & 95.15 & 79.26 & 82.67 \\ \hline
      Model 4 & 5 HL & 95.14 & 81.71 & 81.01 \\ \hline
      Model 5 & 5 HL + weighted loss ($\omega = 5$) &  95.67 & 79.05 & 86.39 \\ \hline
    \end{tabular}}
    \caption{Performances \label{table:1ddetectionrates}}
    }

\end{table}

Going forward, we select model 4 as it performs well and the resulting size of the weights matrices per layer is significantly smaller than model 3. It is debatable whether the differences between models 3, 4 and 5 are statistically significant. Furthermore, we note that the weighted loss function did not improve the chosen performance metrics.

\subsubsection{Linear Advection}
\label{sub:adv1d}
Consider a linear advection equation with $a \in \mathbb{R}$:

\begin{equation}
\label{eq:linadv}
\frac{\partial}{\partial t}u + a\frac{\partial}{\partial x}u = 0
\end{equation}
and periodic boundary conditions.

\textbf{Case of a Gaussian pulse:}
We consider the following initial condition:
\begin{equation}
\label{eq:gaussian}
u_0(x) = 1 + 3\exp(-100(x-0.5)^2) \quad (x,t)\in[0,1]\times \mathbb{R}^+
\end{equation}
with advection velocity $a = 1$.

\begin{figure}
\centering
\includegraphics[scale=0.45]{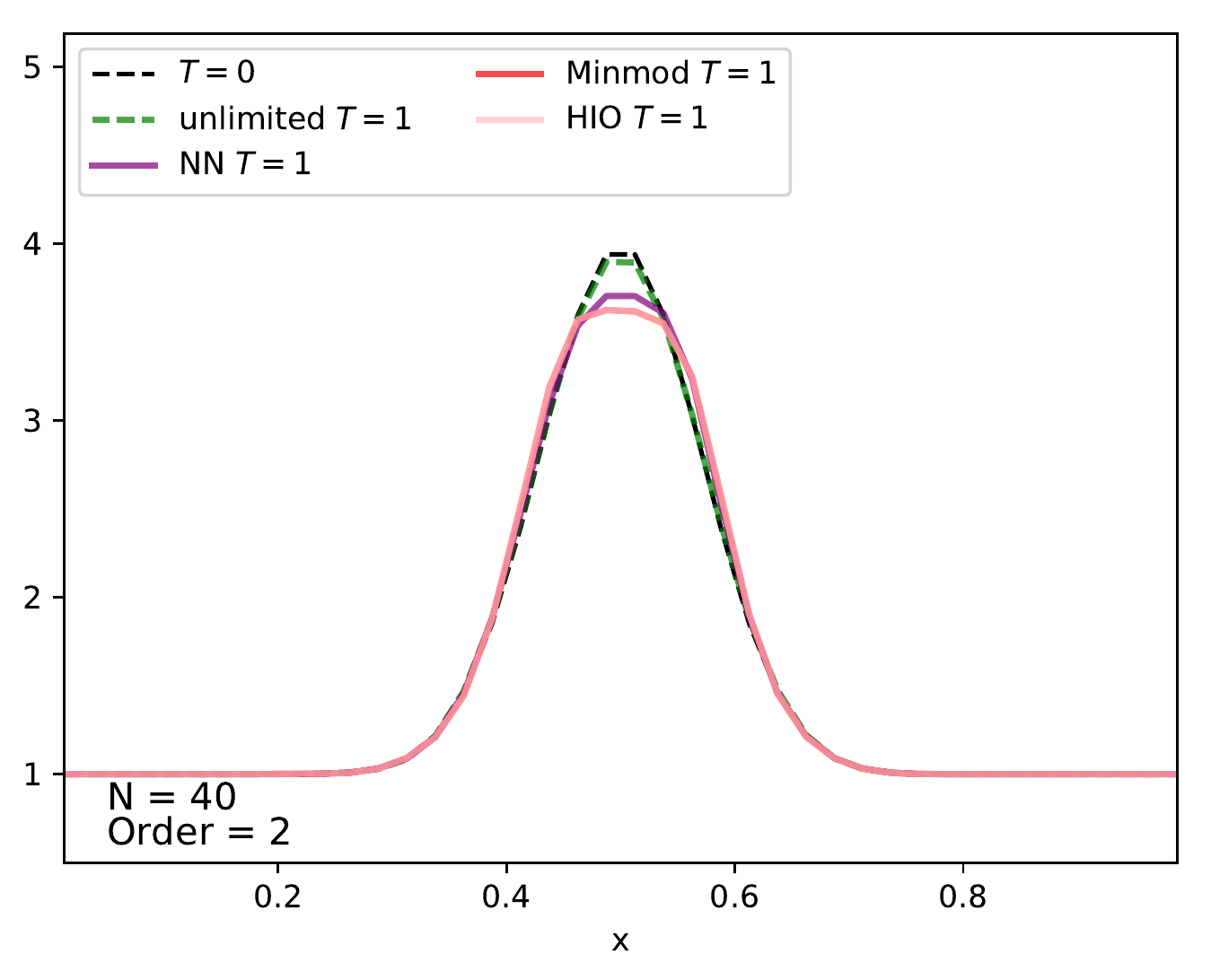}
\includegraphics[scale=0.45]{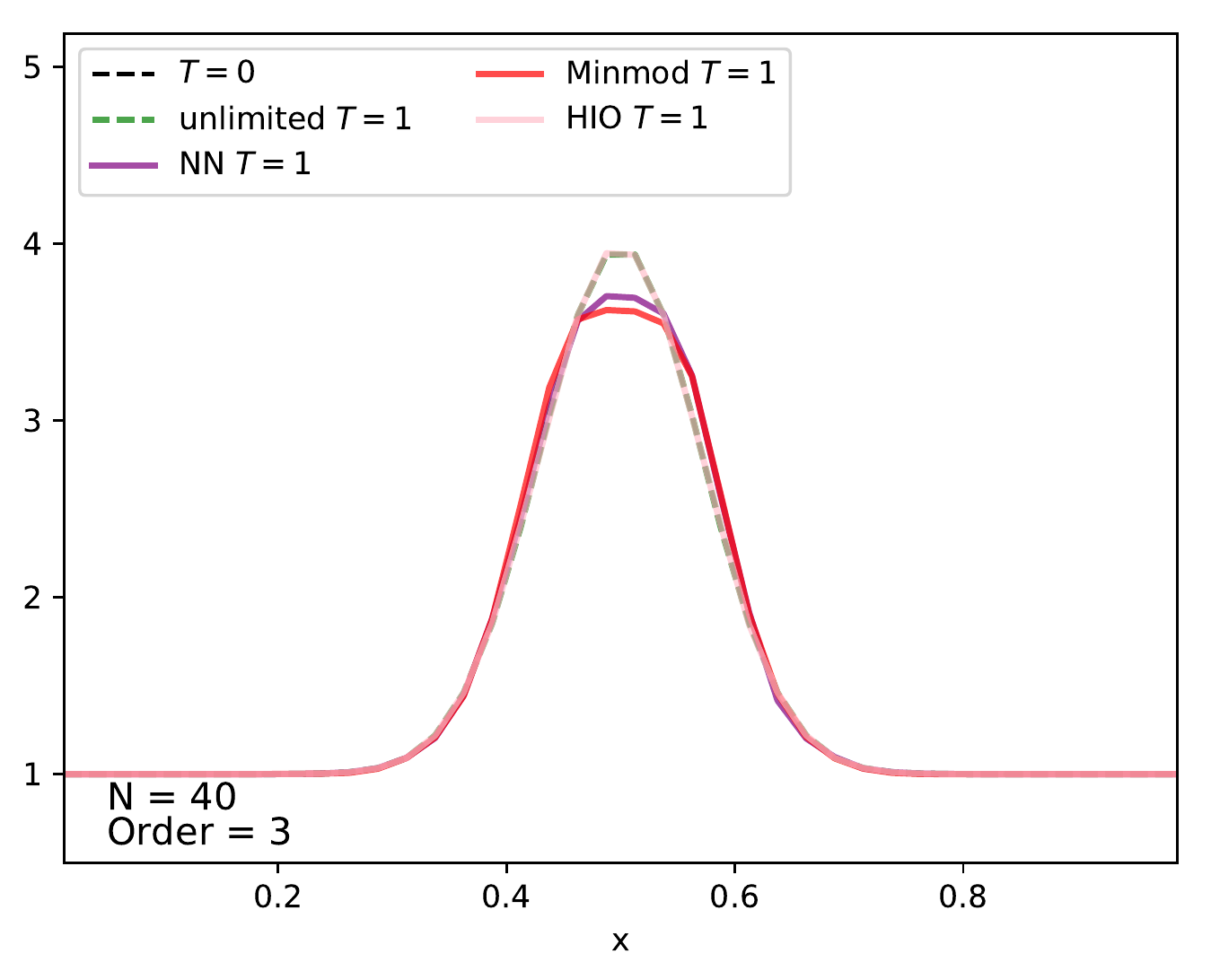}\\
\includegraphics[scale=0.45]{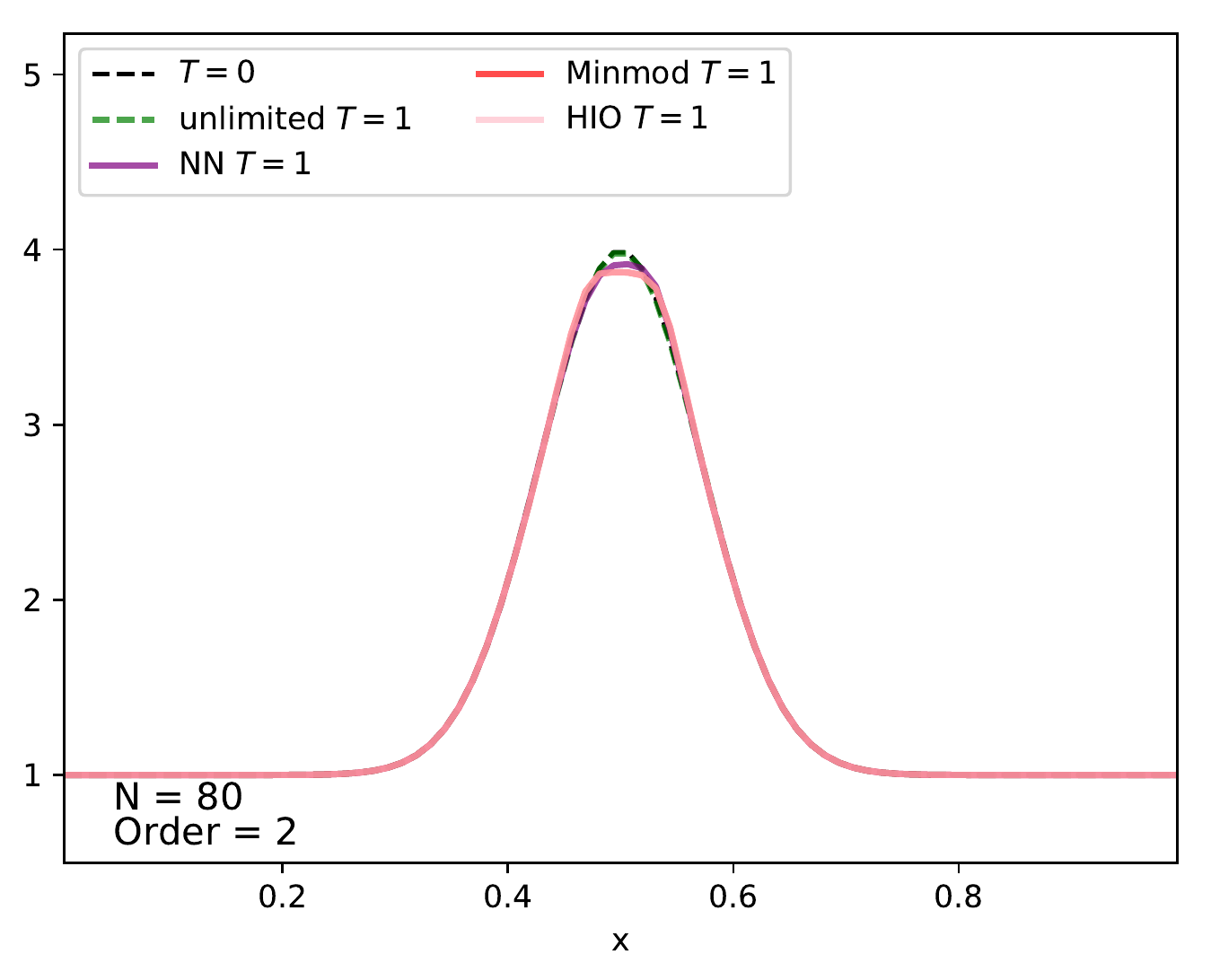}
\includegraphics[scale=0.45]{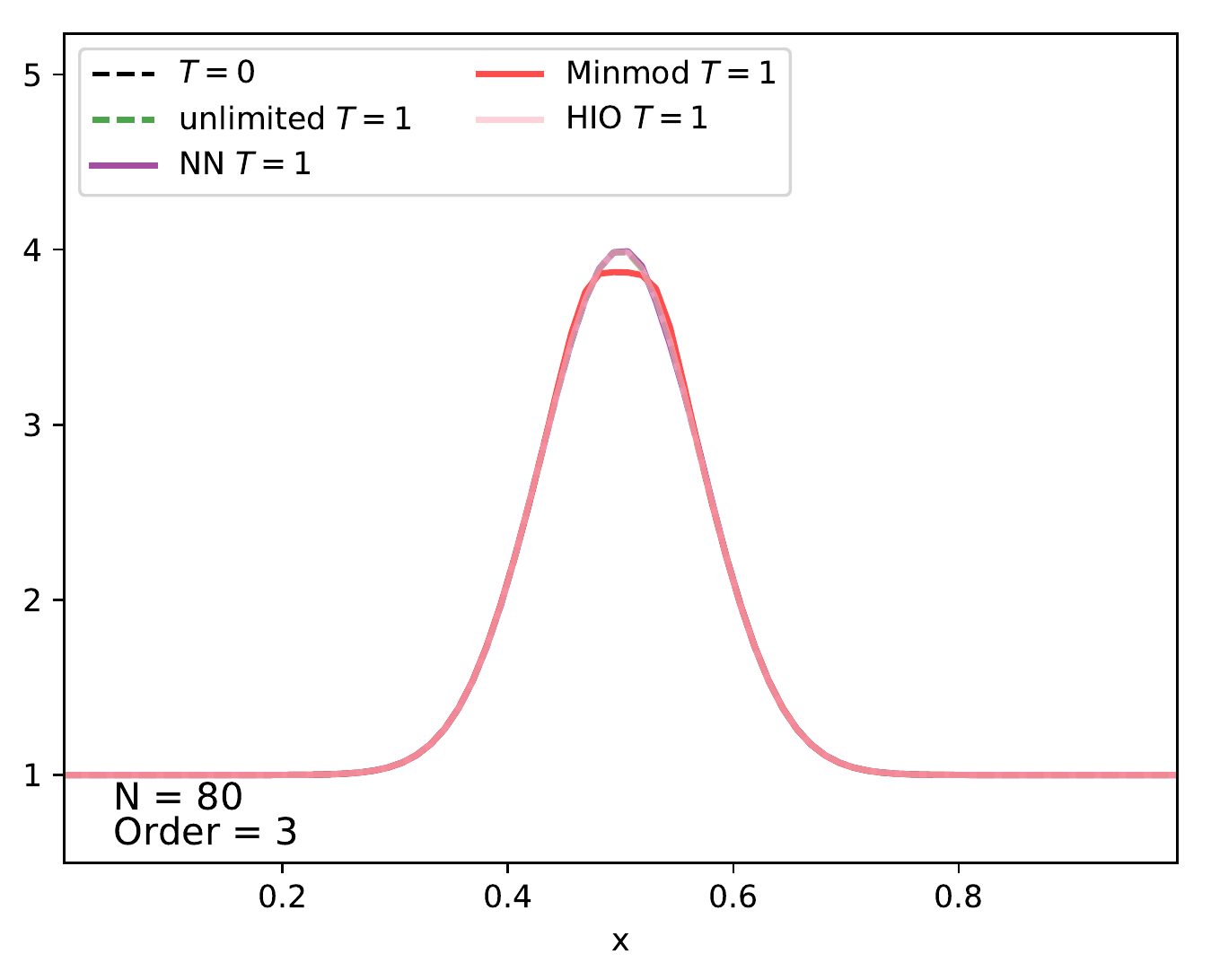}
\caption{Maxima clipping of the Gaussian pulse \ref{eq:gaussian} after one full crossing for approximation order of 2 (left) and 3 (right) for grid size $N = 40$ and $N = 80$.\label{fig:clippinggaussian}}
\end{figure}

The convergence is shown in table \ref{tab:1dgaussianconv} after one full crossing for orders $2$ and $3$. Furthermore, in figure \ref{fig:clippinggaussian} we show how the maxima is clipped using different methods for a grid size of $N=40$ and $N=80$. We note that the \textit{Minmod limiter} clips the maximum value of the solution. The \textit{hierarchical high order limiter} (denoted as HIO throughout this section) behaves as Minmod for the second order case, but for third order it does not limit the solution. The learned limiter (denoted as NN) behaves similarly for both orders, but it is to note that it seems to outperform the HIO limiter for the second order case. It makes sense that the performance does not depend (as much as the HIO limiter) on the order of the method, as we train the algorithm with only nodal information. 

\begin{table}
  \small
  \centering
  \caption{$L^1$ error for one crossing of the Gaussian pulse \eqref{eq:gaussian} using different limiters}\label{tab:1dgaussianconv}
  \subcaption{order 2}
    \begin{tabular}{ c | c| c| c | c }
      \hline \hline
      \textbf{N} & \textbf{No Limiter}& \textbf{MinMod} & \textbf{HIO} &	\textbf{NN} \\ \hline \hline 
20 & 3.04E-02 $\mid$ 0.0 &  7.72E-02 $\mid$ 0.0 &  7.72E-02 $\mid$ 0.0 &  7.66E-02 $\mid$ 0.0 \\ \hline 
40 & 3.51E-03 $\mid$ 0.6 &  1.82E-02 $\mid$ 0.4 &  1.82E-02 $\mid$ 0.4 &  1.34E-02 $\mid$ 0.4 \\ \hline 
60 & 1.15E-03 $\mid$ 3.1 &  7.00E-03 $\mid$ 2.1 &  7.00E-03 $\mid$ 2.1 &  5.06E-03 $\mid$ 2.5 \\ \hline 
80 & 5.60E-04 $\mid$ 3.0 &  3.54E-03 $\mid$ 2.2 &  3.54E-03 $\mid$ 2.2 &  2.27E-03 $\mid$ 2.5 \\ \hline 
100 & 3.28E-04 $\mid$ 2.9 &  2.06E-03 $\mid$ 2.2 &  2.06E-03 $\mid$ 2.2 &  3.28E-04 $\mid$ 2.5 \\ \hline 

\end{tabular}
  \\
  \subcaption{order 3}
\begin{tabular}{ c | c| c| c | c }
      \hline \hline
      \textbf{N} & \textbf{No Limiter}& \textbf{MinMod} & \textbf{HIO} &	\textbf{NN} \\ \hline \hline 
20 & 1.26E-03 $\mid$ 0.0 &  4.35E-02 $\mid$ 0.0 &  5.00E-03 $\mid$ 0.0 &  4.32E-02 $\mid$ 0.0 \\ \hline 
40 & 1.53E-04 $\mid$ 1.1 &  1.10E-02 $\mid$ 0.5 &  8.43E-04 $\mid$ 0.9 &  8.75E-03 $\mid$ 0.5 \\ \hline 
60 & 4.45E-05 $\mid$ 3.0 &  4.37E-03 $\mid$ 2.0 &  2.43E-04 $\mid$ 2.6 &  2.37E-03 $\mid$ 2.3 \\ \hline 
80 & 1.93E-05 $\mid$ 3.0 &  2.26E-03 $\mid$ 2.1 &  9.87E-05 $\mid$ 2.7 &  3.36E-04 $\mid$ 2.6 \\ \hline 
100 & 1.01E-05 $\mid$ 3.0 &  1.35E-03 $\mid$ 2.1 &  4.67E-05 $\mid$ 2.8 &  2.04E-04 $\mid$ 3.3 \\ \hline 

\end{tabular}
\end{table}

\textbf{Case of a smooth pulse and square hat:}
The following initial conditions contain a smooth Gaussian pulse and a hat function:

\begin{equation}
\label{eq:gausshat}
u_0(x) = \left\{
\begin{array}{ll}
      2, & |x-0.7| \leq 0.1 \\
	  1+\exp\left(-\frac{(x-0.25)^2}{2\times0.05^2}\right), & \mbox{otherwise} \\
\end{array} 
\right.  \quad (x,t)\in[0,1]\times \mathbb{R}^+,
\end{equation}
again with advection velocity $a = 1$.

\begin{figure}
\centering
\includegraphics[width=0.45\textwidth]{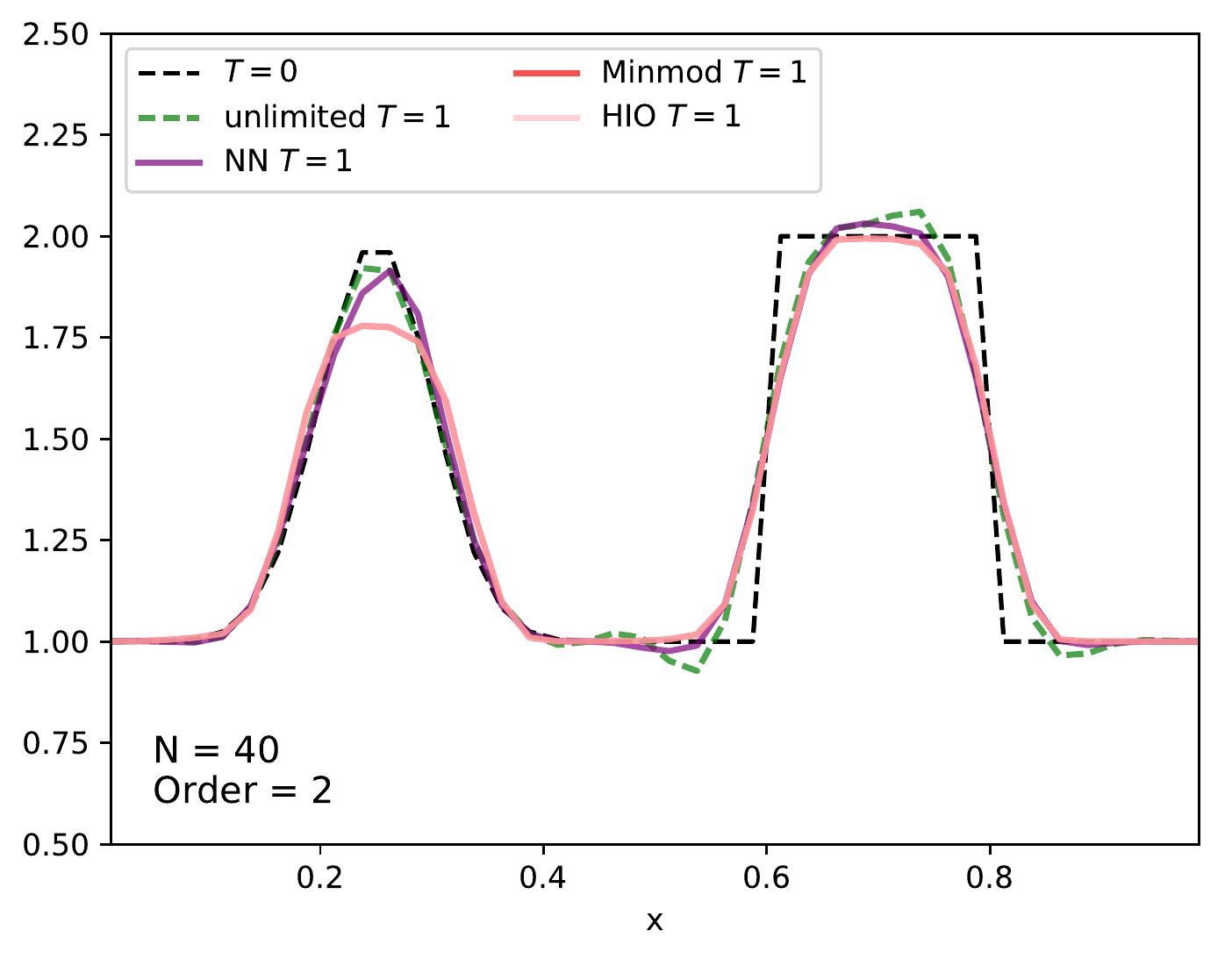}
\includegraphics[width=0.45\textwidth]{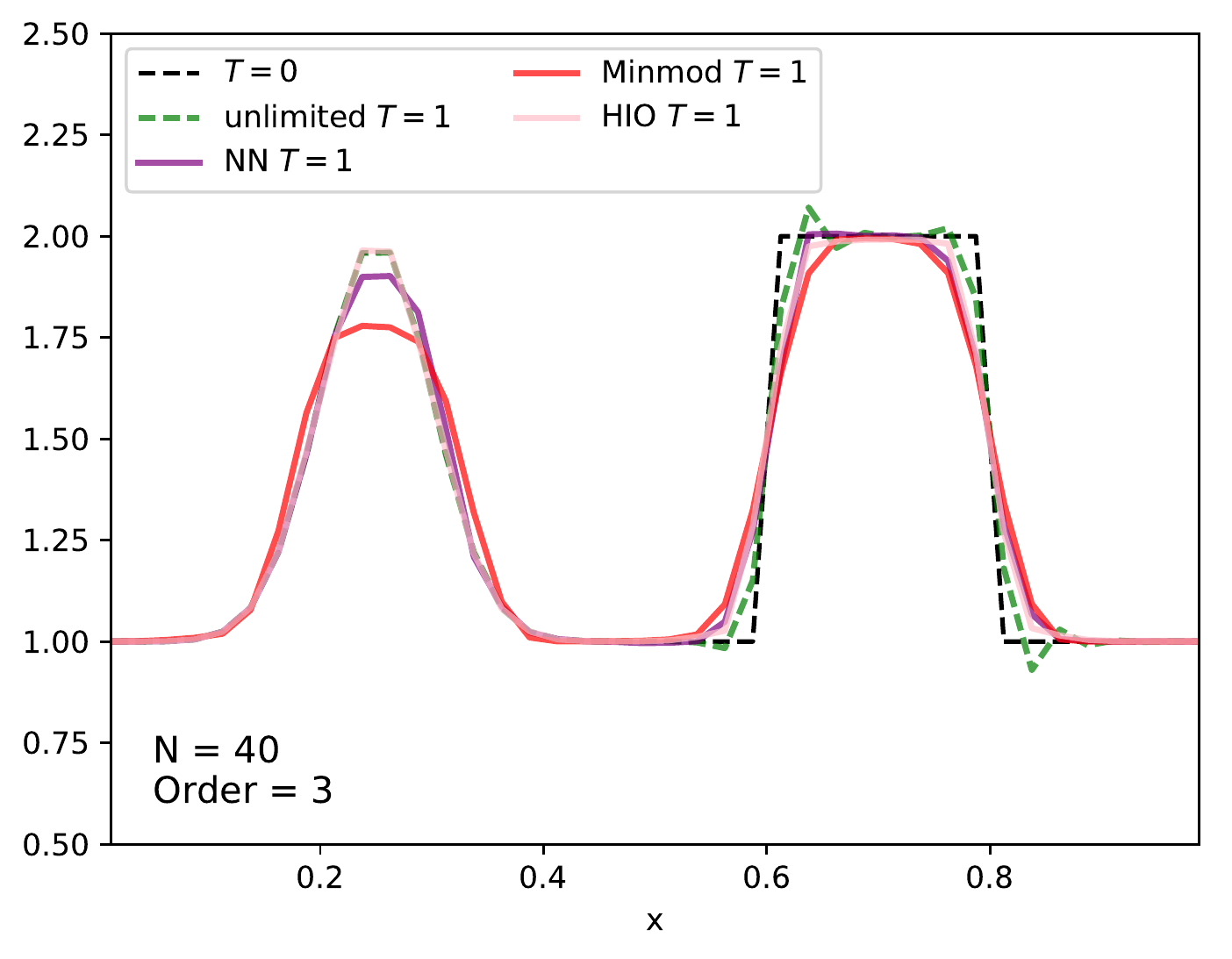}\\
\includegraphics[width=0.45\textwidth]{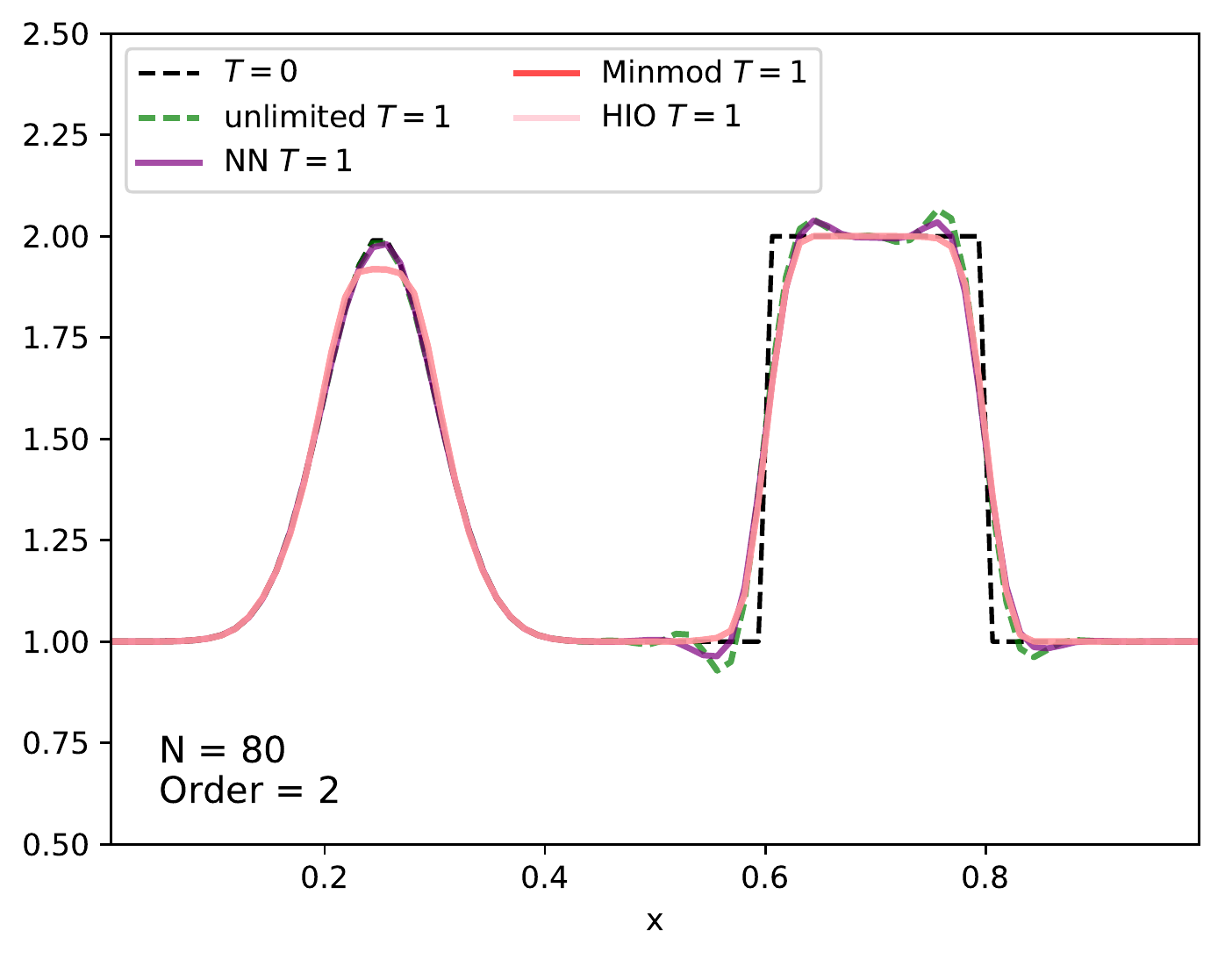}
\includegraphics[width=0.45\textwidth]{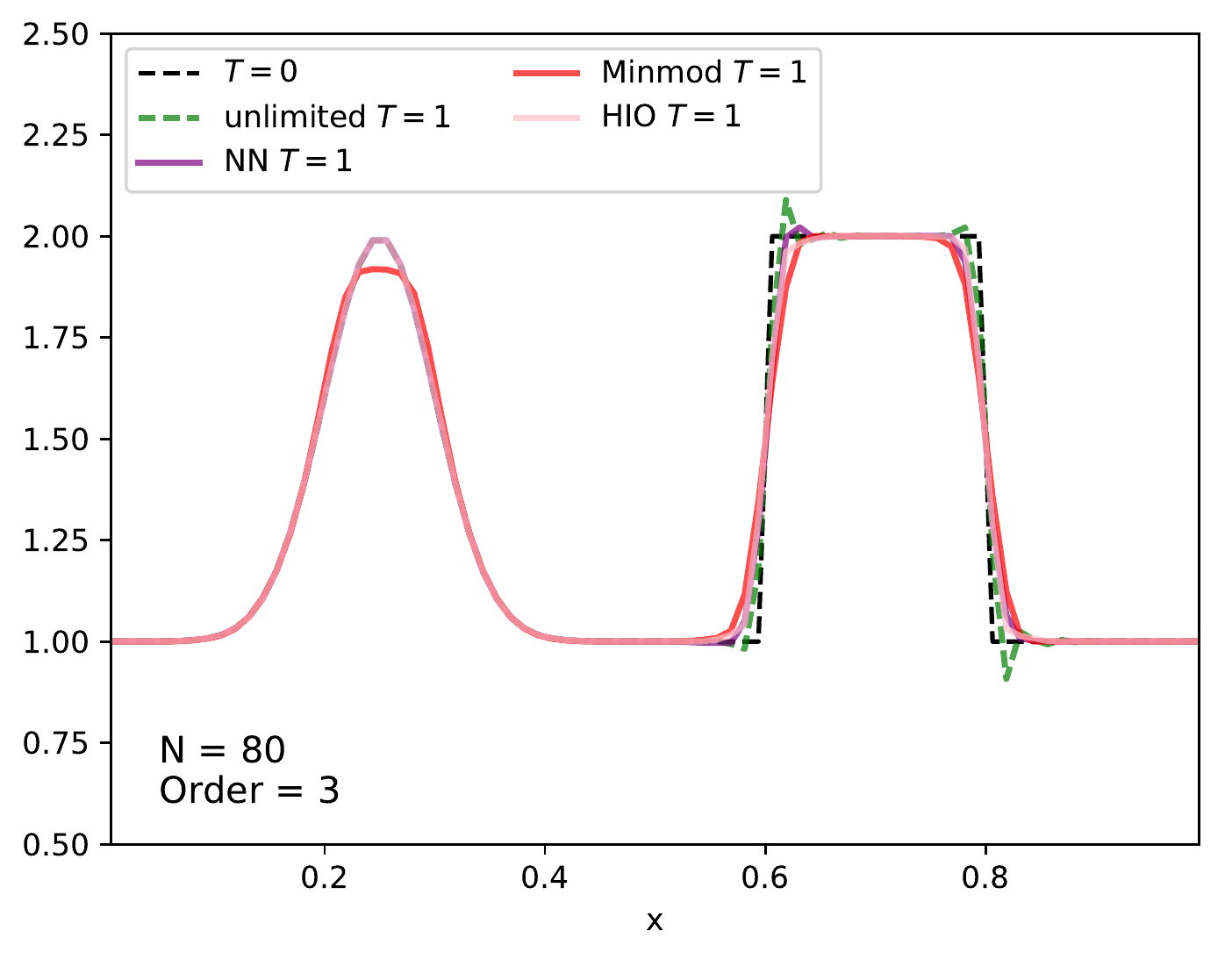}\\
\caption{Maxima clipping of the Gaussian and hat pulses \ref{eq:gausshat} after one full crossing for approximation order of 2 (left) and 3 (right) for grid size $N = 40$ and $N = 80$.\label{fig:zoomhat}}
\end{figure}

The convergence is shown in table \ref{tab:1dgaussianhatconv} after one full crossing for orders $2$ and $3$. Furthermore, in figure \ref{fig:zoomhat} we show how the different limiters perform, for a grid size of $N=40$ and $N=80$.

\begin{table}
  \small
  \centering
  \caption{$L^1$ error for one crossing of the Gaussian pulse and hat function \eqref{eq:gausshat} using different limiters}\label{tab:1dgaussianhatconv}
  \subcaption{order 2}
    \begin{tabular}{ c | c| c| c | c }
      \hline \hline
      \textbf{N} & \textbf{No Limiter}& \textbf{MinMod} & \textbf{HIO} &	\textbf{NN} \\ \hline \hline 
20 & 6.38E-02 $\mid$ 0.0 &  7.52E-02 $\mid$ 0.0 &  7.52E-02 $\mid$ 0.0 &  6.73E-02 $\mid$ 0.0 \\ \hline 
40 & 2.73E-02 $\mid$ 0.5 &  3.07E-02 $\mid$ 0.4 &  3.07E-02 $\mid$ 0.4 &  2.72E-02 $\mid$ 0.5 \\ \hline 
60 & 1.87E-02 $\mid$ 1.2 &  1.85E-02 $\mid$ 1.3 &  1.85E-02 $\mid$ 1.3 &  1.72E-02 $\mid$ 1.3 \\ \hline 
80 & 1.40E-02 $\mid$ 1.1 &  1.31E-02 $\mid$ 1.3 &  1.31E-02 $\mid$ 1.3 &  1.31E-02 $\mid$ 1.2 \\ \hline 
100 & 1.13E-02 $\mid$ 1.1 &  1.02E-02 $\mid$ 1.3 &  1.02E-02 $\mid$ 1.3 &  1.07E-02 $\mid$ 1.2 \\ \hline 
\end{tabular}
  \\
  \subcaption{order 3}
\begin{tabular}{ c | c| c| c | c }
      \hline \hline
      \textbf{N} & \textbf{No Limiter}& \textbf{MinMod} & \textbf{HIO} &	\textbf{NN} \\ \hline \hline 
20 & 2.08E-02 $\mid$ 0.0 &  4.40E-02 $\mid$ 0.0 &  2.57E-02 $\mid$ 0.0 &  3.30E-02 $\mid$ 0.0 \\ \hline 
40 & 1.08E-02 $\mid$ 0.6 &  1.67E-02 $\mid$ 0.5 &  1.04E-02 $\mid$ 0.6 &  1.28E-02 $\mid$ 0.6 \\ \hline 
60 & 7.40E-03 $\mid$ 0.9 &  1.02E-02 $\mid$ 1.4 &  6.79E-03 $\mid$ 1.3 &  7.43E-03 $\mid$ 1.4 \\ \hline 
80 & 5.66E-03 $\mid$ 0.9 &  7.27E-03 $\mid$ 1.3 &  5.12E-03 $\mid$ 1.2 &  5.36E-03 $\mid$ 1.4 \\ \hline 
100 & 4.62E-03 $\mid$ 0.9 &  5.69E-03 $\mid$ 1.3 &  4.11E-03 $\mid$ 1.2 &  4.32E-03 $\mid$ 1.3 \\ \hline 

\end{tabular}
\end{table}

\subsubsection{Euler equation}
\label{sub:euler1d}
Now we consider the 1-dimensional Euler equations, which describe the behavior of an inviscid flow. This system of equations describe the evolution of a density $\rho$, a velocity $v$, a pressure $p$ and total energy $E$.

\begin{align}
\label{eq:euler}
\frac{\partial}{\partial t}\rho + \frac{\partial}{\partial x}(\rho v) &= 0 \\
\frac{\partial}{\partial t}(\rho v) + \frac{\partial}{\partial x}(\rho v^2 + p) &= 0 \\
\frac{\partial}{\partial t}E + \frac{\partial}{\partial x}((E+p)v) &= 0.
\end{align}

The system is closed with equation of state for an ideal gas:

\[ \rho e = p/(\gamma - 1),\]
where $e = E-\frac{1}{2}\rho v^2$ is the internal energy.

The NN limiter is applied sequentially for each variable. Empirically, better results were observed when applying the limiter to primitive variables.

\textbf{Case of Sod shock tube:}
We consider the standard Sod shock tube test, given by the initial conditions: 
\begin{equation}
\label{eq:sod}
(\rho, v, p)(x,0) = \left\{
\begin{array}{ll}
      (1.0,0.0,1.0) & 0.0 < x \leq 0.5 \\
      (0.125, 0.0, 0.1) & 0.5 < x < 1.0  \\
\end{array} 
\right.  \quad (x,t)\in[0,1]\times [0.0,0.1],
\end{equation}

and $\gamma = 1.4$ and gradient free boundary conditions.

\begin{figure}
\centering
\includegraphics[width=0.5\textwidth]{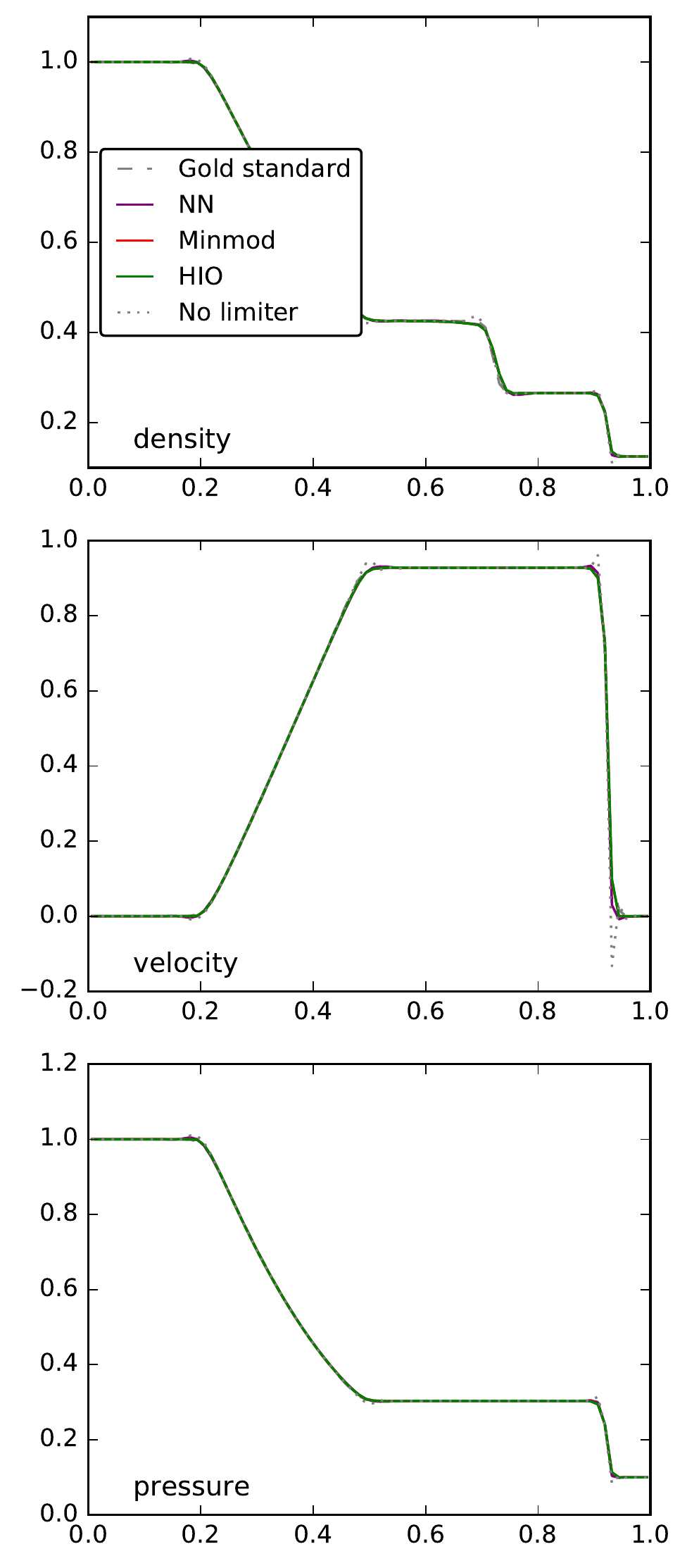}
\caption{Density, velocity and pressure fields for the Sod shock tube \eqref{eq:sod} at $T=0.24$ for grid size $N = 80$ and approximation order 2.
\label{fig:sod2}}
\end{figure}

\begin{figure}
\centering
\includegraphics[width=0.5\textwidth]{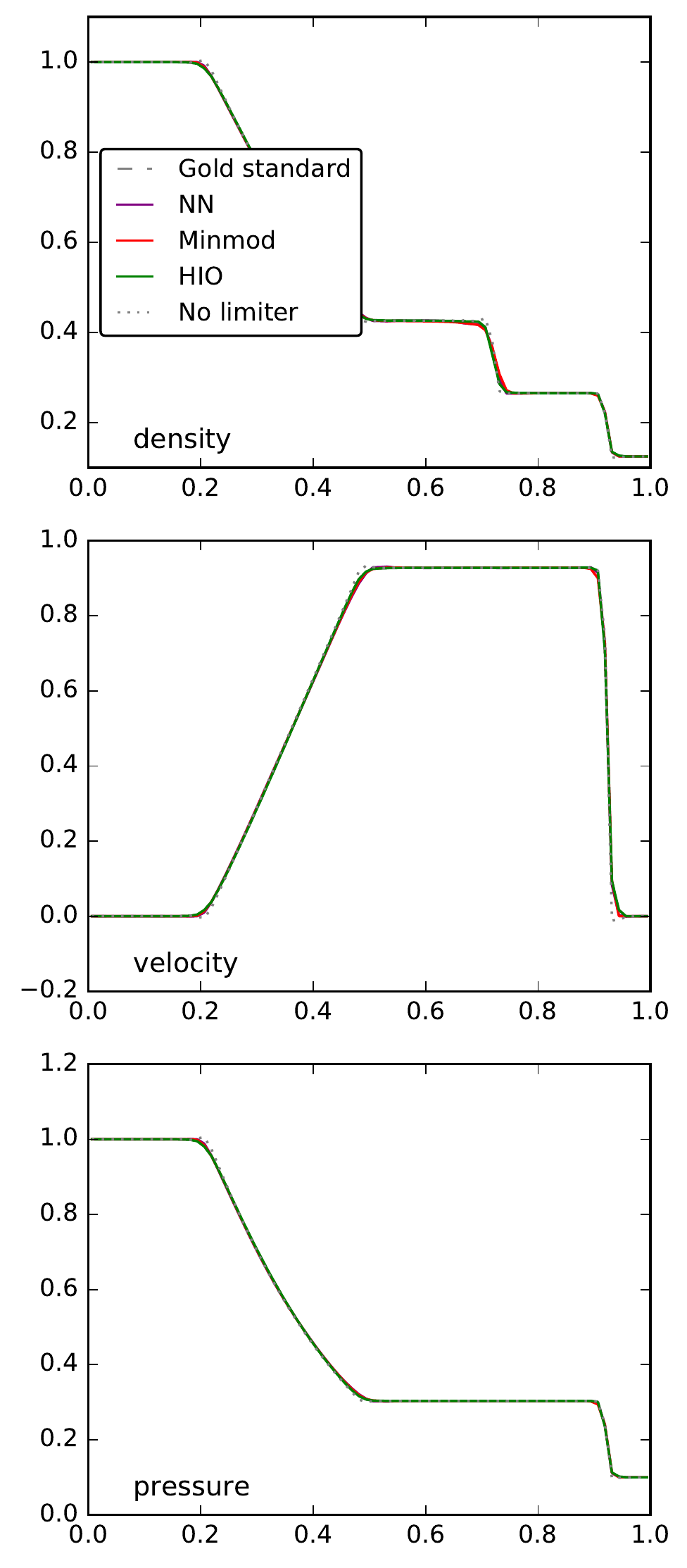}
\caption{Density, velocity and pressure fields for the Sod shock tube \eqref{eq:sod} at $T=0.24$ for grid size $N = 80$ and approximation order 3.\label{fig:sod3}}
\end{figure}

In figures \ref{fig:sod2} and \ref{fig:sod3} we show the comparison between different limiters at $T = 0.24$, for schemes of order 2 and 3 respectively. For the second order results, we note that for the density field (first panel in figure \ref{fig:sod2}), the solution produced by the NN limiter seems oscillation free and similar to the solution of the HIO limiter, but for the velocity and pressure fields (second and third panel in figure \ref{fig:sod2}) some oscillations were not corrected enough, and clearly the behaviour of the HIO limiter is more desirable. This could be due to the fact that the shock detection is done over the primitive variables. This could be improved by performing the shock detection over the characteristic variables. For third order (figure \ref{fig:sod3}) we see that the performances are comparable between the NN limiter and the HIO limiter.

\textbf{Case of Blast wave:}
Next we consider the interacting blast waves test, given by the initial conditions:

\begin{equation}
\label{eq:blast}
(\rho, v, p)(x,0) = \left\{
\begin{array}{ll}
      (1.0,0.0,1000.0) & 0.0 < x \leq 0.1 \\
      (1.0, 0.0, 0.01) & 0.1 < x \leq 0.9  \\
      (1.0, 0.0, 100.0) & 0.9 < x < 1.0  \\
\end{array} 
\right.  \quad (x,t)\in[0,1]\times [0.0,0.038]
\end{equation}
with $\gamma = 1.4$ and reflexive boundary conditions.

\begin{figure}
\centering
\includegraphics[width=0.5\textwidth]{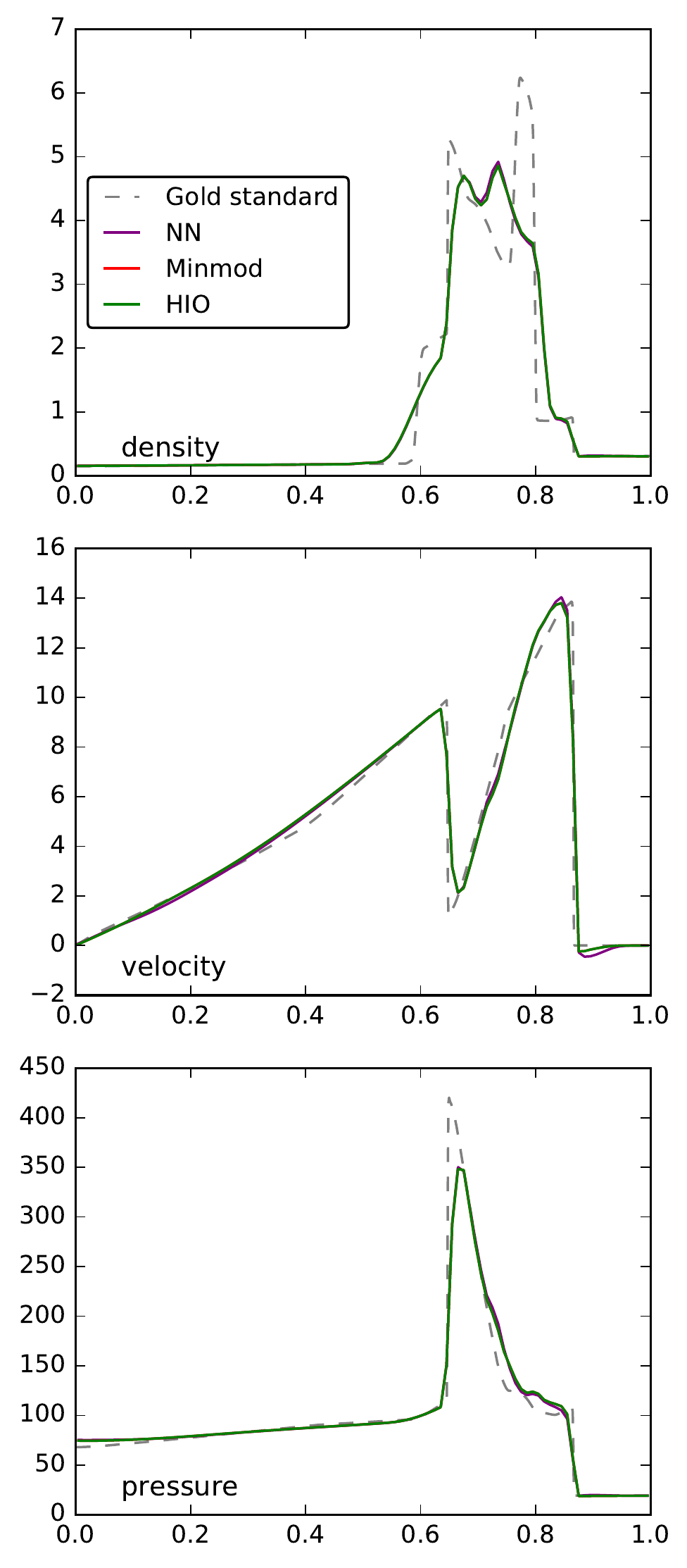}
\caption{Density, velocity and pressure fields of the blast wave interaction \ref{eq:blast} at $T=0.038$ for grid size $N = 100$ and approximation order 2.\label{fig:o2blast}}
\end{figure}

\begin{figure}
\centering{
\includegraphics[width=0.5\textwidth]{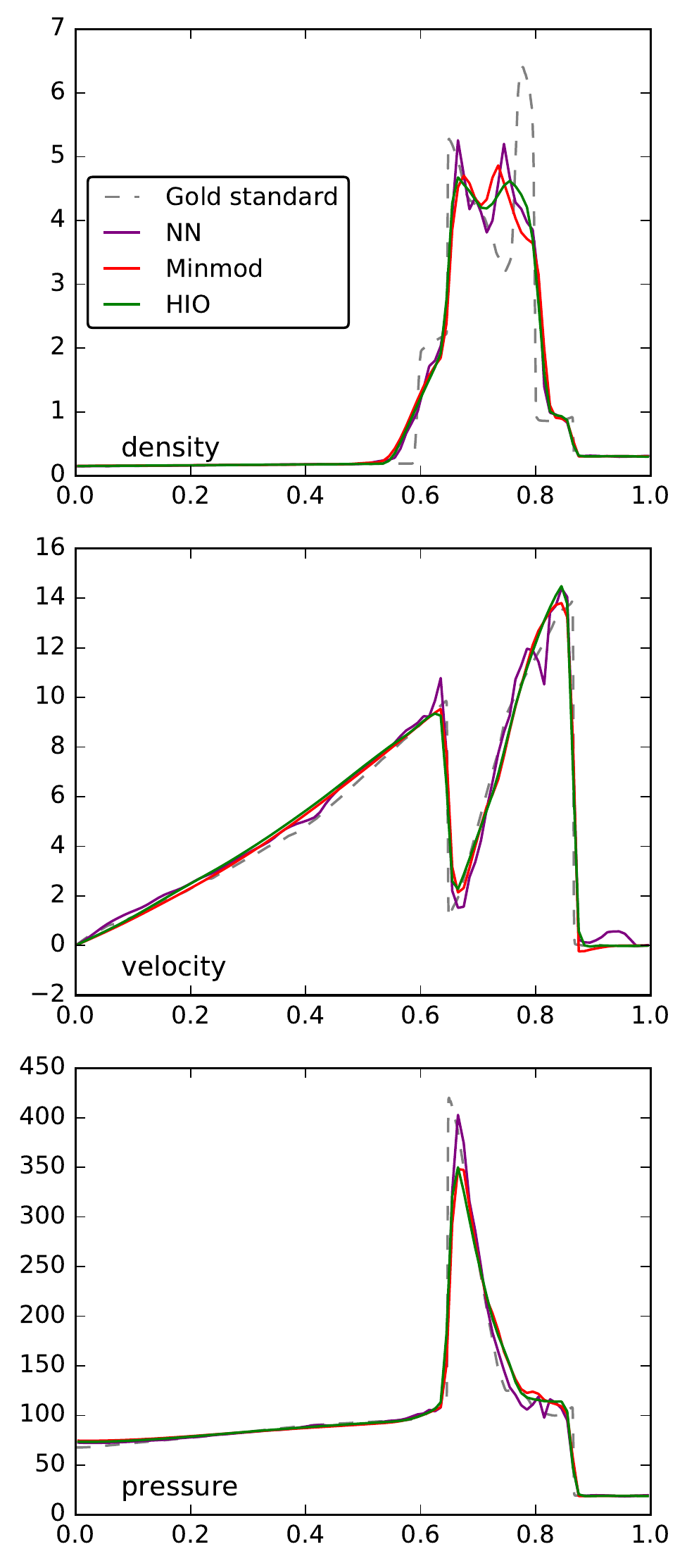}
\caption{Density, velocity and pressure fields of the Blast Wave interaction at $T=0.038$ for grid size $N = 100$ and approximation order 3.\label{fig:o3blast}}}
\end{figure}

In figures \ref{fig:o2blast} and \ref{fig:o3blast} we show the comparison between different limiters at $T = 0.038$ for different orders. The unlimited solution is not shown because the code crashes due to the pressure becoming negative shortly after the start of the simulation for orders higher than 1. The dashed line denotes a high resolution solution, run with $N = 1000$, third order with the HIO limiter. We can note that the NN limiter is not as good at suppressing oscillations as Minmod and the HIO limiter, but stabilises the solution enough to finish the run. Furthermore, we note that the peak is better preserved, which means that it looks like the limiting is less strong than Minmod and the HIO limiter. 

\textbf{Case of Shu-Osher problem:}
Next we consider the Shu-Osher test, given by the initial conditions:

\begin{equation}
\label{eq:shu-osher}
(\rho, v, p)(x,0) = \left\{
\begin{array}{ll}
      (3.857143,2.629369,10.3333) & 0.0 < x \leq \frac{1}{8} \\
      (1+0.2\sin(8\pi x), 0.0, 1.0) & \frac{1}{8} < x \leq 1.0  \\
\end{array} 
\right.  \quad (x,t)\in[0,1]\times [0.0,0.178]
\end{equation}
with $\gamma = 1.4$ and at the boundary, the variables are set to their initial values.

\begin{figure}
\centering{
\includegraphics[width=0.5\textwidth]{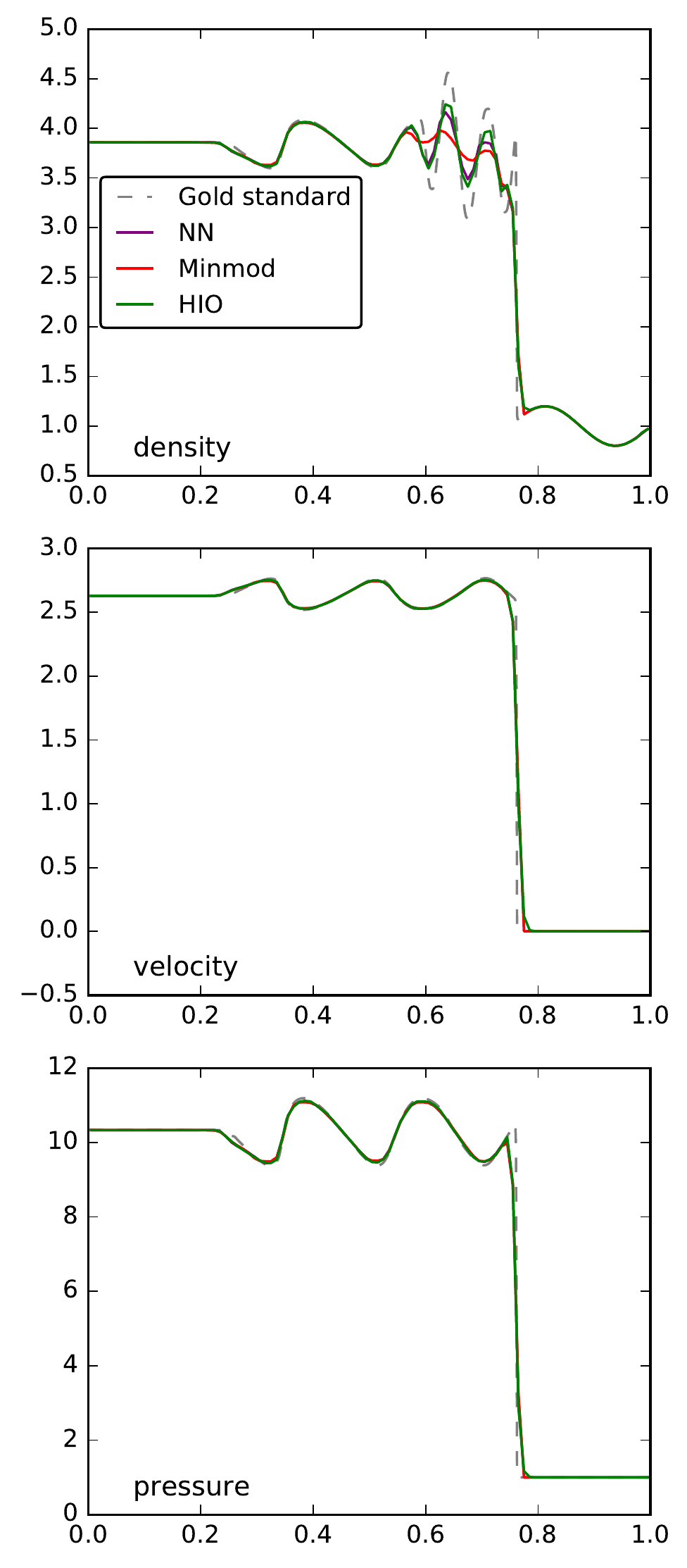}
\caption{Density, velocity and pressure fields of the Shu-Osher problem at $T=0.019$ for grid size $N = 100$ and approximation order 3.\label{fig:o3shu}}}
\end{figure}

In figure \ref{fig:o3shu} we show the comparison between different limiters at $T = 0.178$ for $3^{rd}$ order. The dashed line denotes a high resolution solution, run with $N = 1000$, third order with the HIO limiter. We can note that the NN limiter does not clip the solution as much as the Minmod, but it performs slighly worse than the HIO limiter. 

\subsubsection{Transfer learning}
\label{subsec:transfer1d}
In this section, we show the performance of the neural network limiter applied to residual distribution scheme (RDS), and we compare with a state of the art limiting technique, MOOD\cite{mood}. Other stabilisation strategies for RDS require parameter tuning which are problem dependent.

\begin{figure}
\centering{
\includegraphics[width=0.5\textwidth]{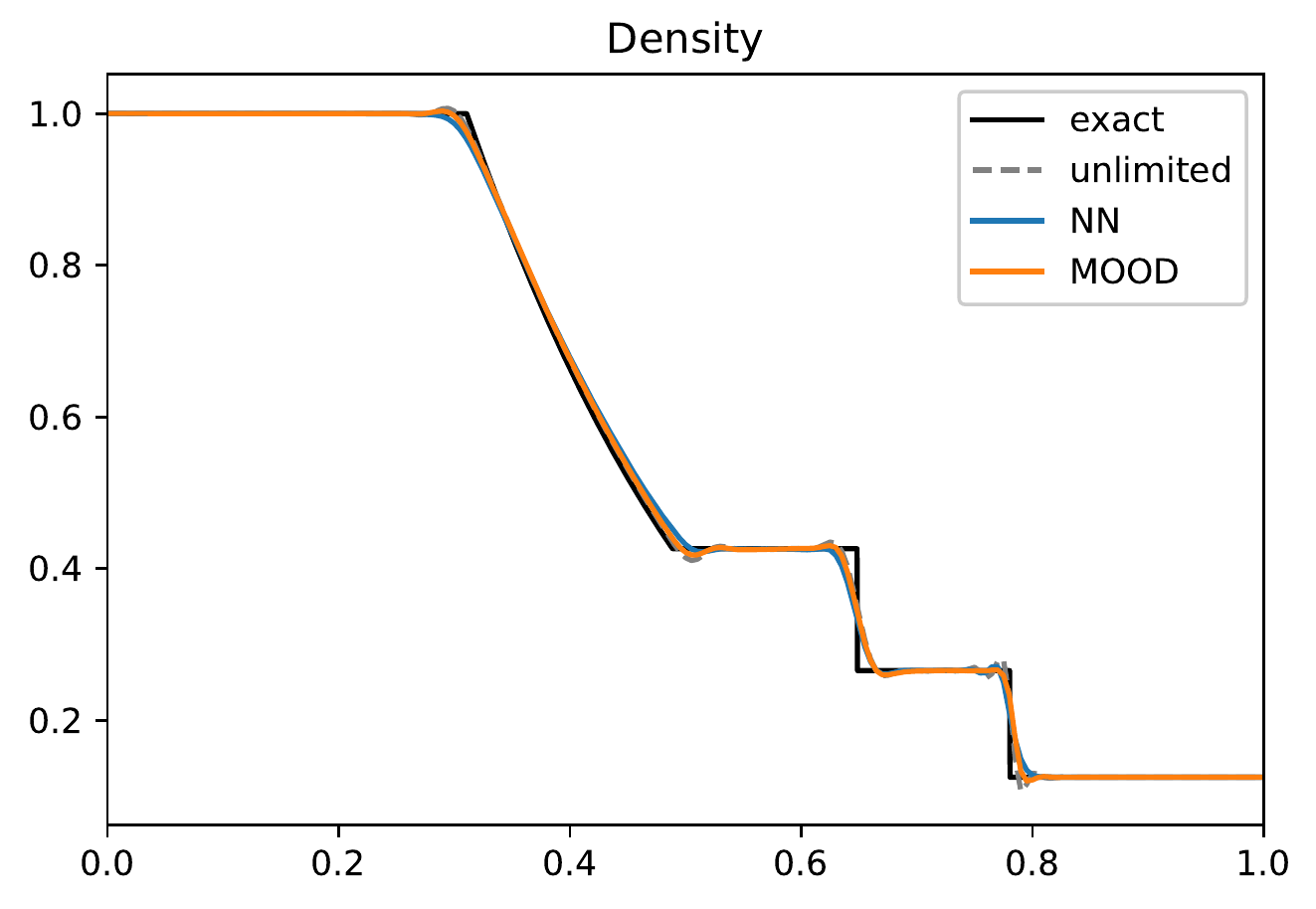}\\
\includegraphics[width=0.5\textwidth]{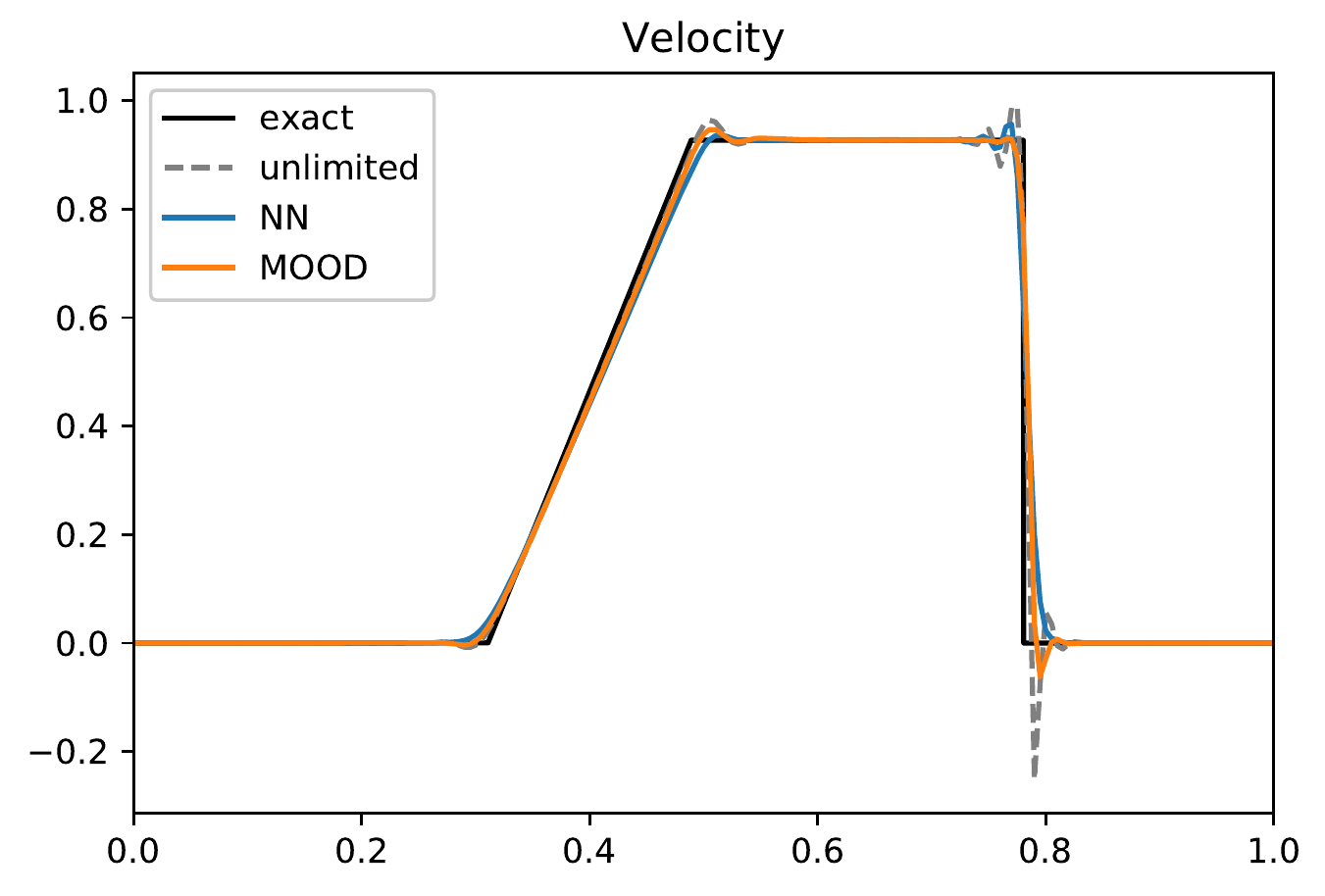}\\
\includegraphics[width=0.5\textwidth]{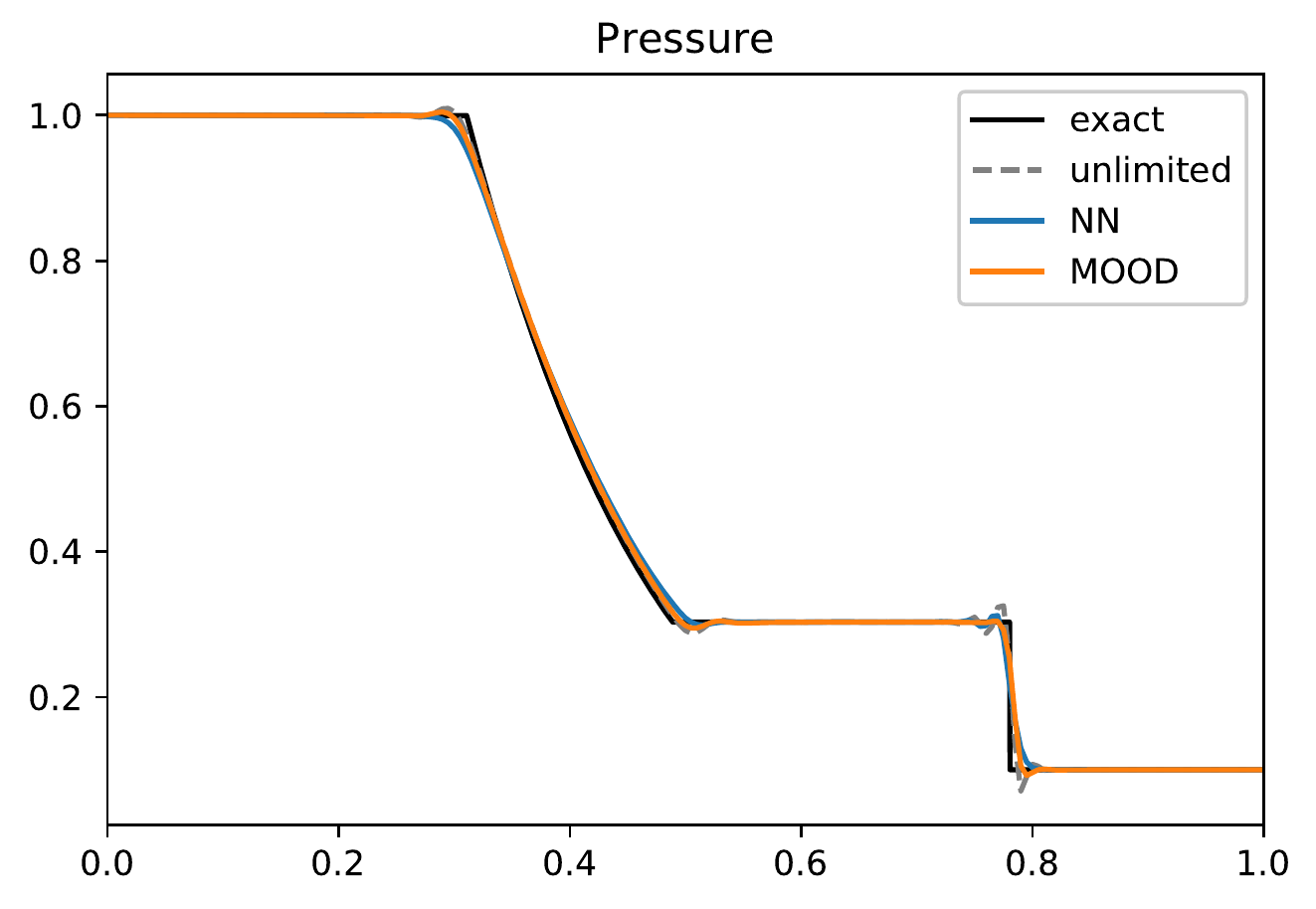}
\caption{Density, velocity and pressure fields of the sod shock.\label{fig:sod1drd}}}
\end{figure}

\textbf{Case of Sod shock:} The initial conditions are given as in \eqref{eq:sod}. A qualitative result is shown in figure \ref{fig:sod1drd}. It can be noted that the NN limiter is slightly more diffusive than MOOD, but that it seems to detect solution undershoots more accurately.

\begin{figure}
\centering{
\includegraphics[width=0.5\textwidth]{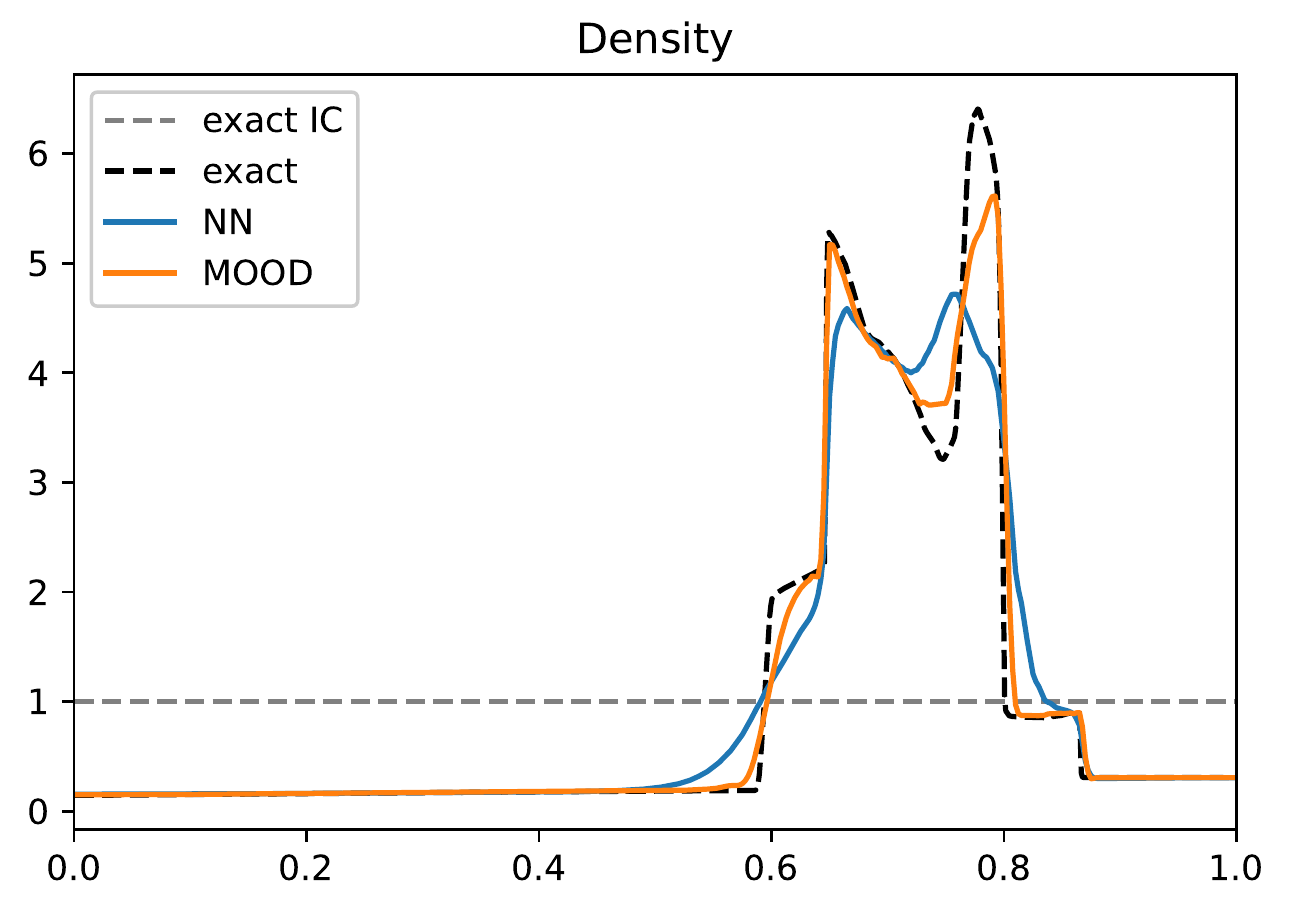}\\
\includegraphics[width=0.5\textwidth]{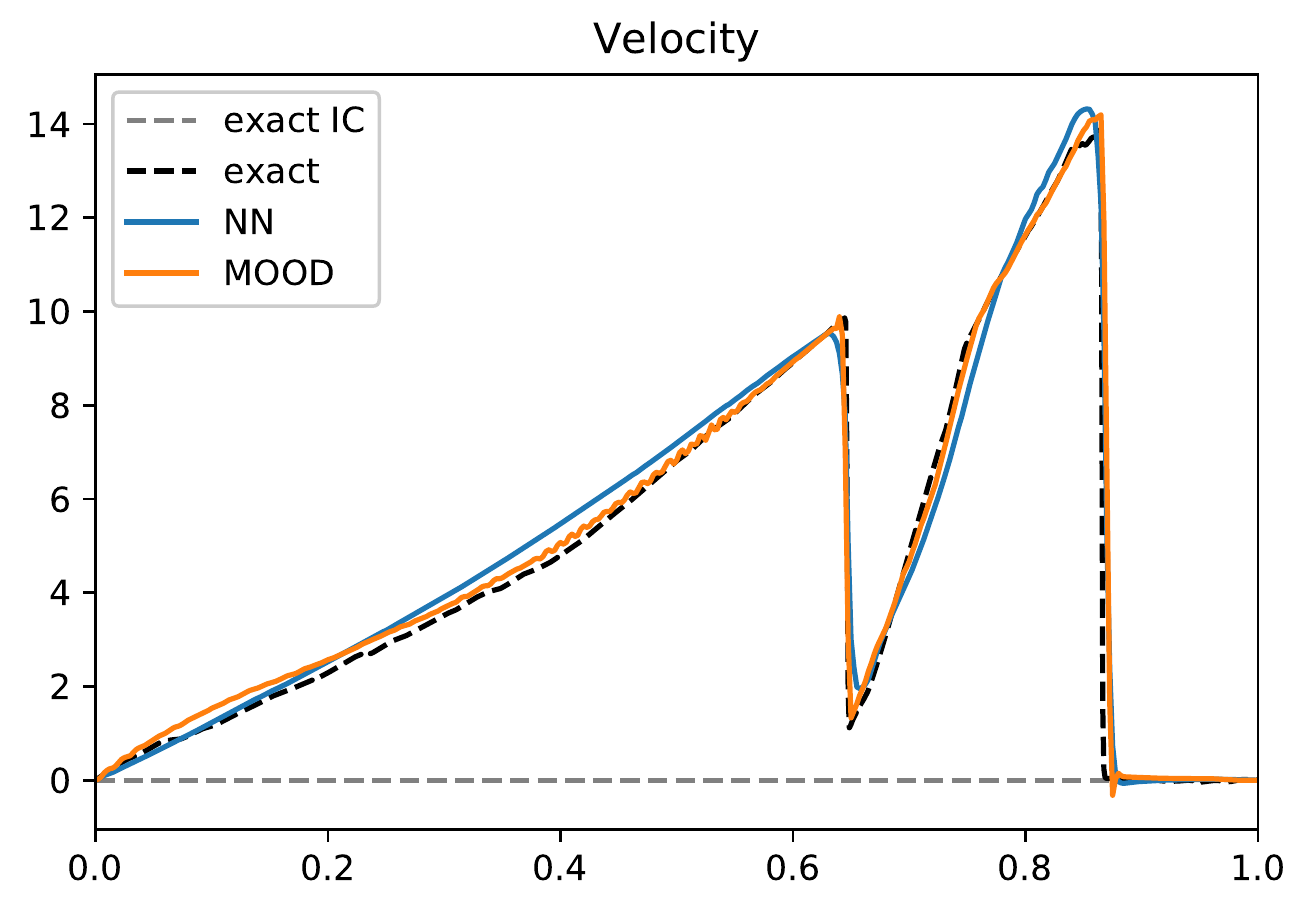}\\
\includegraphics[width=0.5\textwidth]{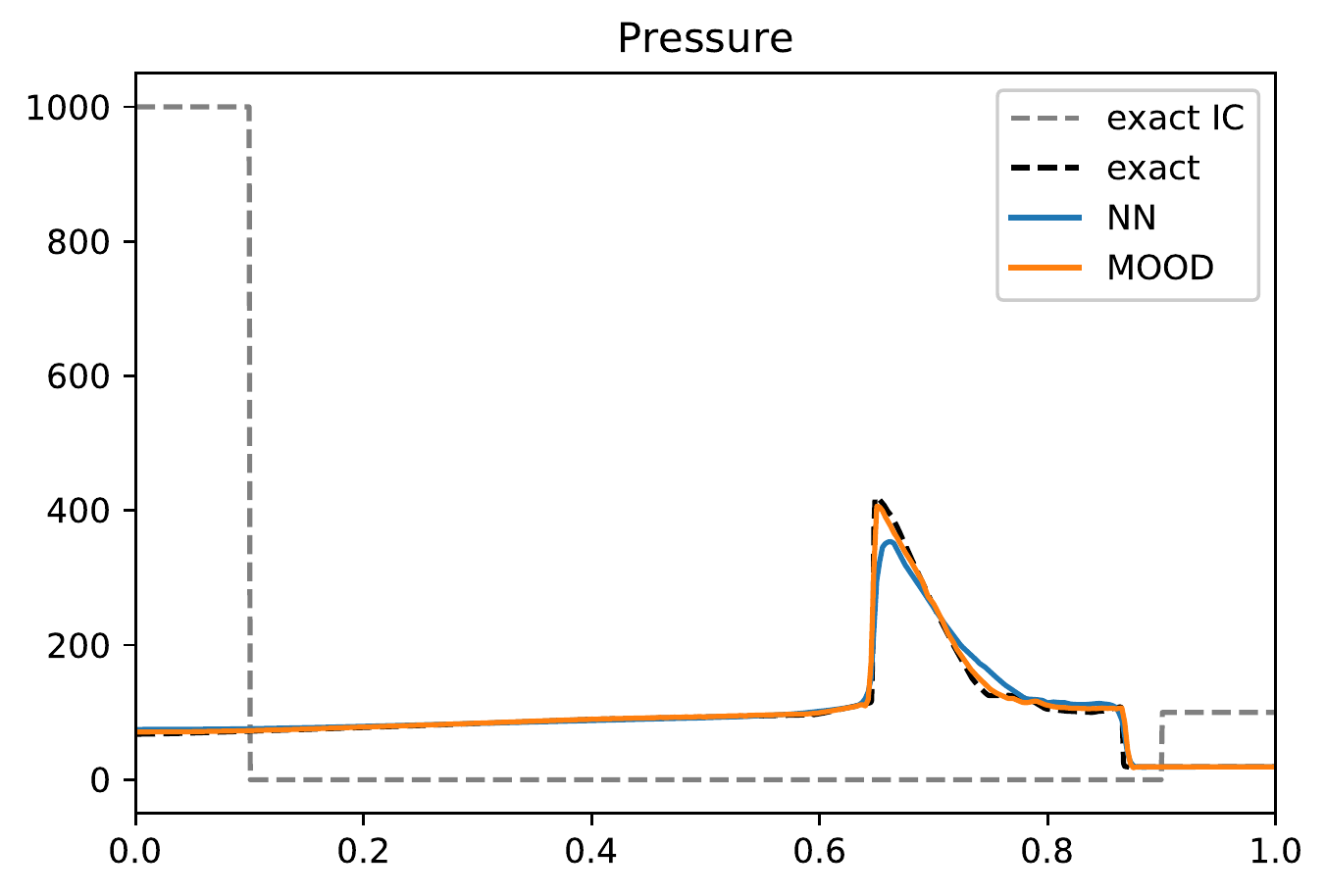}
\caption{Density, velocity and pressure fields of the blast-wave test case.\label{fig:blastrd}}}
\end{figure}

\textbf{Case of Blast-wave:} The initial conditions are given as in \eqref{eq:blast}. A qualitative result is shown in figure \ref{fig:blastrd}. In this example, one can see that MOOD is significantly less diffusive than the NN limiter.

\section{Numerical experiments: 2-dimensional problems}
\label{sec:results2d}

Similarly to section \ref{sec:1ddetection}, we first train a set of neural networks varying the number of neurons and layers, and pick the one with best recall/precision score. The analysis can be found in the repository \cite{repository}. In this section, we compare the performance of the black box limiter (denoted as NN) with the Minmod limiter and hierarchical high order limiter (denoted as HIO) through the $L^1$ error norm. We perform some tests for the linear advection equation and Euler system of equations. The initial conditions are chosen as different from the ones used for the training. Finally, in section \ref{sec:transfer2d}, we show the results for the transfer to residual distribution, for structured and unstructured meshes.

\subsection{Linear Advection}
\label{sub:adv}
Consider a linear advection equation with $\vec{a} \in \mathbb{R}^2$:

\begin{equation}
\label{eq:linadv}
\frac{\partial}{\partial t}u + \vec{a}\cdot \nabla u  = 0
\end{equation}
and periodic boundary conditions.

\textbf{Smooth initial condition:} We consider the following initial conditions, which contain a smooth function:

\begin{equation}
\label{eq:smoothring}
u_0(x) = 1 + \sin(2\pi r)^{10}  \quad (x,t)\in[0,1]\times \mathbb{R}^+,
\end{equation}
with advection velocity $\vec{a} = (1,1)$ and periodic boundary conditions.

\begin{table}[htbp]
  \small
  \centering
  \caption{$L^1$ error for one crossing of the smooth ring \eqref{eq:smoothring} using different limiters}\label{tab:2dsmoothring}
  \subcaption{order 2}
  %\subfloat[order 2]{
    \begin{tabular}{ c | c| c| c | c }
      \hline \hline
      \textbf{N} & \textbf{No Limiter}& \textbf{MinMod} & \textbf{HIO} &	\textbf{NN} \\ \hline \hline 
$16^2$ & 6.74E-02 $\mid$ 0.0 & 1.26E-01  $\mid$ 0.0 &  1.23E-01 $\mid$ 0.0 & 6.79E-02 $\mid$ 0.0 \\ \hline 
$32^2$ & 2.00E-02 $\mid$ 0.5 &  4.81E-02 $\mid$ 0.4 &   4.14E-02 $\mid$ 0.4 & 1.97E-02 $\mid$ 0.5 \\ \hline 
$64^2$ & 3.45E-03 $\mid$ 1.7 &  1.31E-02 $\mid$ 1.4 &  9.11E-03 $\mid$ 1.6 & 3.45E-03 $\mid$ 1.8 \\ \hline
$128^2$ & 5.42E-04 $\mid$ 2.1 &  2.75E-03  $\mid$ 1.6 &  1.59E-03 $\mid$ 1.9 & 5.41E-04 $\mid$ 2.1 \\ \hline
\end{tabular}
\bigskip
  %}\\
  %\subfloat[order 3]{
  \subcaption{order 3}
\begin{tabular}{ c | c| c| c | c }
    \hline \hline
    \textbf{N} & \textbf{No Limiter}& \textbf{MinMod} & \textbf{HIO} &	\textbf{NN} \\ \hline \hline 
        $16^2$ & 6.92E-03 $\mid$ 0.0 & 1.26E-01  $\mid$ 0.0 &  1.17E-01 $\mid$ 0.0 & 1.88E-02 $\mid$ 0.0 \\ \hline 
        $32^2$ & 4.58E-04 $\mid$ 0.9 &   4.80E-02 $\mid$ 0.4 &   2.92E-02 $\mid$ 0.4 & 7.91E-04 $\mid$ 0.7  \\ \hline 
        $64^2$ & 7.93E-05 $\mid$ 3.9 &   1.31E-02  $\mid$ 1.4 &  6.85E-03 $\mid$ 2.0 & 1.08E-04 $\mid$ 4.6 \\ \hline
        $128^2$ & 4.36E-05 $\mid$ 3.2 & 2.75E-03  $\mid$ 1.6 &  1.20E-03 $\mid$ 2.0 & 5.15E-05 $\mid$ 3.7 \\ \hline
\end{tabular}
%  }
\end{table}

\textbf{Case of smooth pulse and square hat:} We consider the following initial conditions, which contain a smooth Gaussian pulse and a hat function, defined in $(\vec{x},t)\in[0,1]^2\times \mathbb{R}^+$,

\begin{equation}
\label{eq:gausshat2d}
u_0(x) = 
\begin{cases}
      2, & (|\vec{x}-0.25|,|y-0.5|) \leq (0.1,0.1) \\
      1 + \exp\left(-100(||\vec{x}-\vec{x}_1||^2)\right),& x \geq 0.5 \\
      1, & \mbox{otherwise}
\end{cases} 
,
\end{equation}
again with advection velocity $\vec{a} = (1,1)$, $\vec{x}_1=(0.75,0.5)$.

\begin{table}[htbp]
  \small
  \centering
  \caption{$L^1$ error for one crossing of the Gaussian pulse and hat function \eqref{eq:gausshat2d} using different limiters}\label{tab:2dgausshatconv}
%  \subfloat[order 2]{
\subcaption{order 2}
    \begin{tabular}{ c | c| c| c | c }
      \hline \hline
      \textbf{N} & \textbf{No Limiter}& \textbf{MinMod} & \textbf{HIO} &	\textbf{NN} \\ \hline \hline 
$16^2$ & 3.04E-02 $\mid$ 0.0 & 5.37E-02 $\mid$ 0.0 &  4.94E-02 $\mid$ 0.0 & 3.20E-02 $\mid$ 0.0 \\ \hline 
$32^2$ & 2.13E-02 $\mid$ 0.6 &  2.63E-02 $\mid$ 0.5 &  2.63E-02 $\mid$ 0.5 & 2.12E-02 $\mid$ 0.6  \\ \hline 
$64^2$ & 1.12E-02 $\mid$ 0.5 &  1.16E-02 $\mid$ 1.0 &  1.01E-02 $\mid$ 0.9 & 1.11E-03 $\mid$ 0.6 \\ \hline
$128^2$ & 8.33E-03 $\mid$ 0.7 &  7.74E-03 $\mid$ 1.1 &  7.74E-03 $\mid$ 1.1 & 8.24E-03 $\mid$ 0.8 \\ \hline
\end{tabular}
\bigskip
  %}\\
  %\subfloat[order 3]{
  \subcaption{order 3}
\begin{tabular}{ c | c| c| c | c }
      \hline \hline
      \textbf{N} & \textbf{No Limiter}& \textbf{MinMod} & \textbf{HIO} &	\textbf{NN} \\ \hline \hline 
$16^2$ & 1.17E-02 $\mid$ 0.0 &  5.37E-02 $\mid$ 0.0 &  5.27E-02 $\mid$ 0.0 & 2.09E-02 $\mid$ 0.0 \\ \hline 
$32^2$ & 7.61E-03 $\mid$ 0.8 &  2.63E-02 $\mid$ 0.5 &   1.63E-02 $\mid$ 0.5 & 1.05E-02 $\mid$ 0.7  \\ \hline
$64^2$ & 4.63E-03 $\mid$ 0.6 &  1.16E-02 $\mid$ 1.0 &  7.58E-03 $\mid$ 1.7 & 5.73E-03 $\mid$ 1.0 \\ \hline
$128^2$ & 2.91E-03 $\mid$ 0.7 &  7.74E-03 $\mid$ 1.1  &  4.21E-03 $\mid$ 1.4 & 3.37E-03 $\mid$ 0.9 \\ \hline
\end{tabular}
  %}
\end{table}

\subsection{Euler equation}
\label{sub:euler}
Now we consider the 2-dimensional Euler equations, which describe the behaviour of an inviscid flow. This system of equations describe the evolution of a density $\rho$, a velocity vector $\vec{v}=(v_1,v_2)$, a pressure $p$ and total energy $E$.

\begin{align}
%\label{eq:nn-euler}
\partial_t \rho + \nabla \cdot (\rho \vec{v}) &= 0\\
\partial_t \rho \vec{v} + \nabla \cdot \left(\rho \vec{v}\otimes\vec{v} + p\mathcal{I}_{3}\right) &= 0\\
\partial_t E + \nabla \cdot \left( E + p\right)\vec{v} &= 0.
\end{align}

The system is closed equation of state for an ideal gas:

\[ \rho e = p/(\gamma - 1),\]
where $e = E-\frac{1}{2}\rho \abs{\vec{v}}^2$ is the internal energy.

\textbf{Case of the 2-dimensional Sod shock tube:} We consider the standard Sod shock tube test, given by the initial conditions: 

\begin{equation}
\label{eq:sod}
(\rho, v_x, v_y, p)(x,0) = \left\{
\begin{array}{ll}
      (1.0,0.0,0.0,1.0) & 0.0 < r \leq 0.5 \\
      (0.125, 0.0, 0.0, 0.1) & 0.5 < r < 1.0  \\
\end{array} 
\right.
\end{equation}
where $r=\sqrt{x^2+y^2}$, $(x,y)\in[0,1]^2$, $\gamma = 1.4$ and gradient free boundary conditions.

\begin{figure}
\centering{
\includegraphics[width=0.23\textwidth]{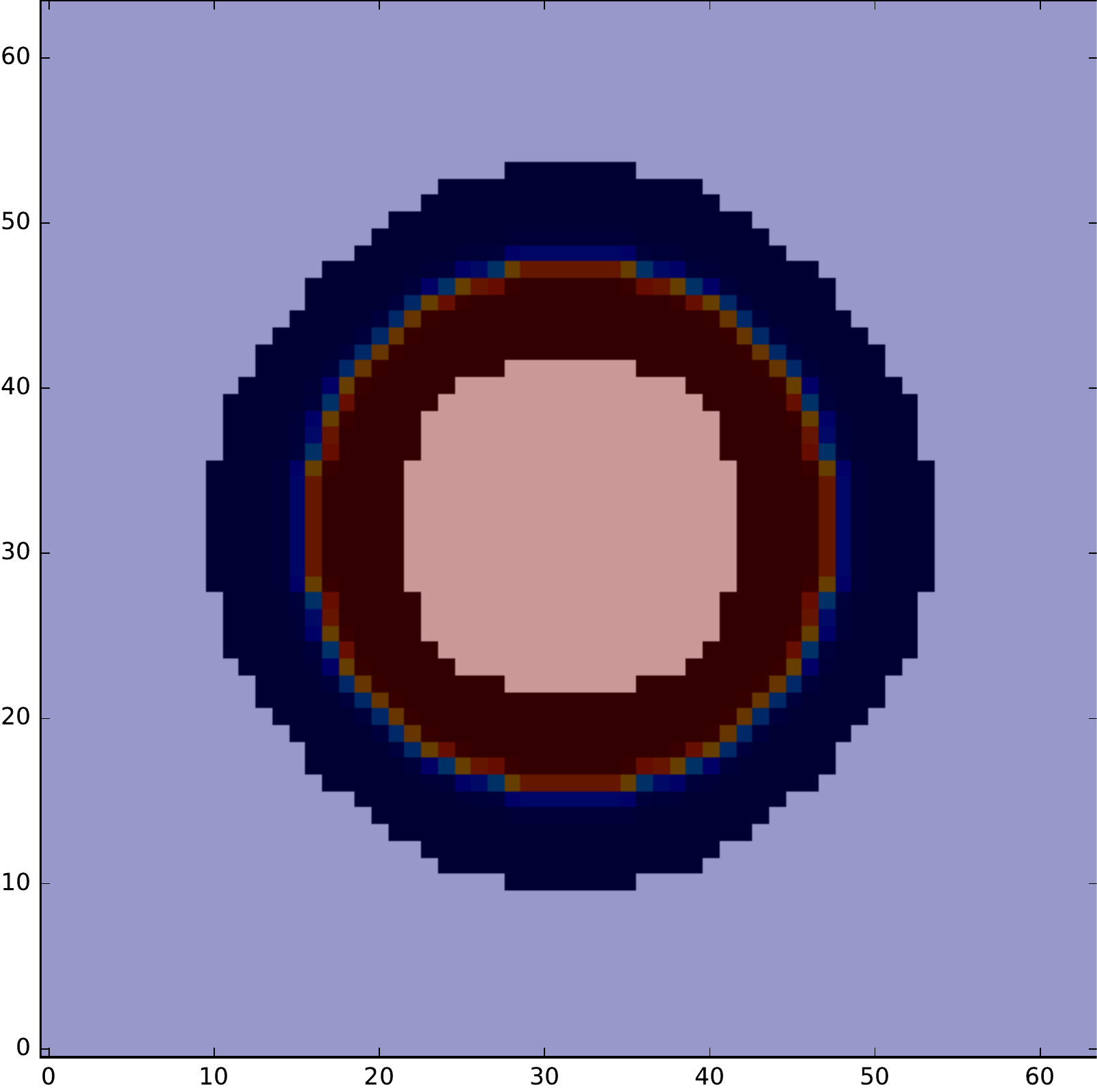}
\includegraphics[width=0.23\textwidth]{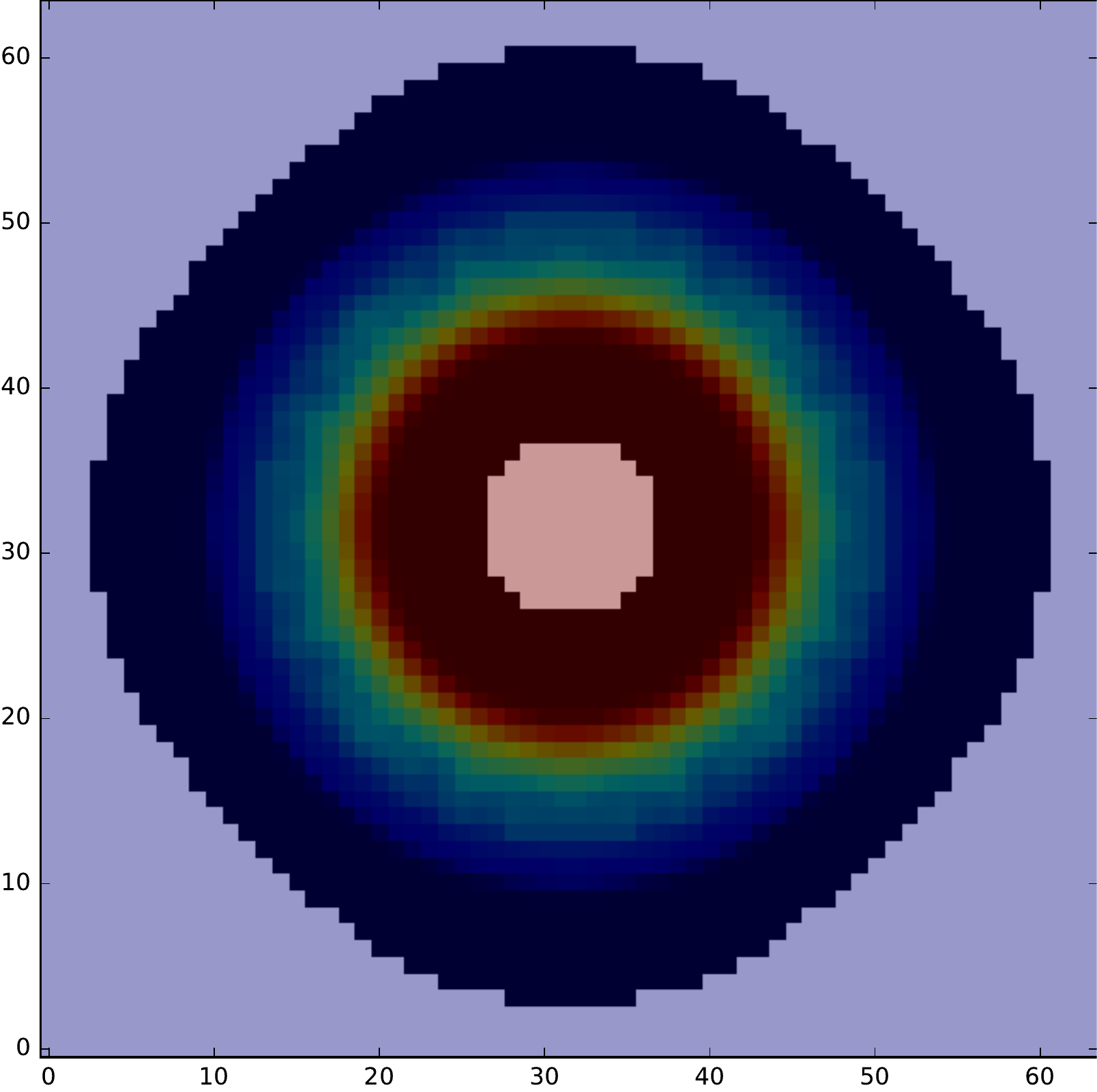}
\includegraphics[width=0.23\textwidth]{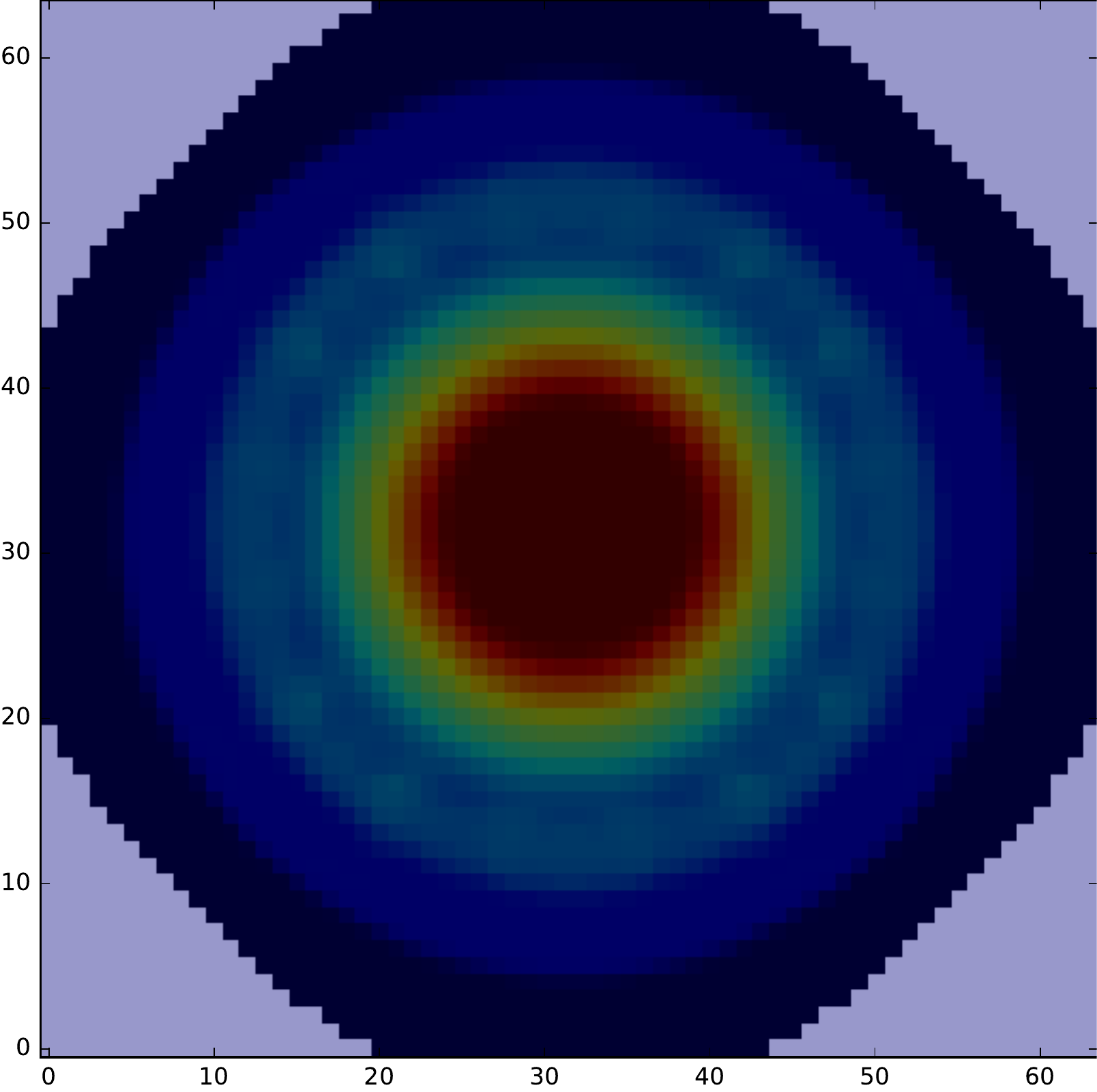}
\includegraphics[width=0.23\textwidth]{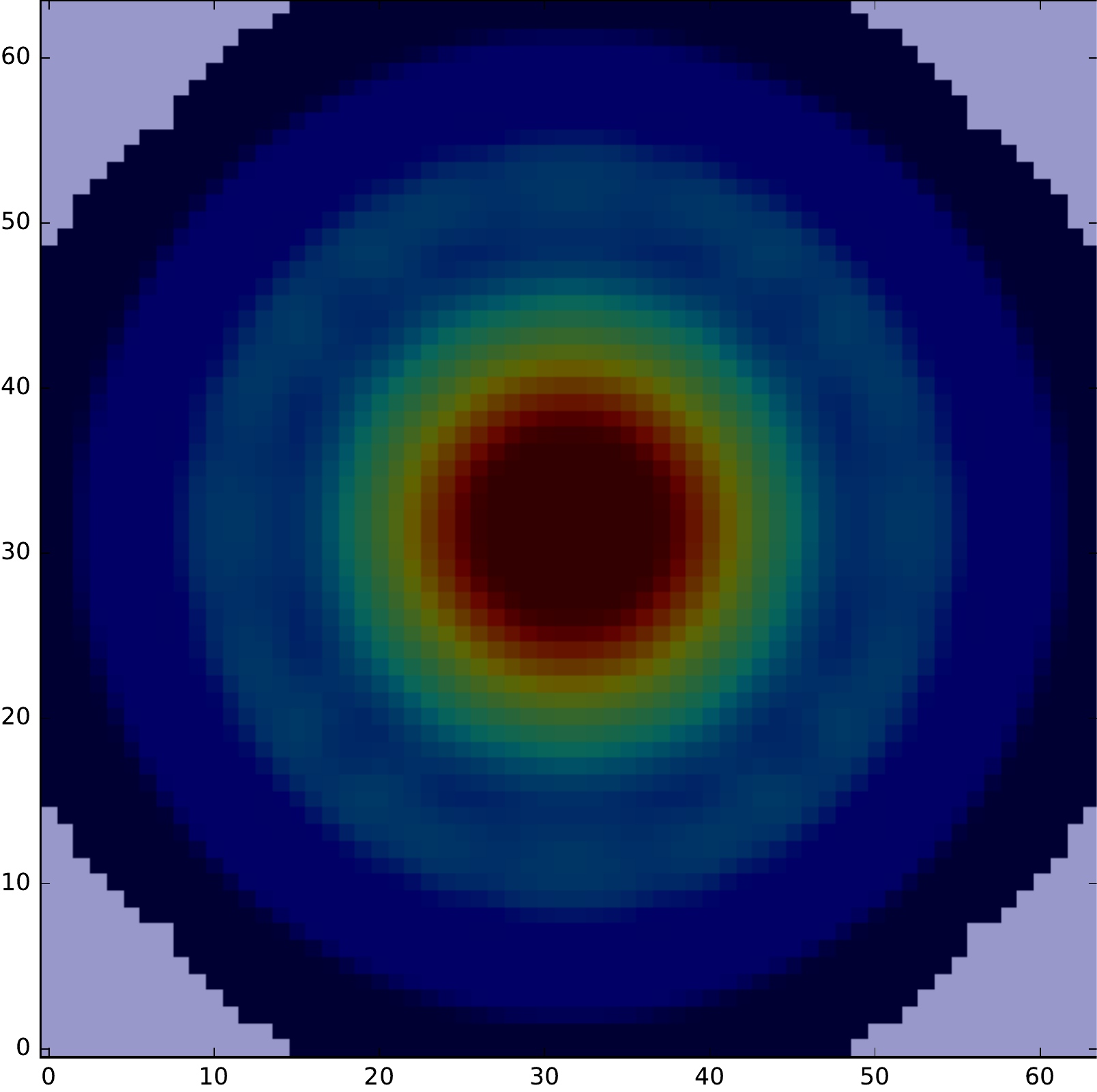}\\
\includegraphics[width=0.23\textwidth]{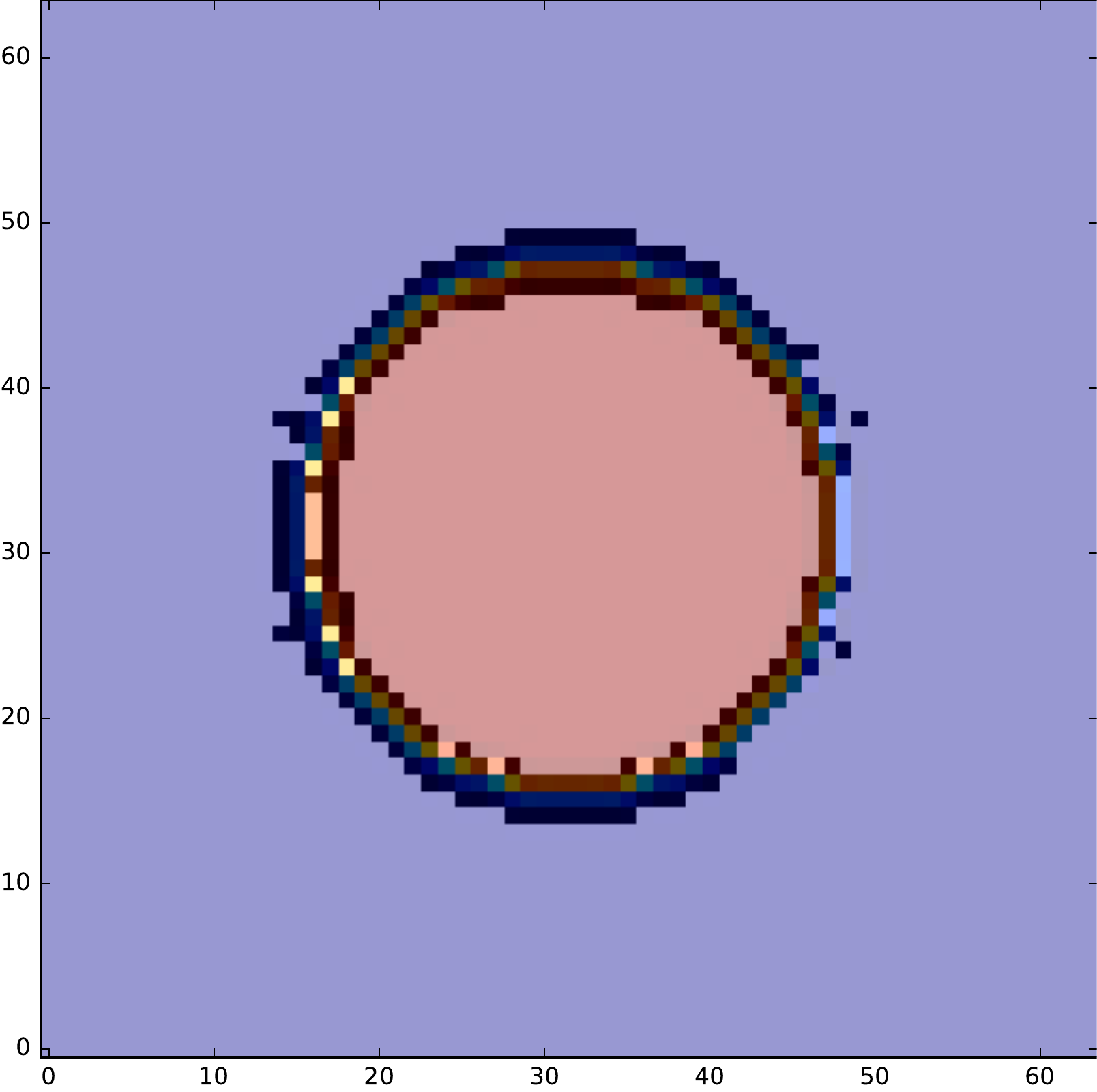}
\includegraphics[width=0.23\textwidth]{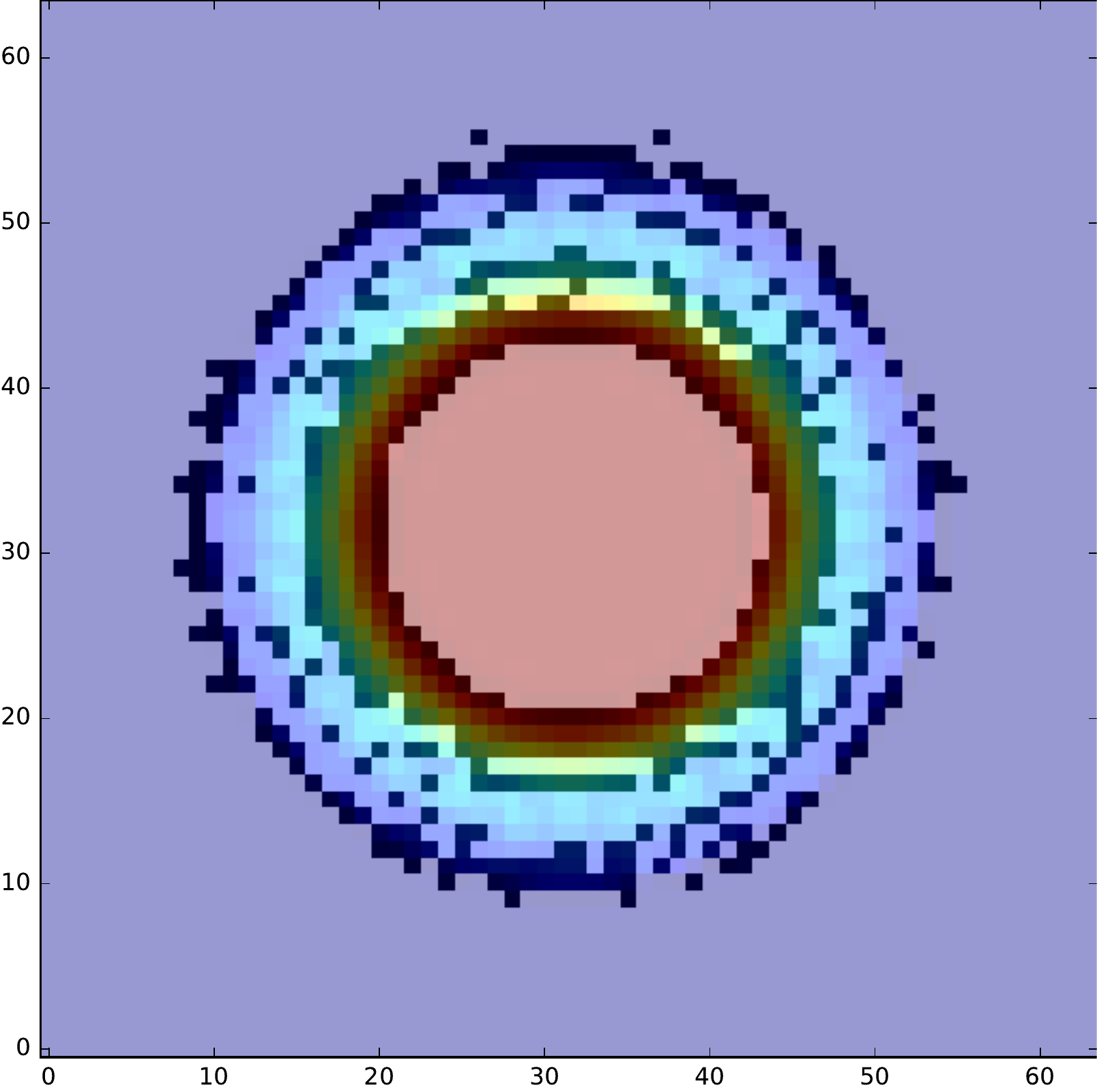}
\includegraphics[width=0.23\textwidth]{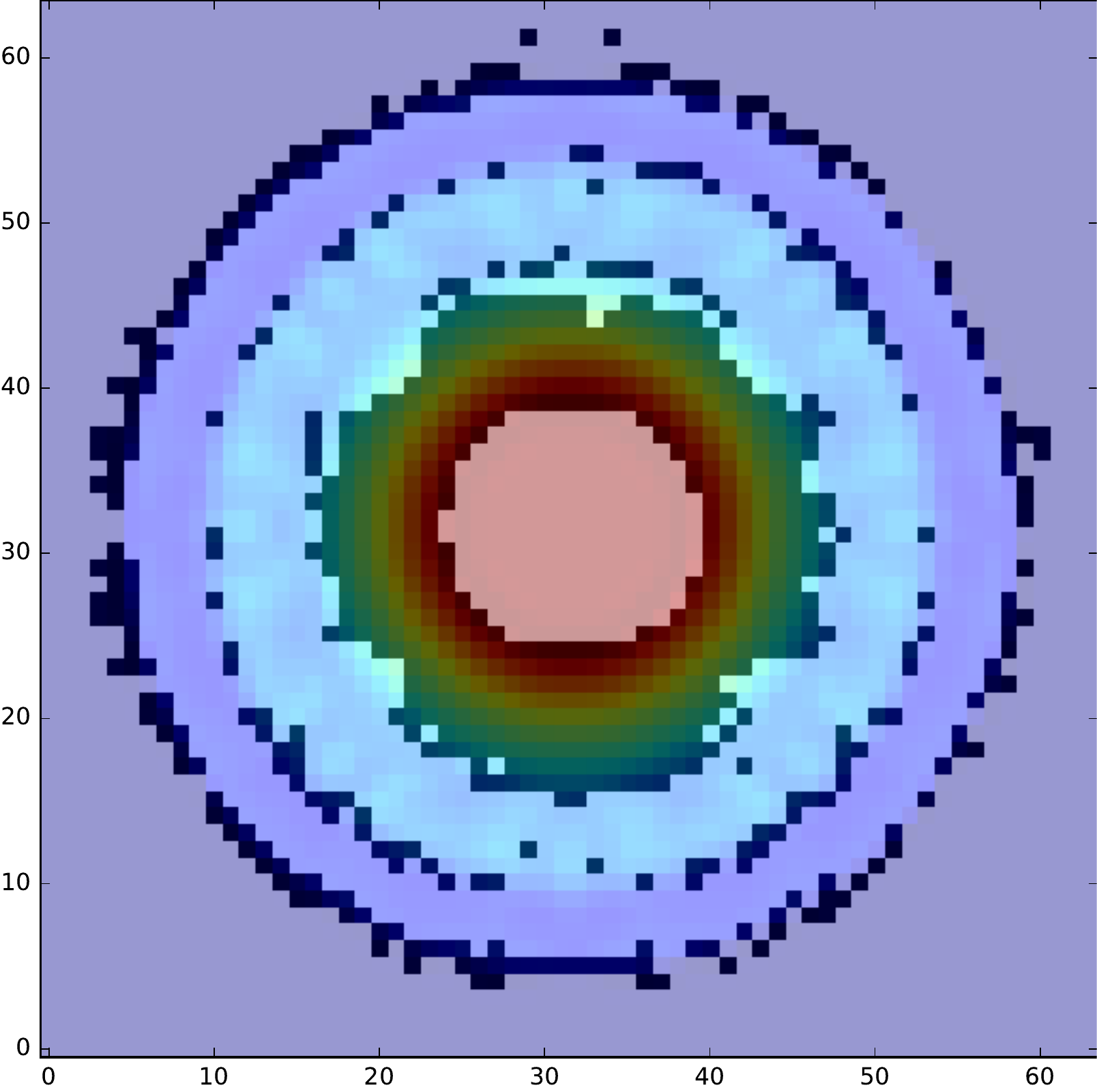}
\includegraphics[width=0.23\textwidth]{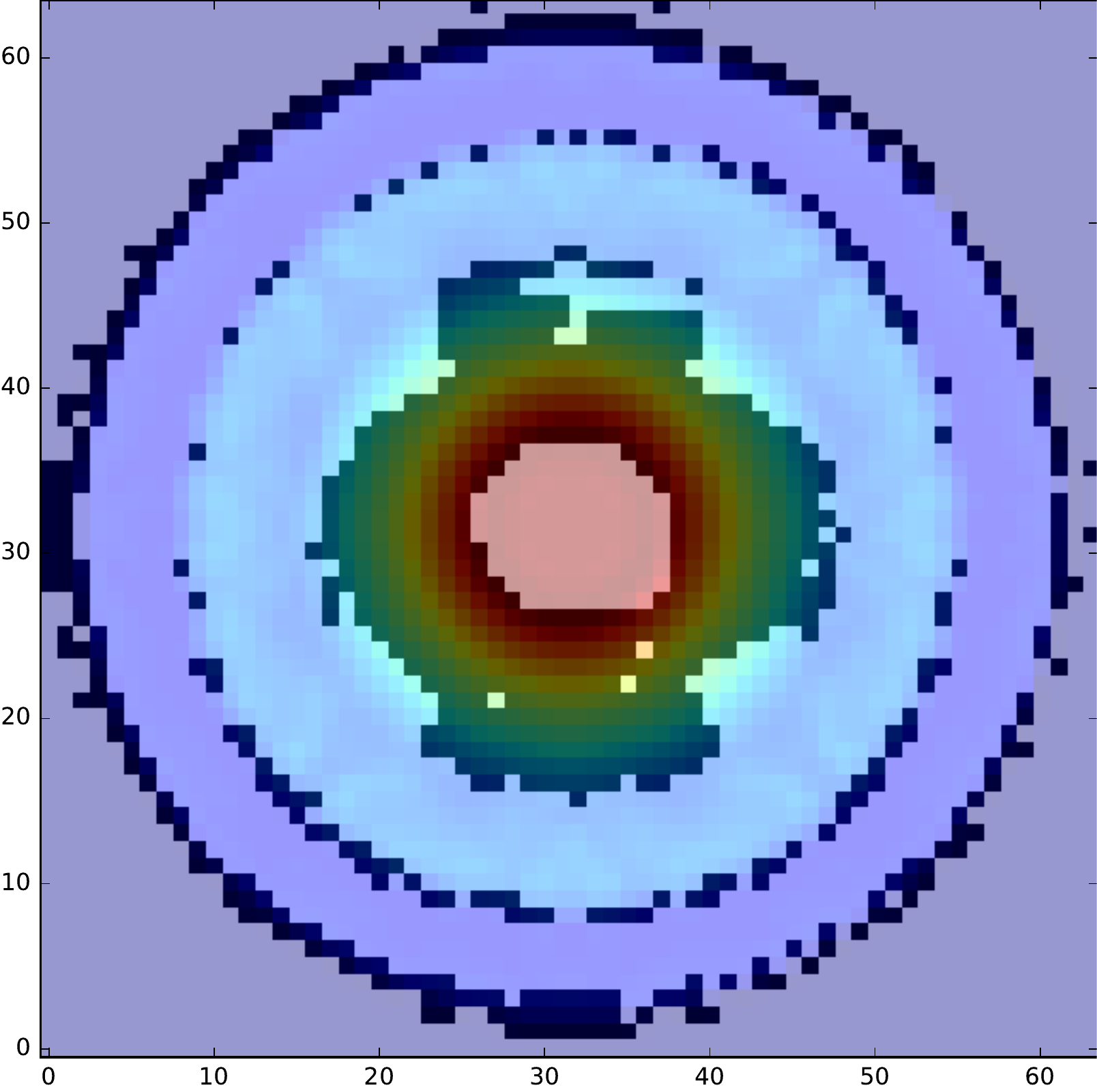}
\caption{Detection comparison between high order limiter (top row) and neural network limiter (bottom row).\label{fig:sod2ddetectioncomparison}}}
\end{figure}

\begin{figure}
\centering{
\includegraphics[width=0.45\textwidth]{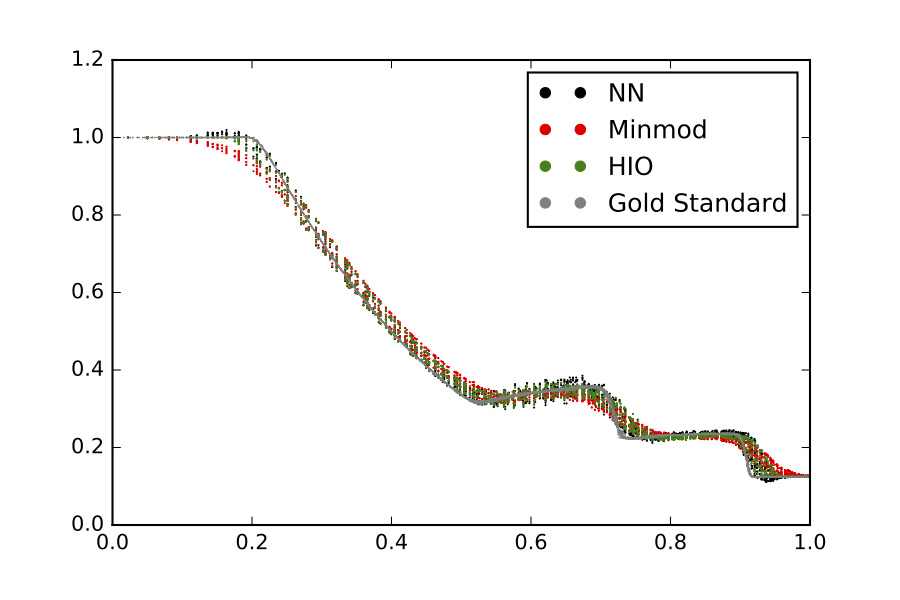}
\includegraphics[width=0.45\textwidth]{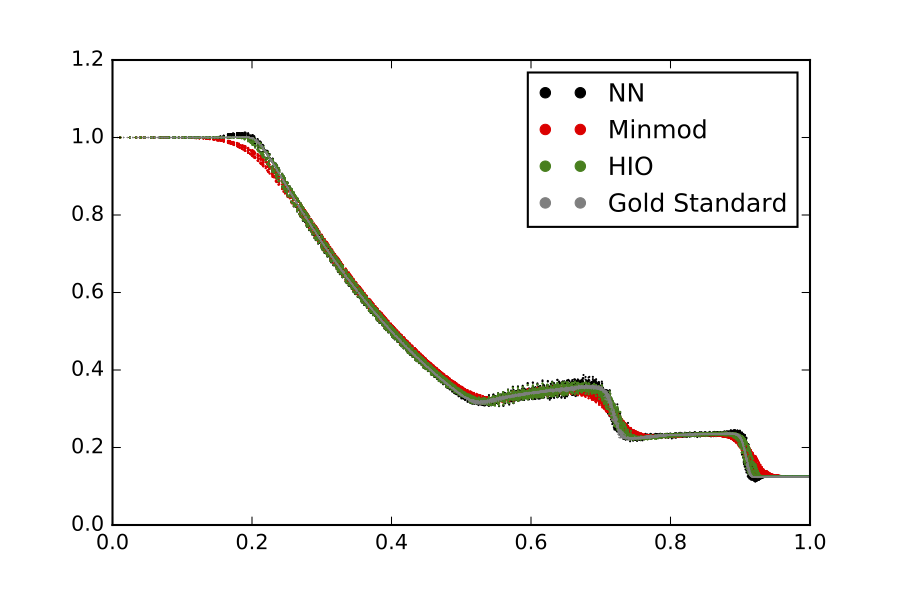}
\caption{Scatter plot of density of the 2-dimensional Sod shock problem at $T=0.24$ for $N_{elem} = 64^2$ and $N_{elem} = 128^2$ and approximation order 2.\label{fig:scatter_sod2d_dg}}}
\end{figure}

\textbf{Case of the Riemann problem 12:} We consider a 2-dimensional Riemann problem (configuration 12) \cite{kurganov}. The initial data are:

\begin{equation}
\label{eq:sod}
(\rho, v_x, v_y, p)(x,0) = \left\{
\begin{array}{ll}
      (1 ,0.7276, 0 , 1 ) & x < 0, y > 0 \\
      (0.8 ,0 ,0 , 1) & x < 0, y < 0 \\
      (0.5313, 0, 0, 0.4) & x > 0, y > 0 \\
      (1,0,0.7276,1) & x > 0, y < 0
\end{array} 
\right.
\end{equation}
for $\gamma = 1.4$ and gradient-free boundary conditions.

\begin{figure}
\centering{
\includegraphics[width=0.30\textwidth]{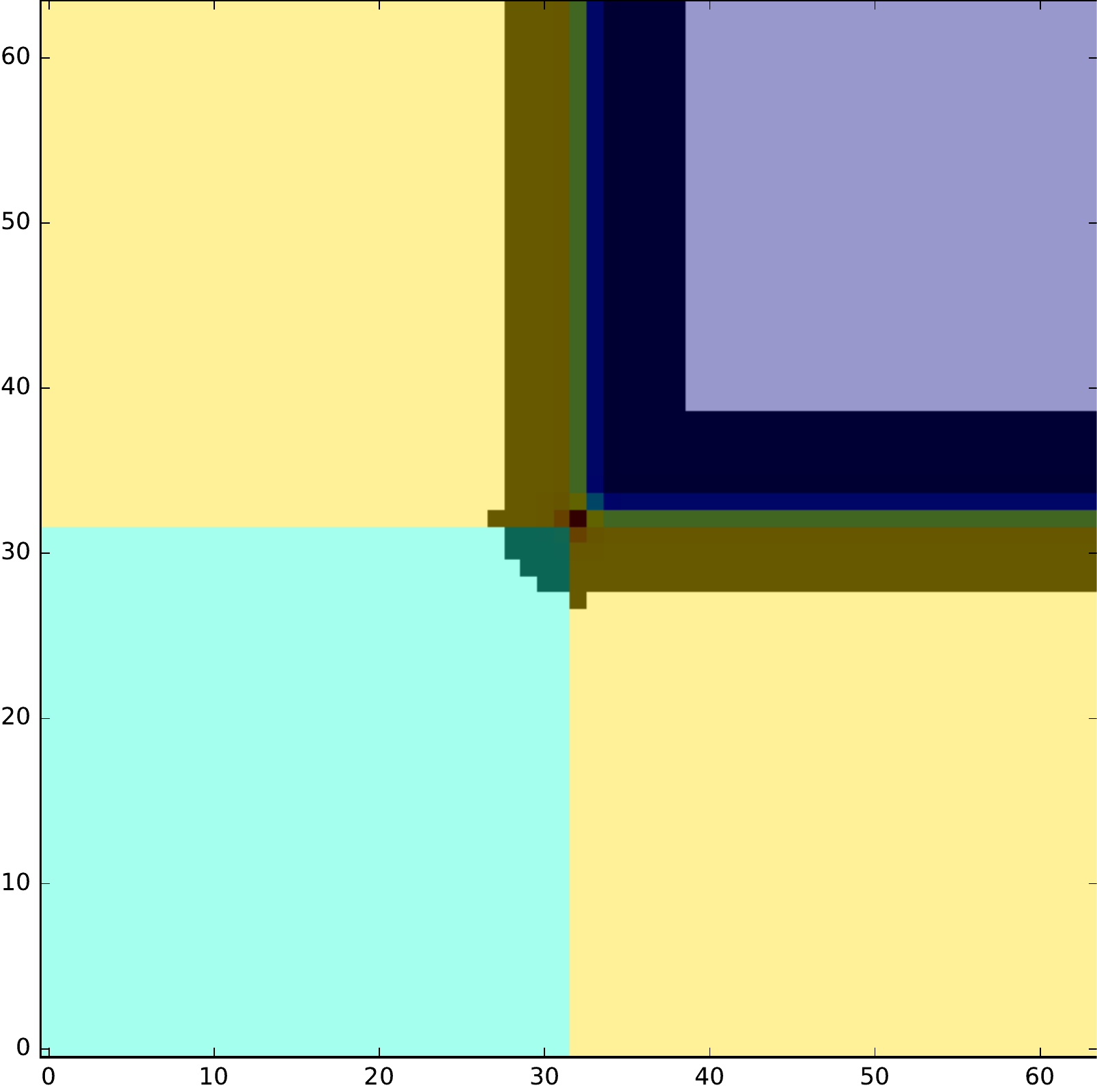}
\includegraphics[width=0.30\textwidth]{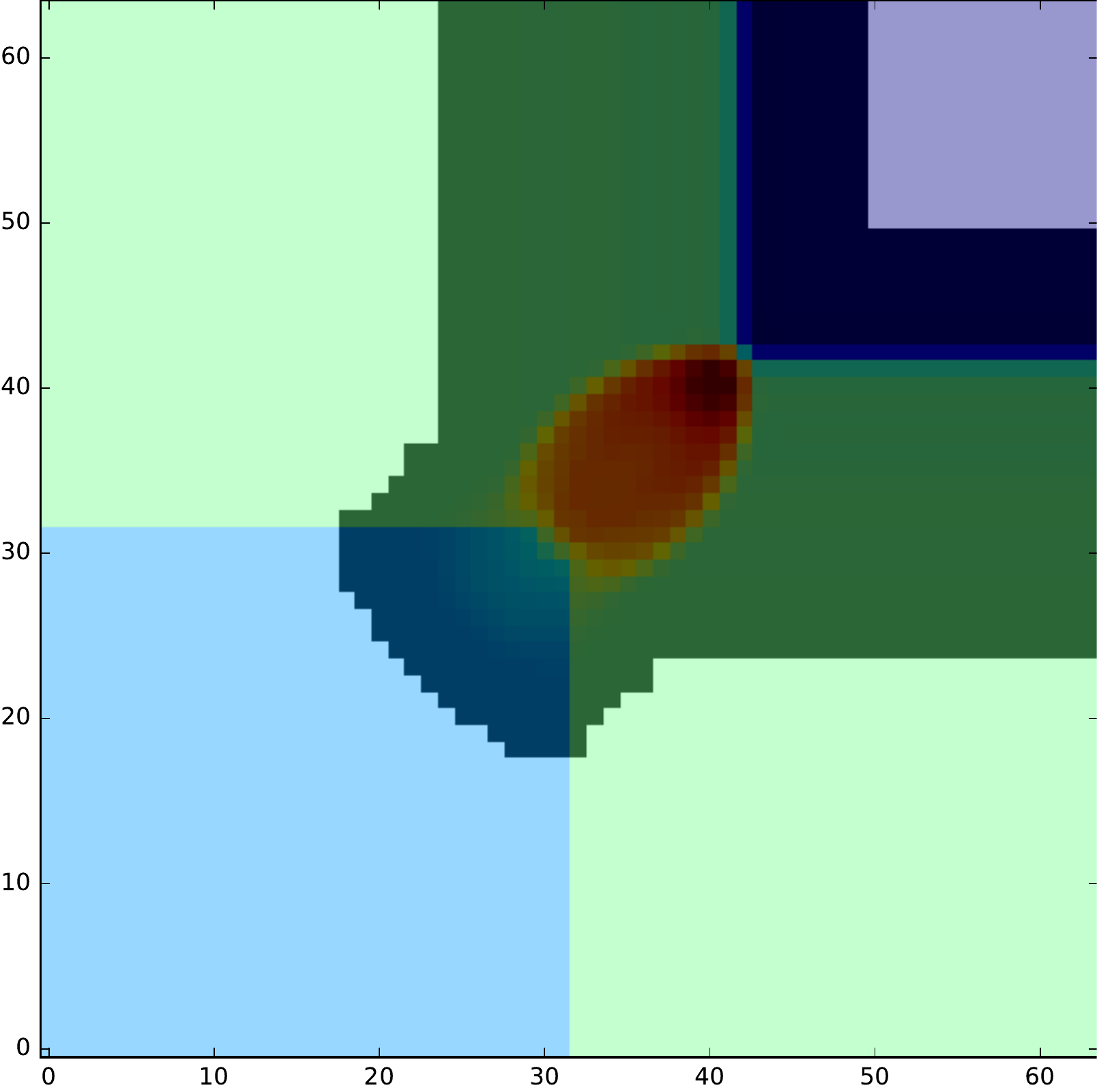}
\includegraphics[width=0.30\textwidth]{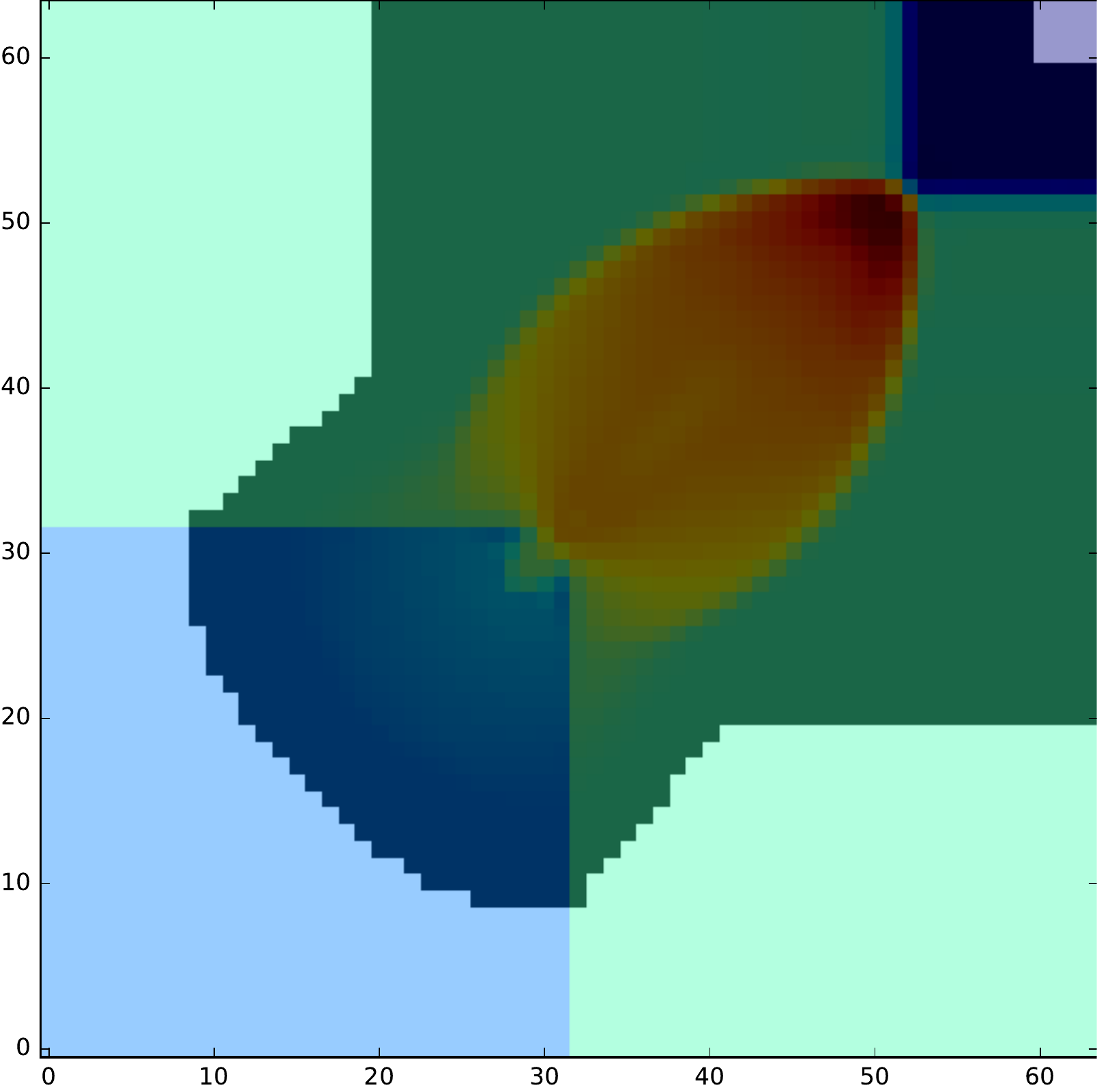}\\
\includegraphics[width=0.30\textwidth]{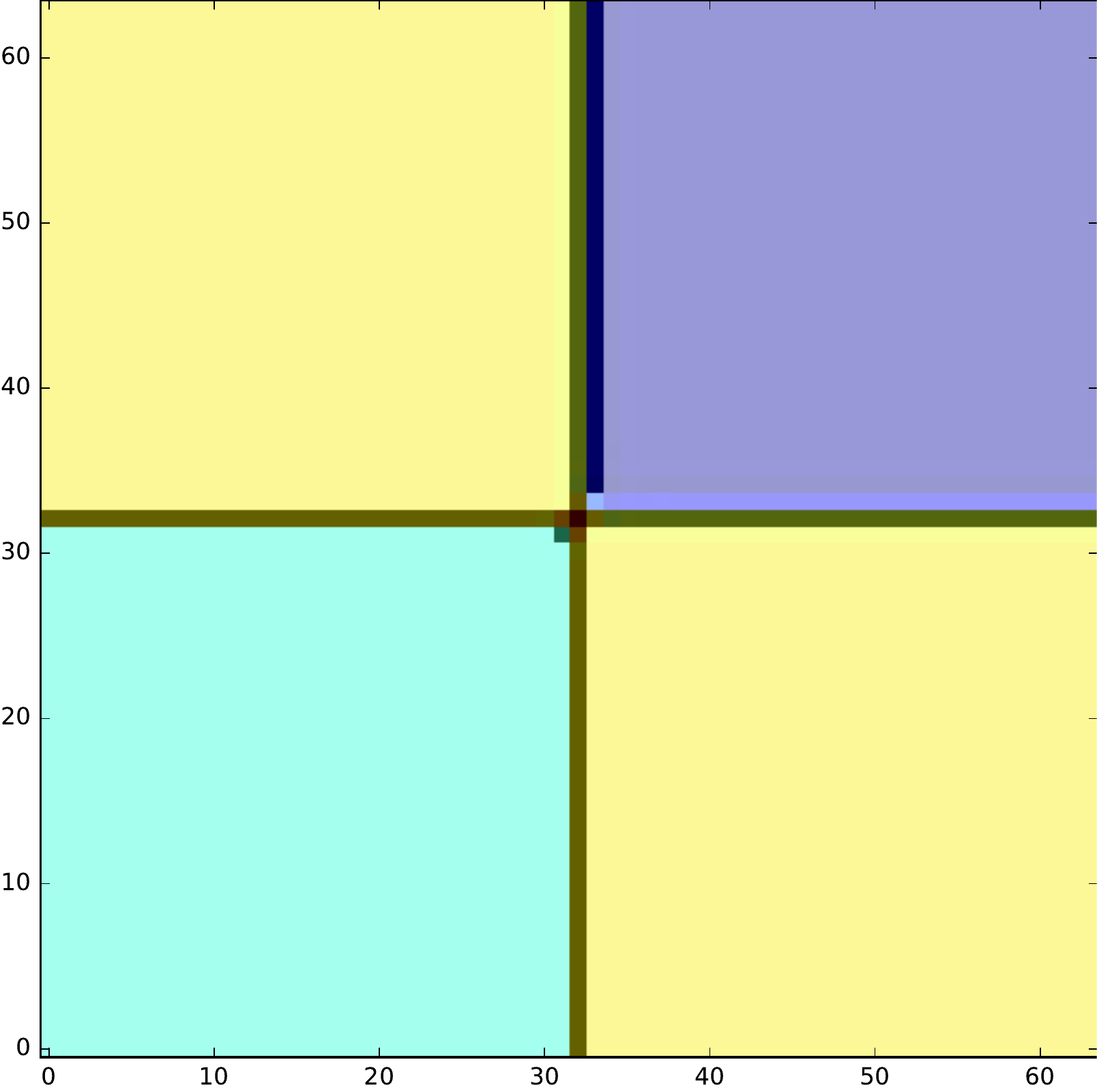}
\includegraphics[width=0.30\textwidth]{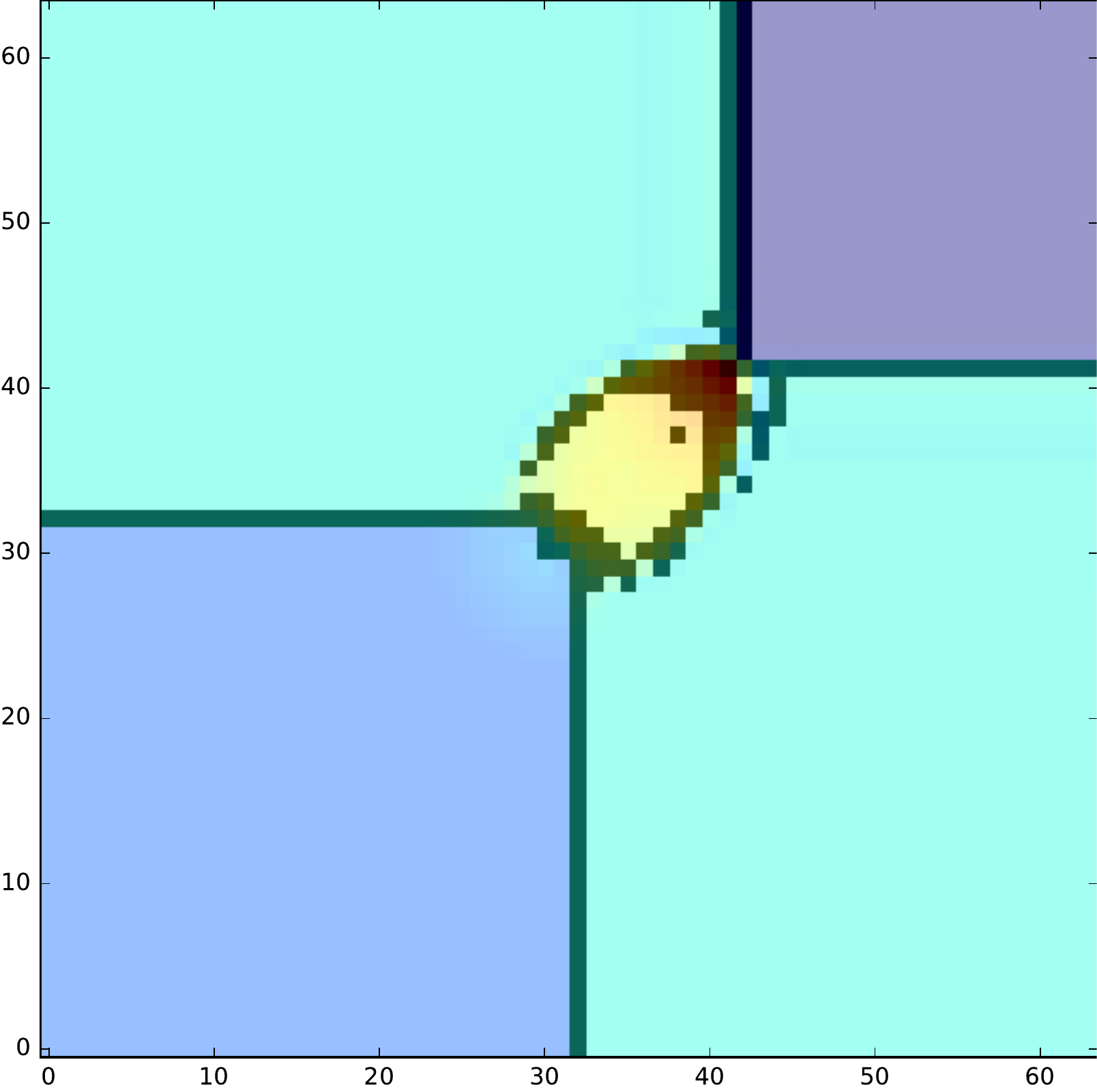}
\includegraphics[width=0.30\textwidth]{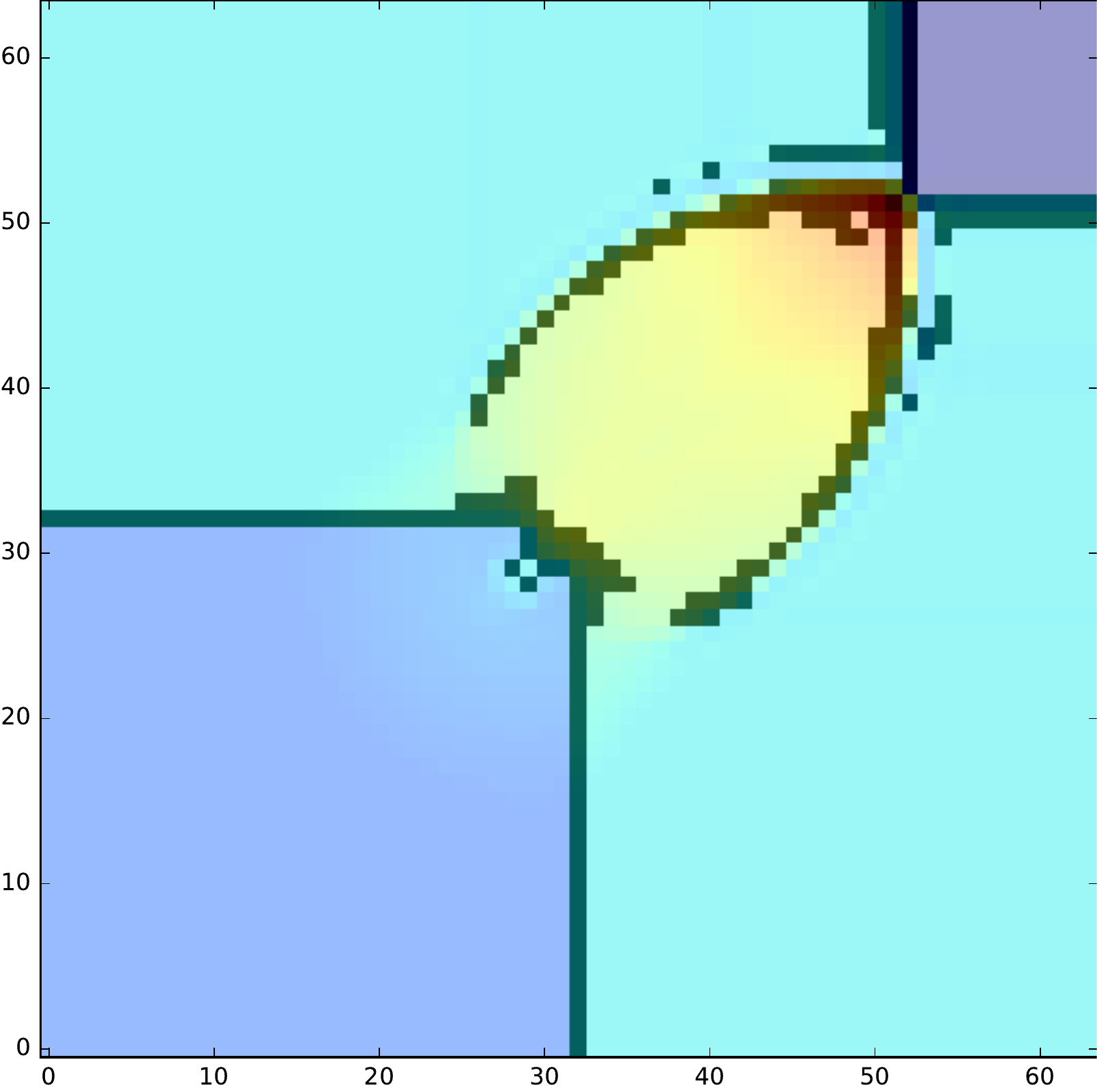}
\caption{Detection comparison between high order limiter (top row) and neural network limiter (bottom row).\label{fig:rieman2ddetectioncomparison}}}
\end{figure}

\subsection{Transfer to residual distribution}
\label{sec:transfer2d}
In this section we report the results of the neural network limiter applied to a residual distribution scheme. In particular, we report on the performance of the neural based limiter when no retraining phase is performed versus the performance of a neural based limiter that is retrained on a reduced dataset using data from numerical runs using the residual distribution scheme.

\subsubsection{2-d Sod shock tube:}

\begin{figure}
\centering{
\includegraphics[width=0.23\textwidth]{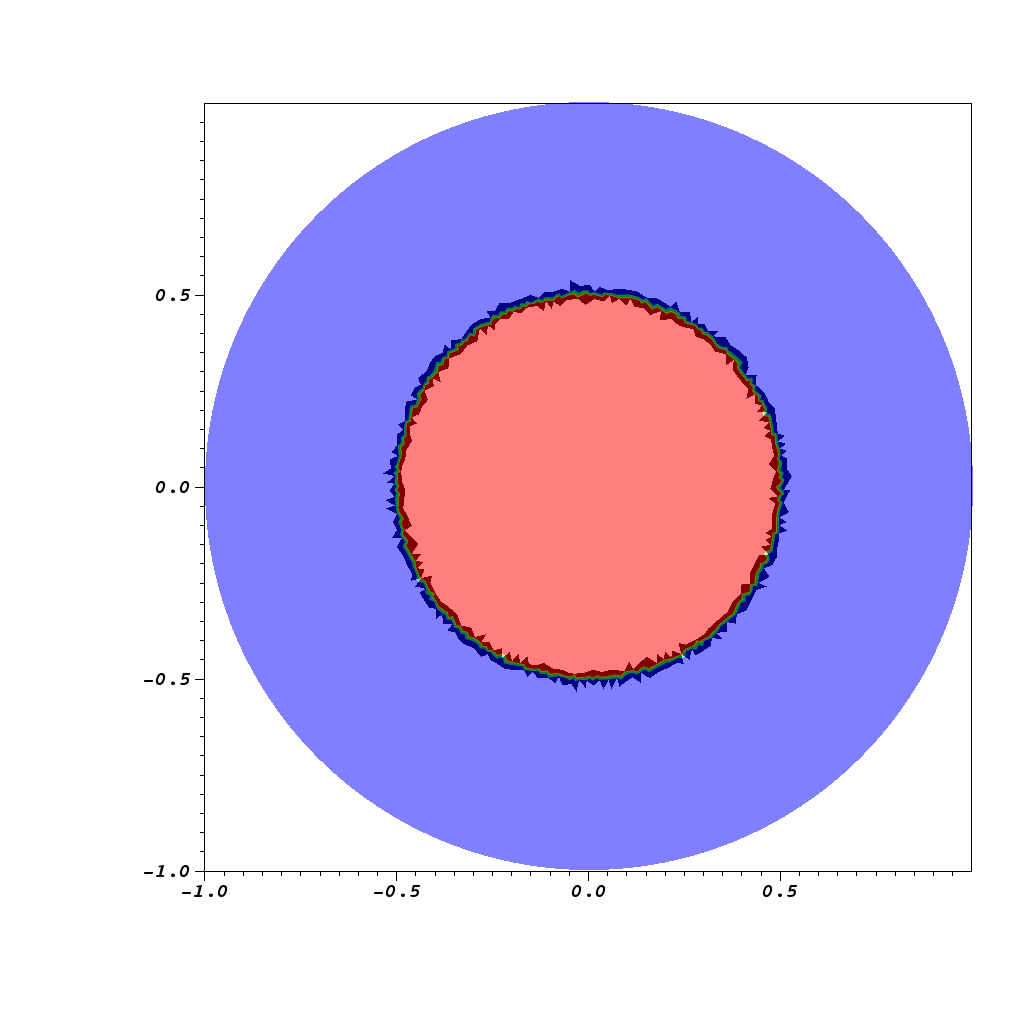}
\includegraphics[width=0.23\textwidth]{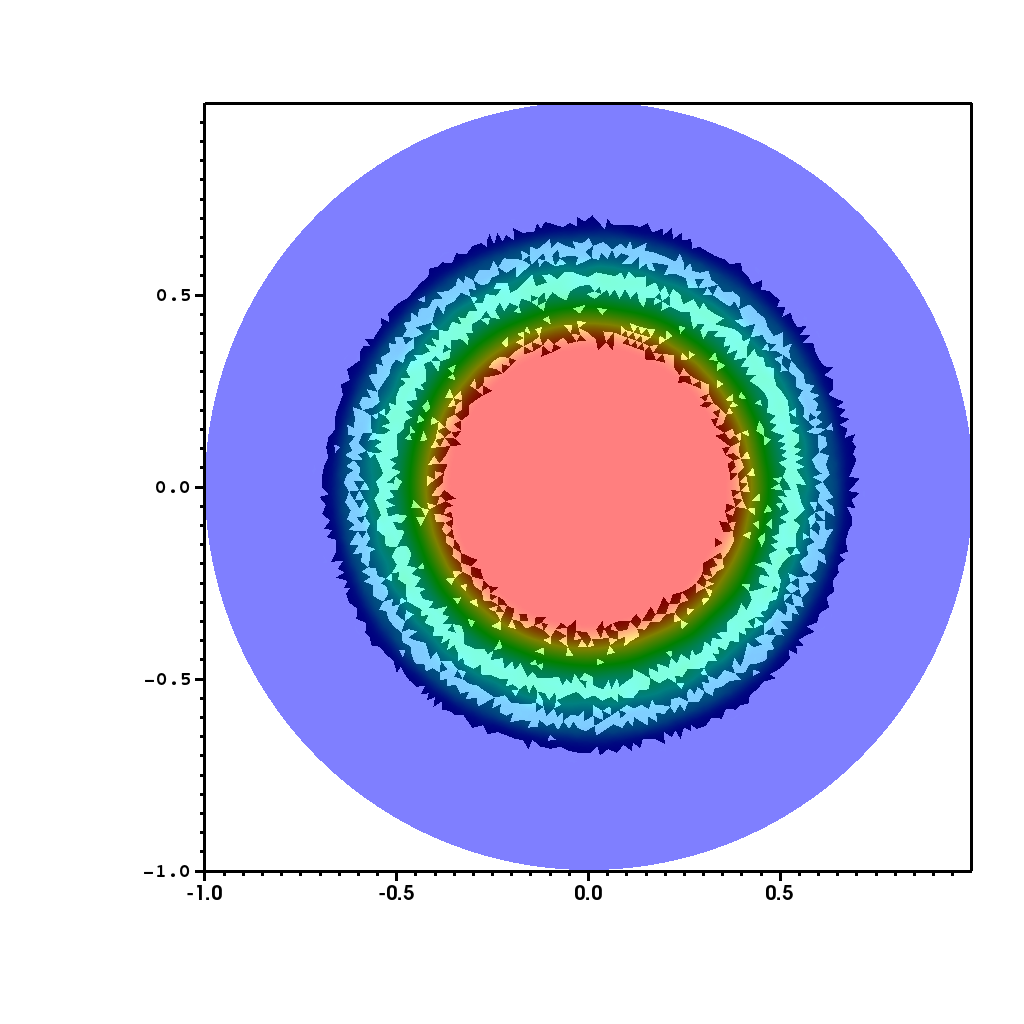}
\includegraphics[width=0.23\textwidth]{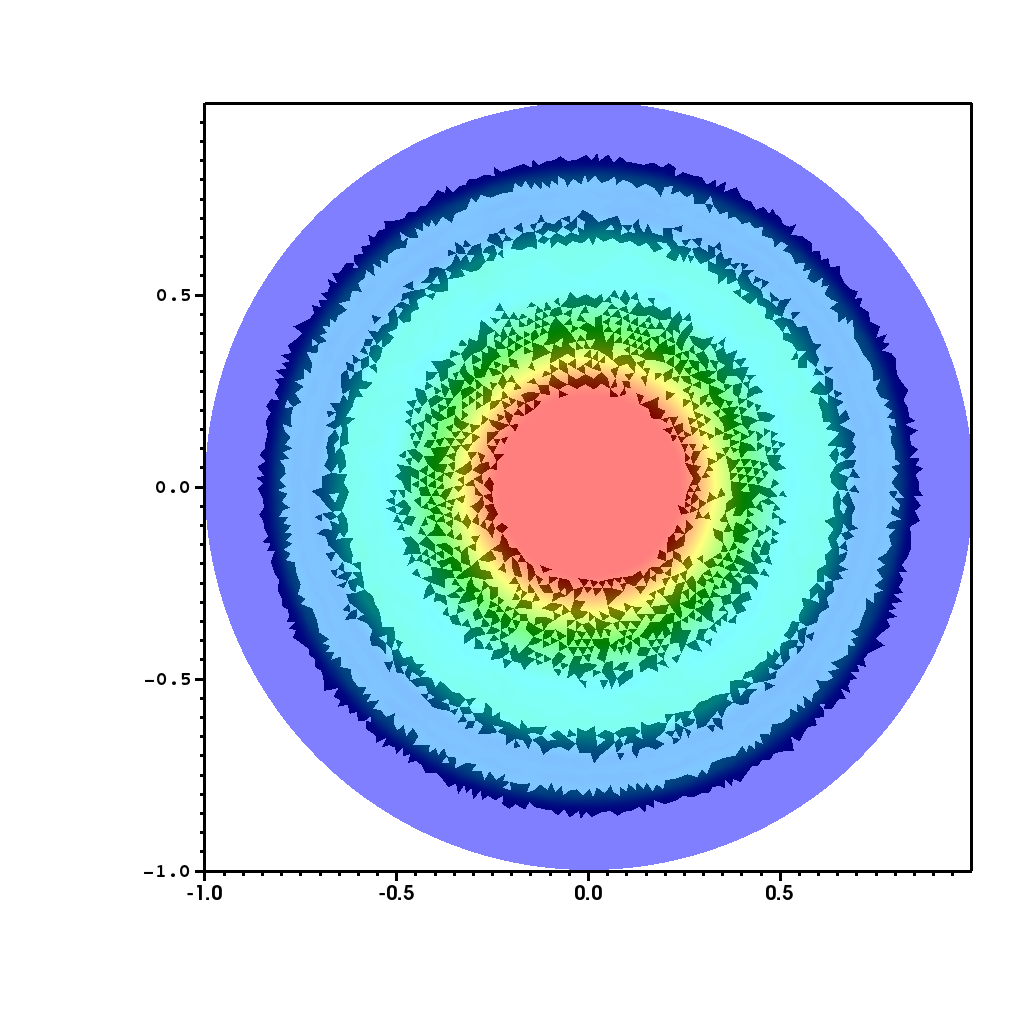}
\includegraphics[width=0.23\textwidth]{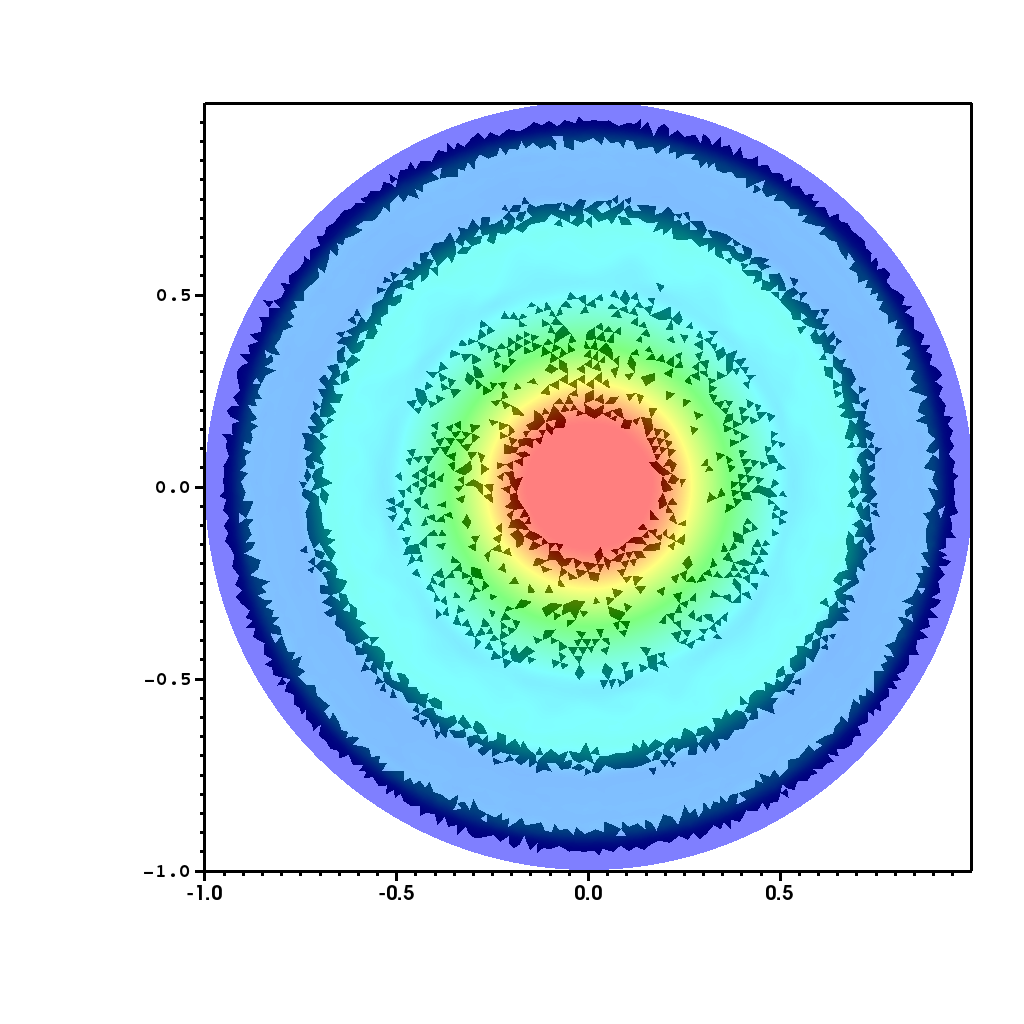}\\
\includegraphics[width=0.23\textwidth]{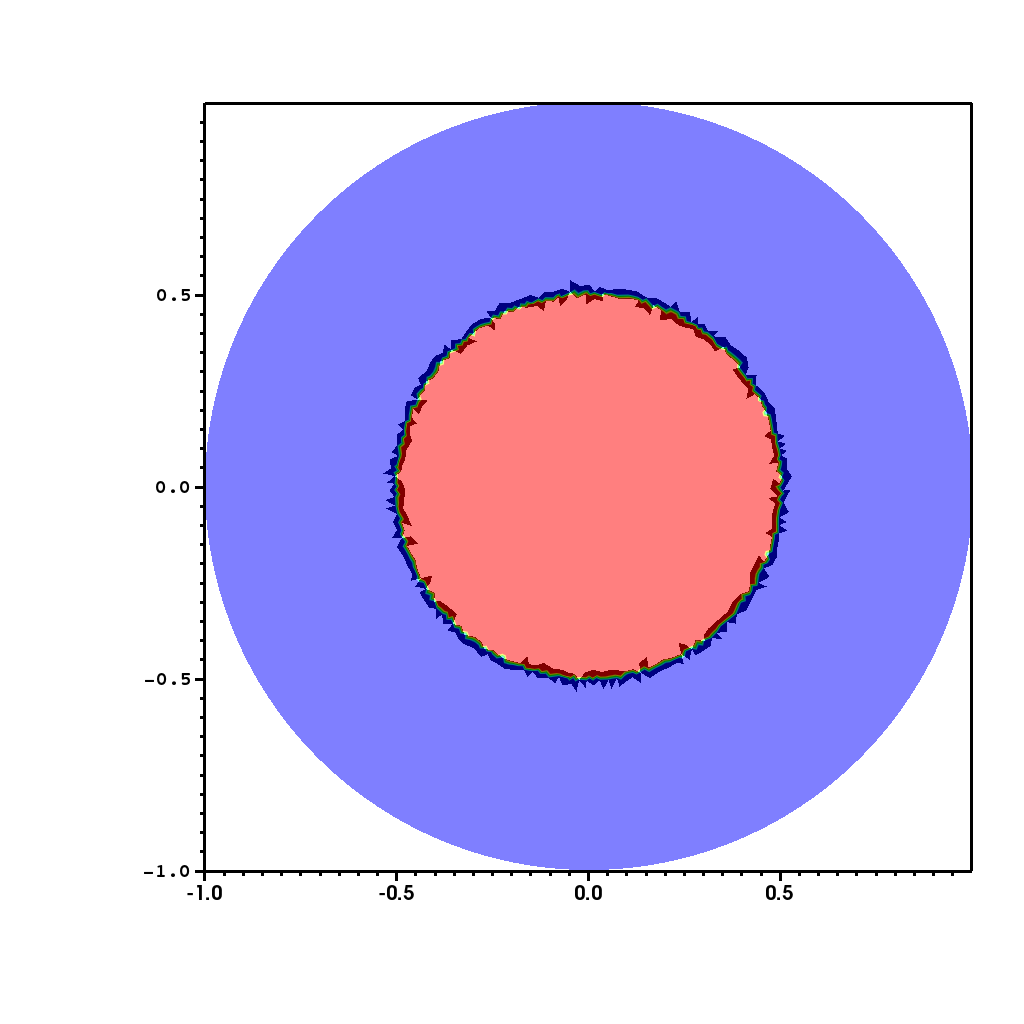}
\includegraphics[width=0.23\textwidth]{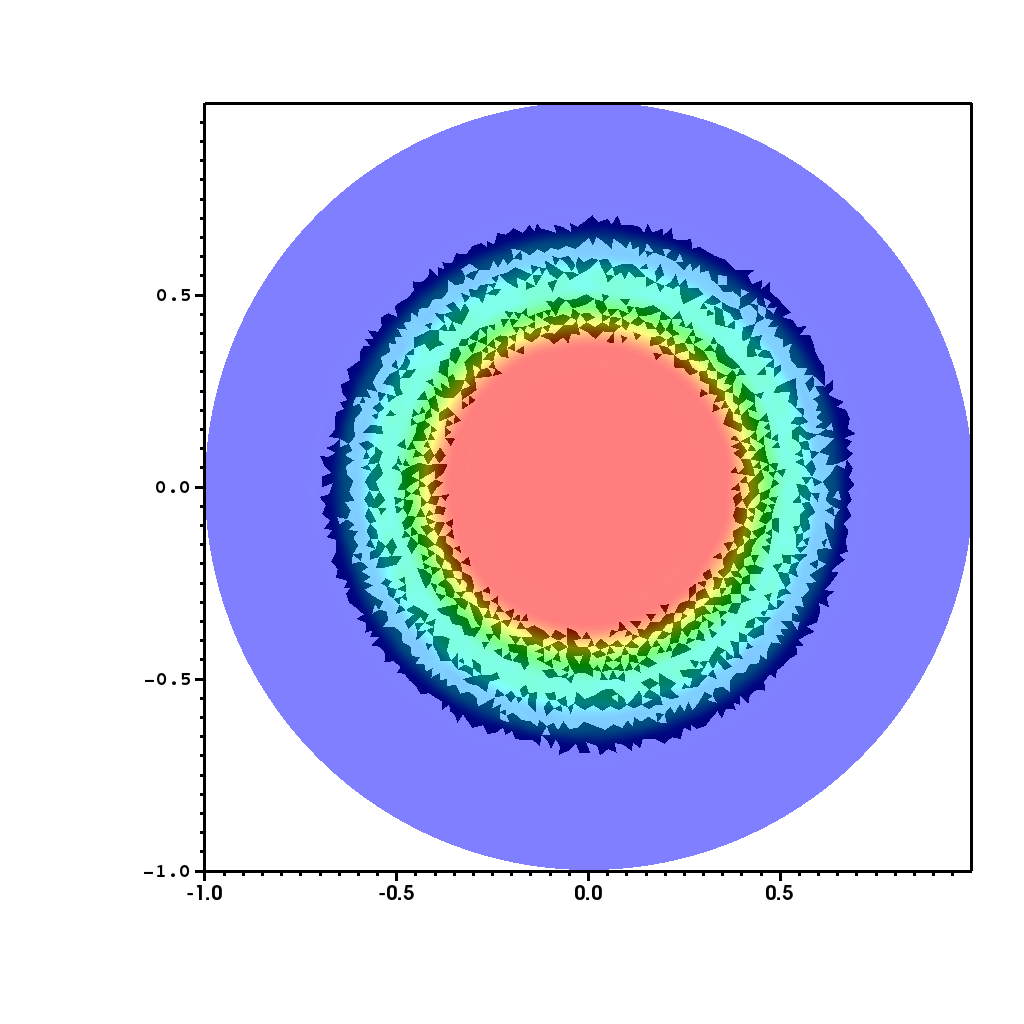}
\includegraphics[width=0.23\textwidth]{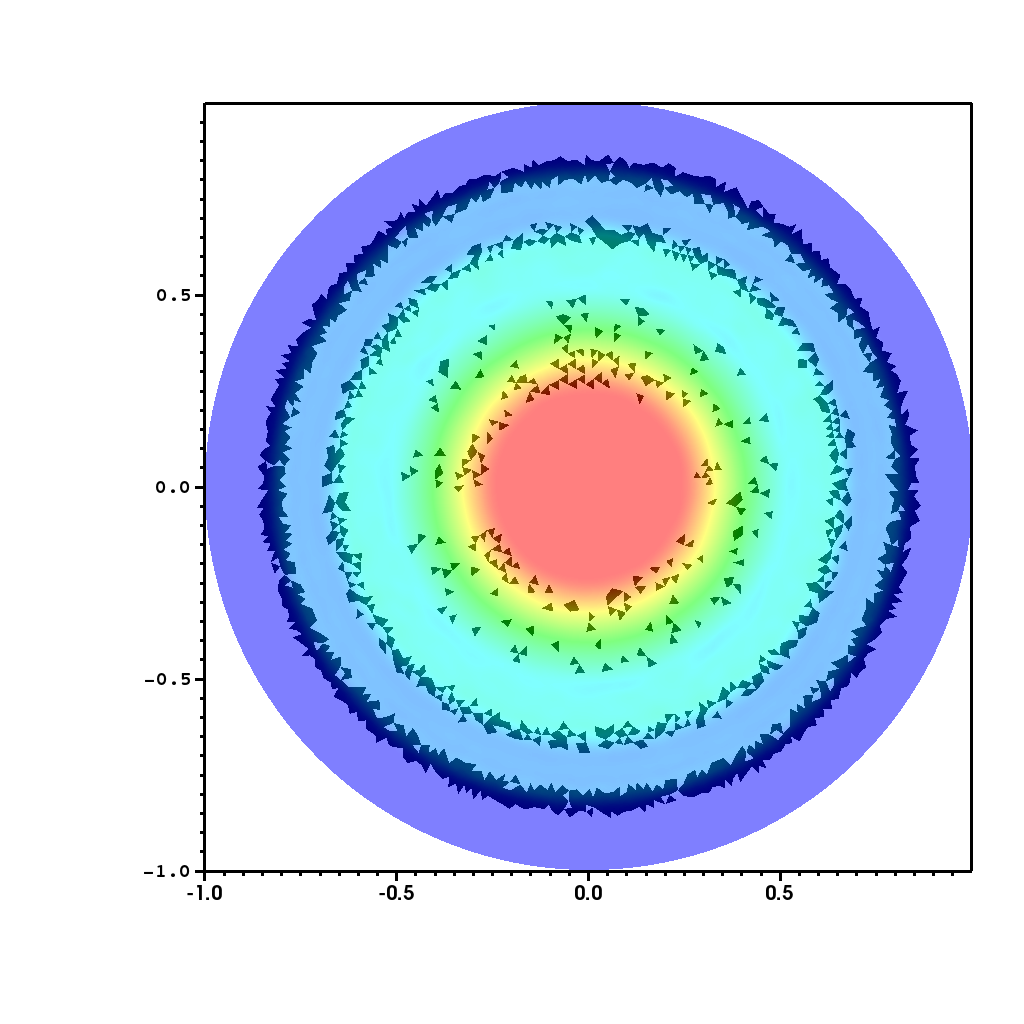}
\includegraphics[width=0.23\textwidth]{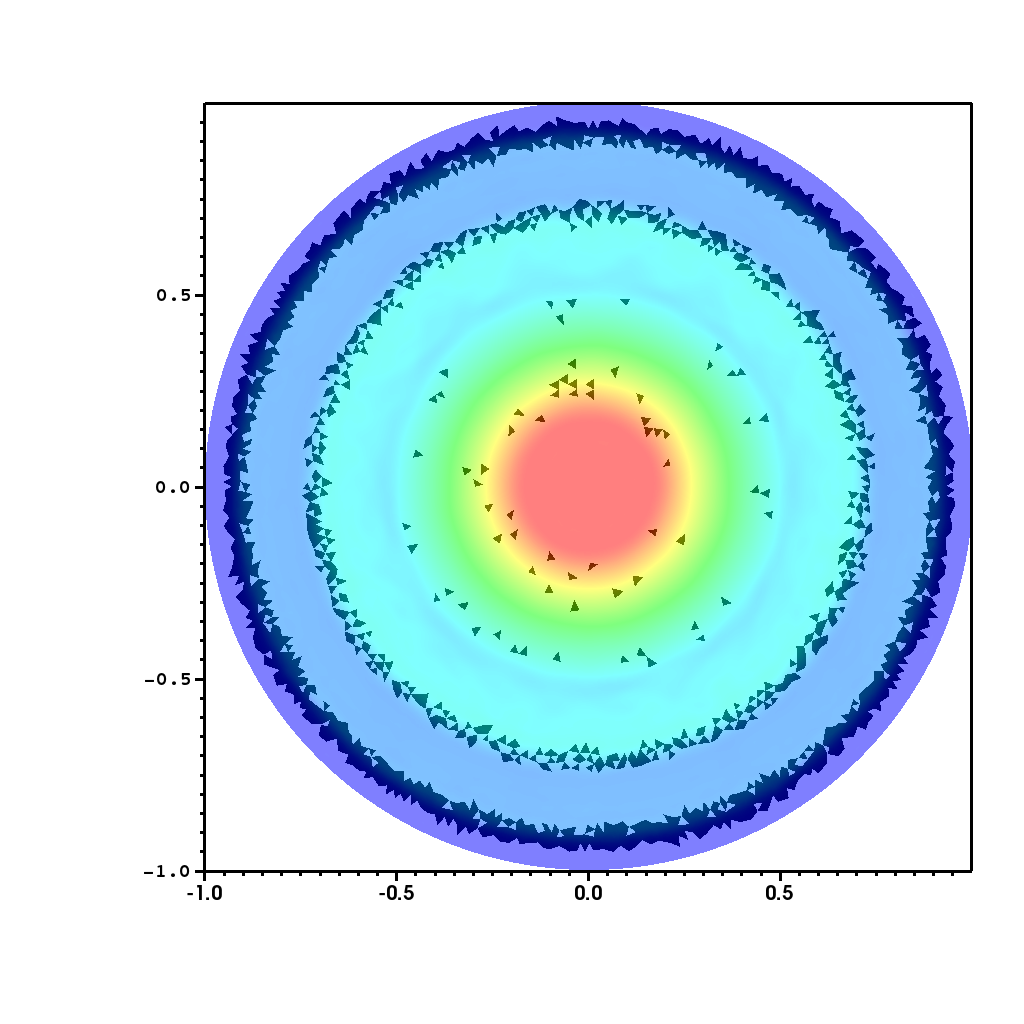}
\caption{Detection comparison between the neural network limiter without retraining (top row) and with retraining (bottom row) on an unstructured mesh.\label{fig:sod2ddetectioncomparisontrian}}}
\end{figure}

\begin{figure}
\centering{
\includegraphics[width=0.23\textwidth]{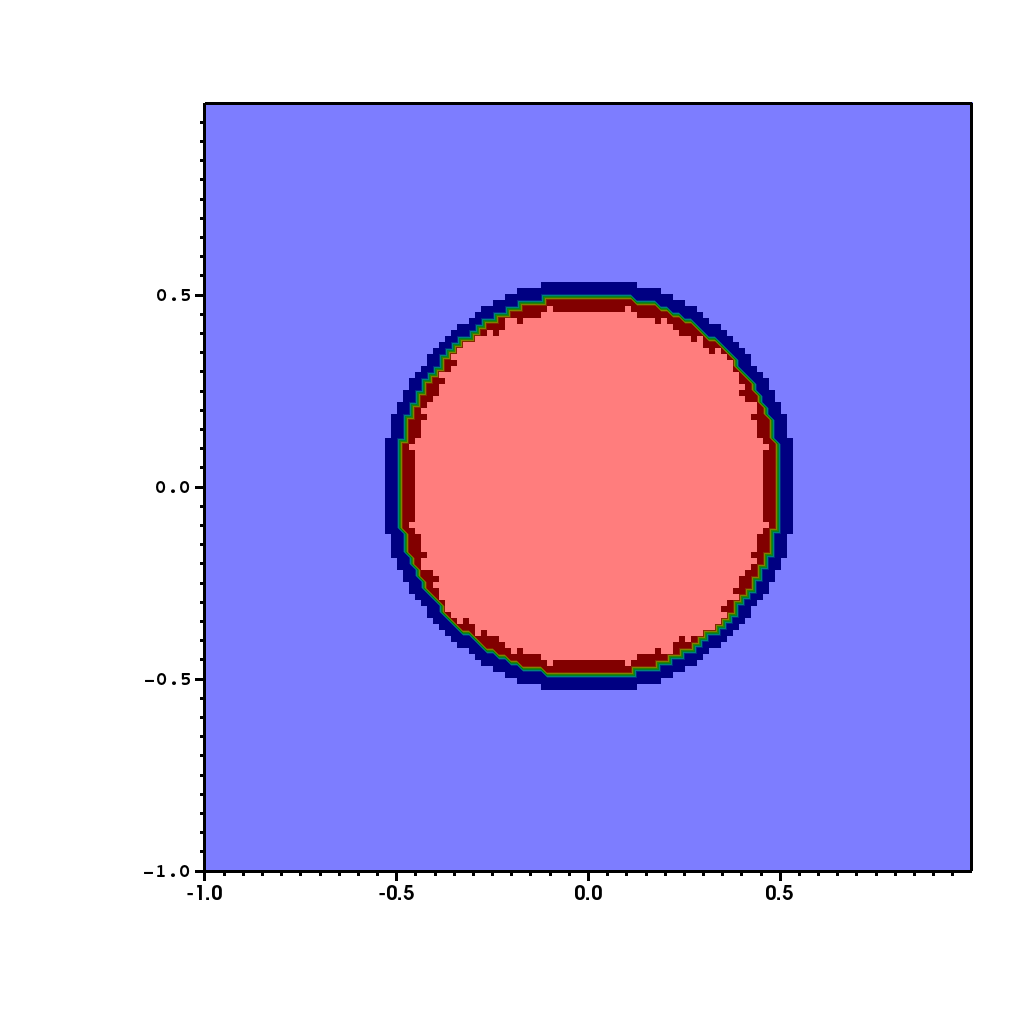}
\includegraphics[width=0.23\textwidth]{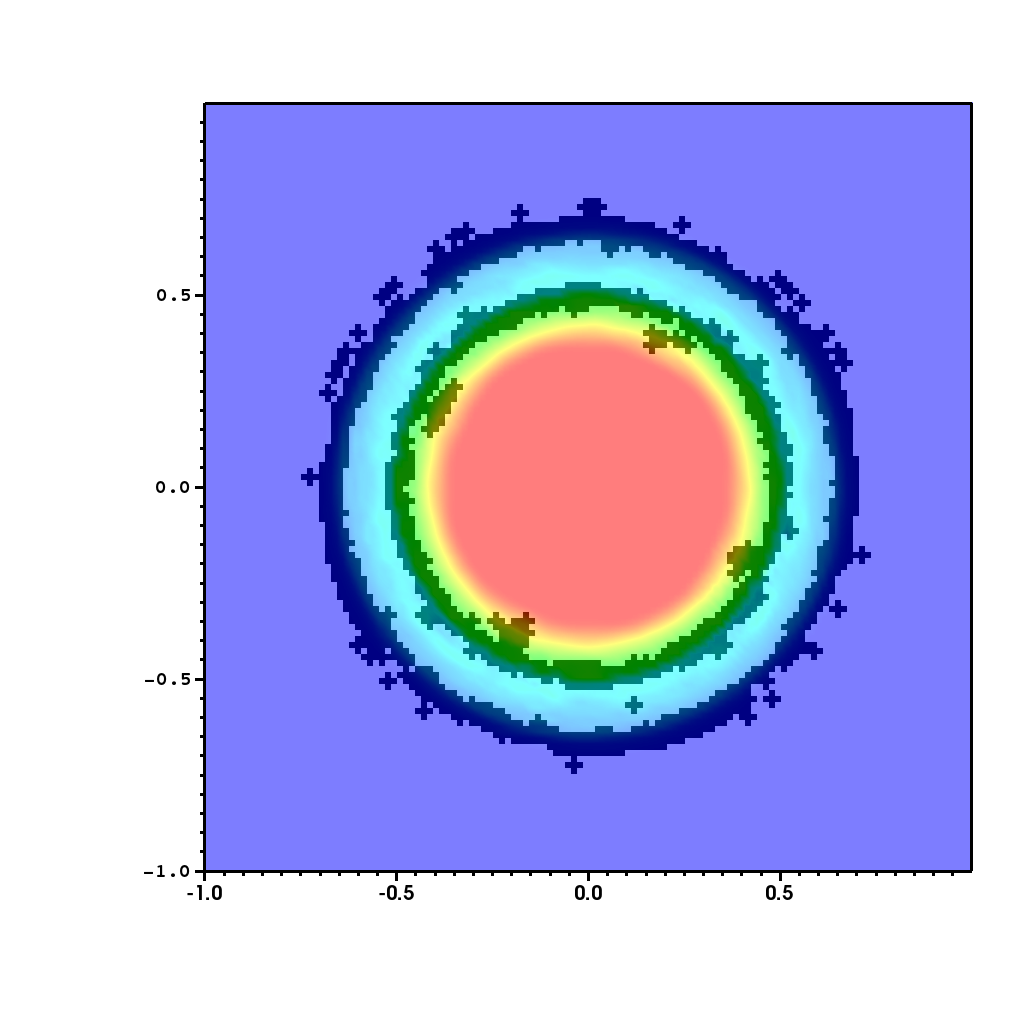}
\includegraphics[width=0.23\textwidth]{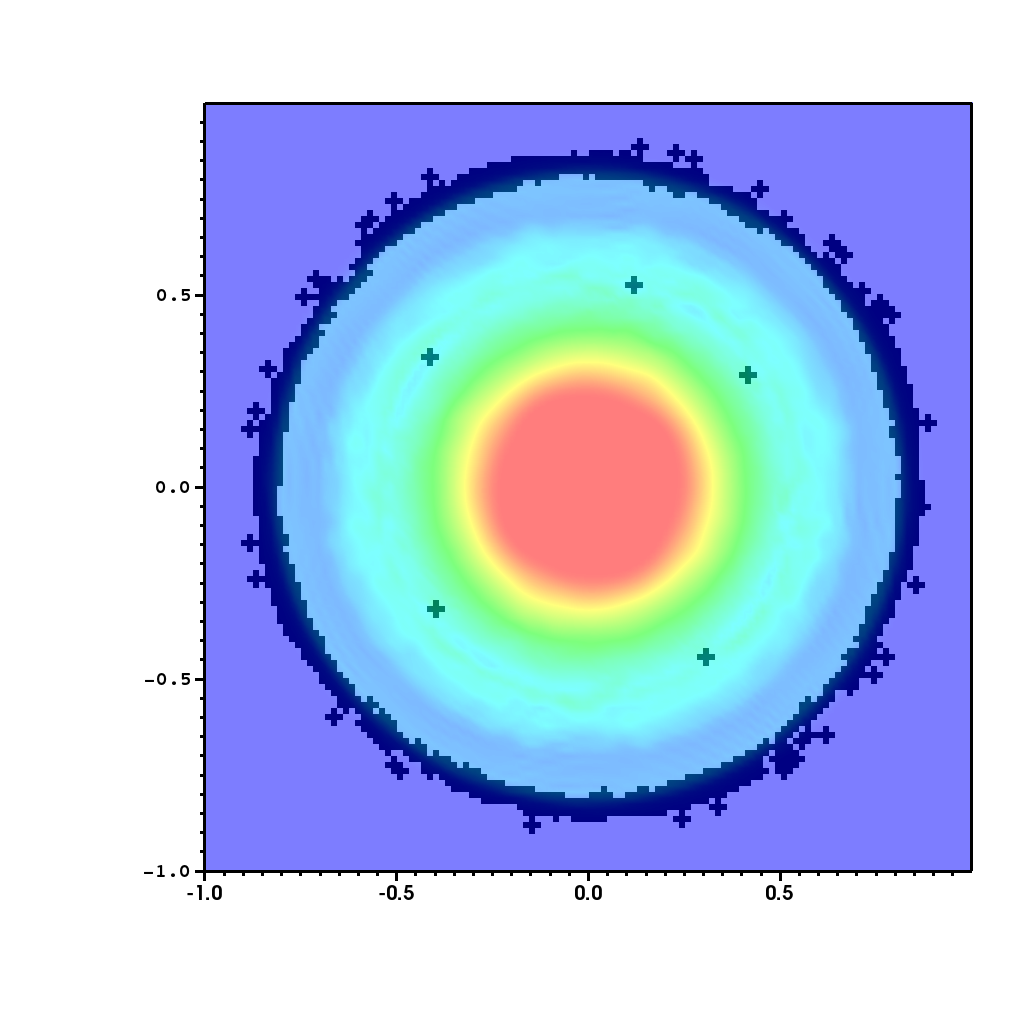}
\includegraphics[width=0.23\textwidth]{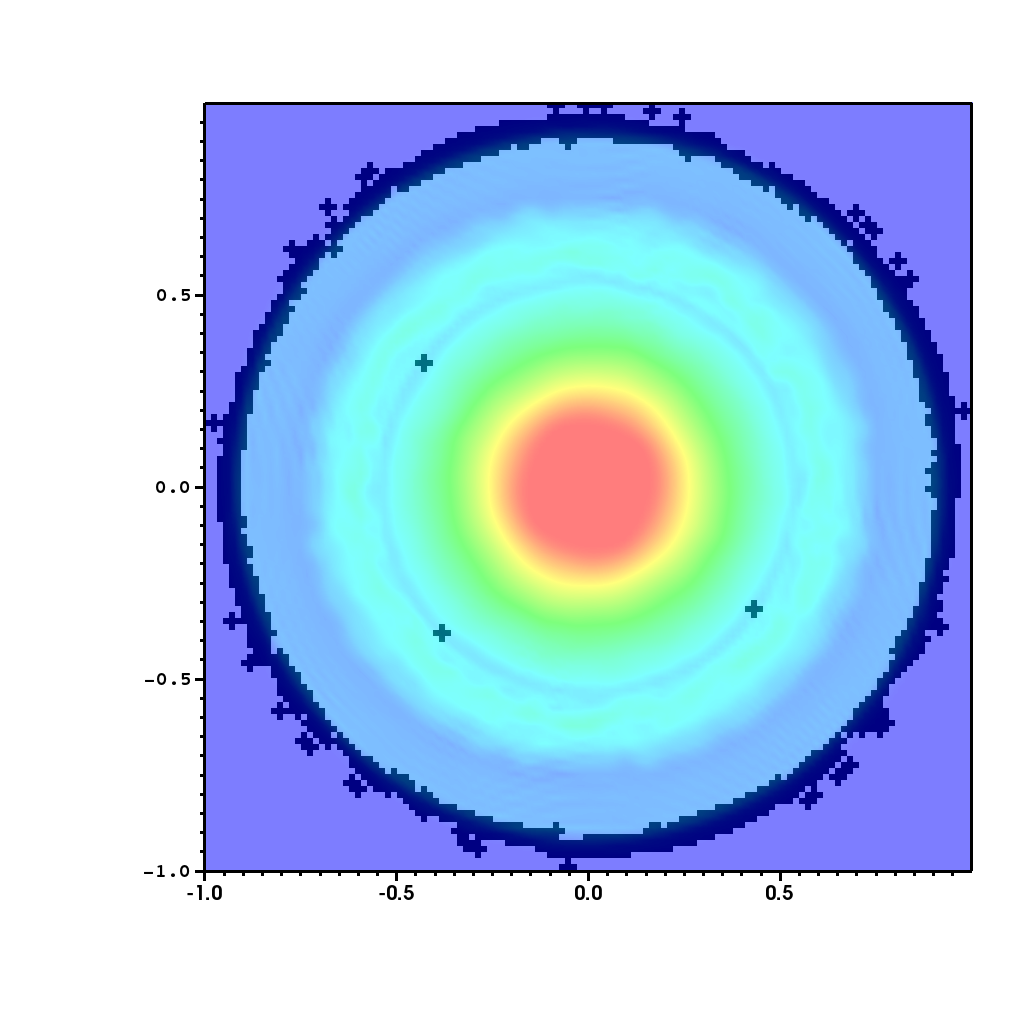}\\
\includegraphics[width=0.23\textwidth]{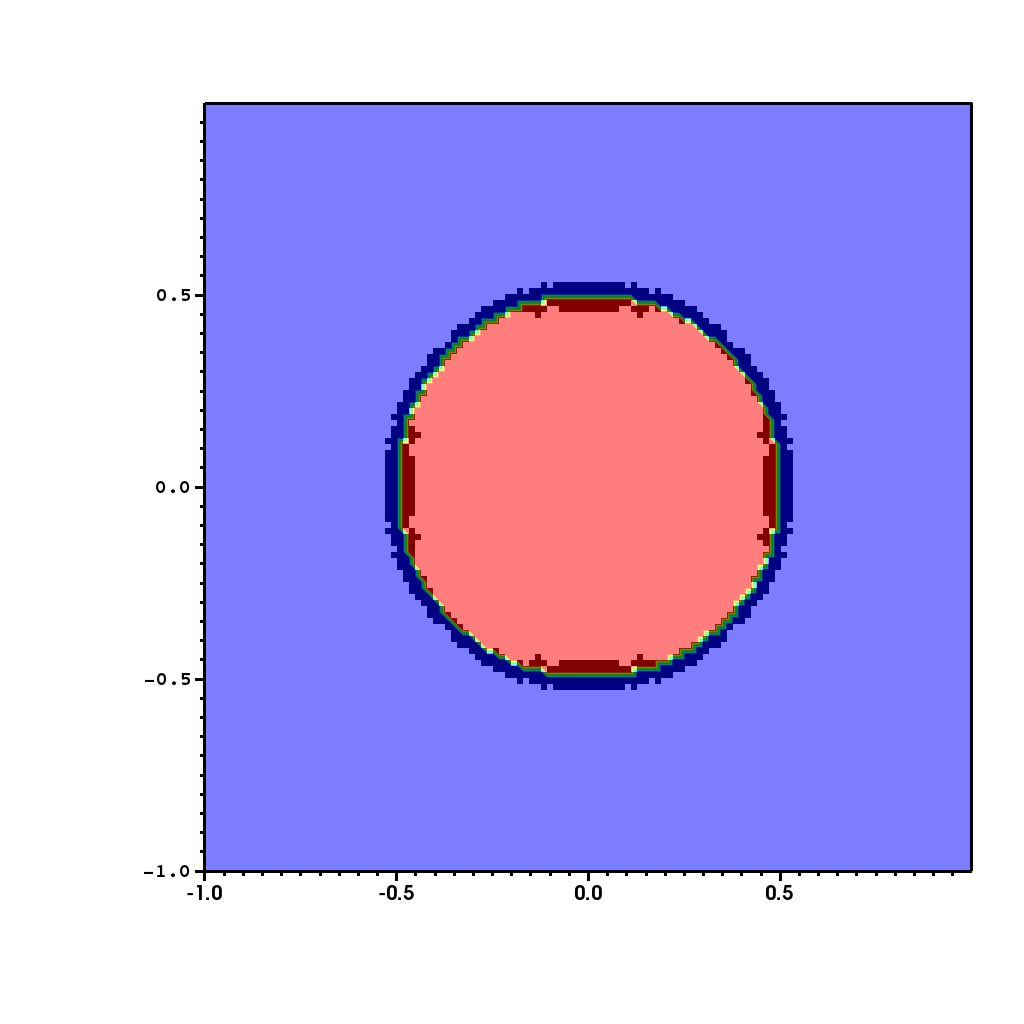}
\includegraphics[width=0.23\textwidth]{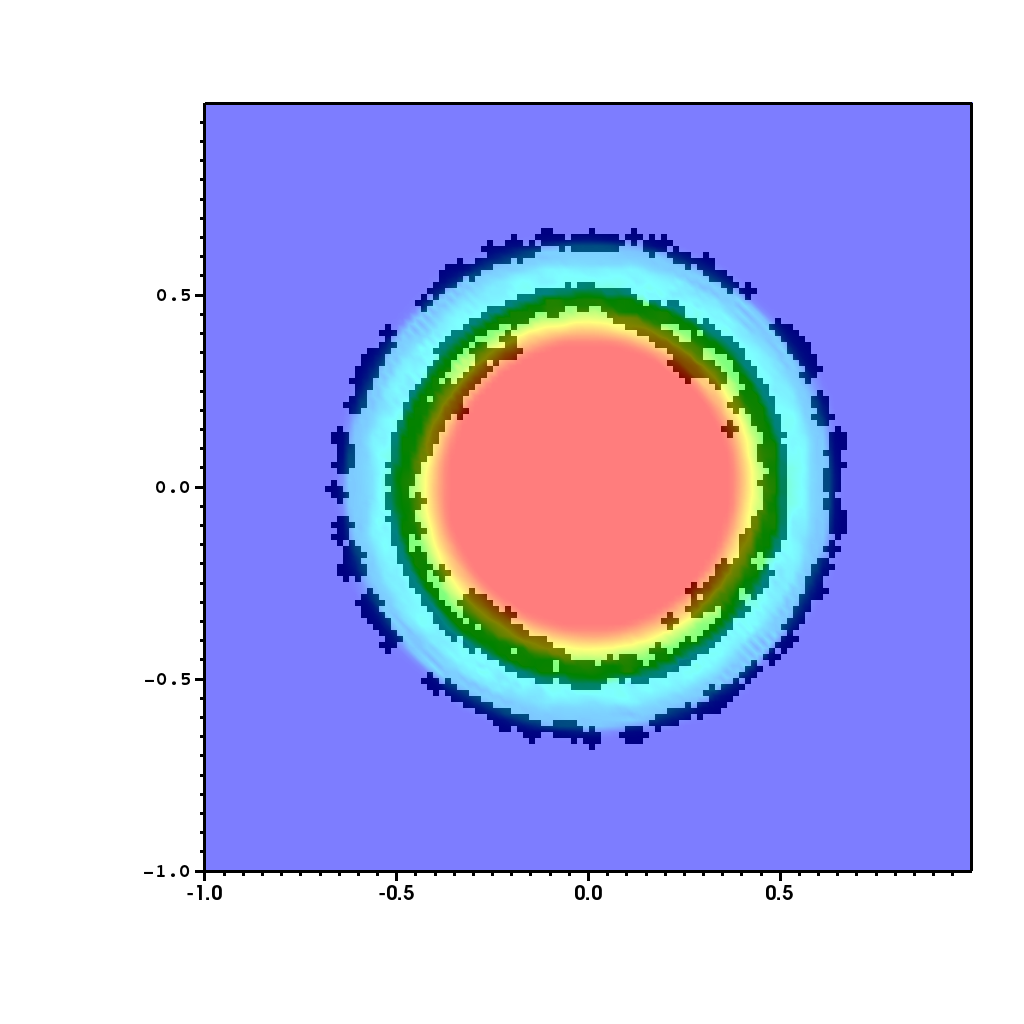}
\includegraphics[width=0.23\textwidth]{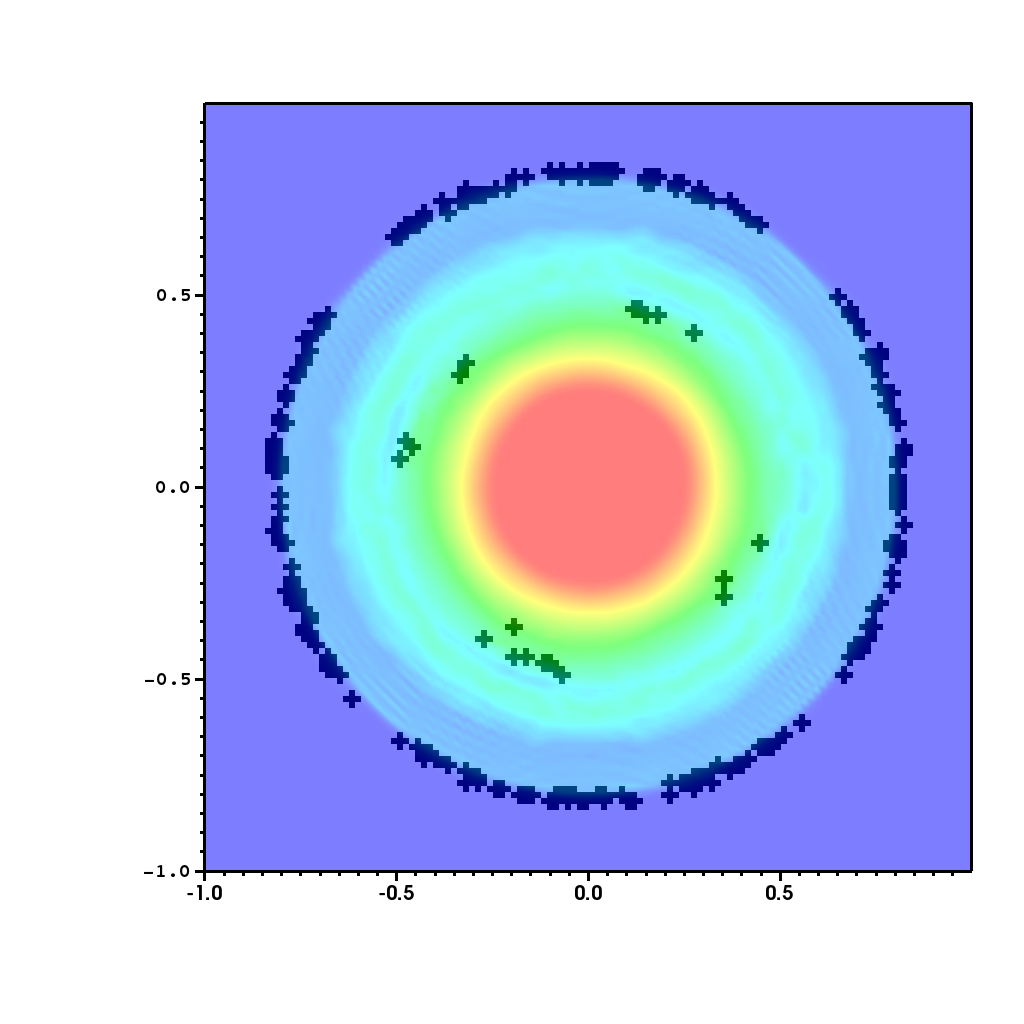}
\includegraphics[width=0.23\textwidth]{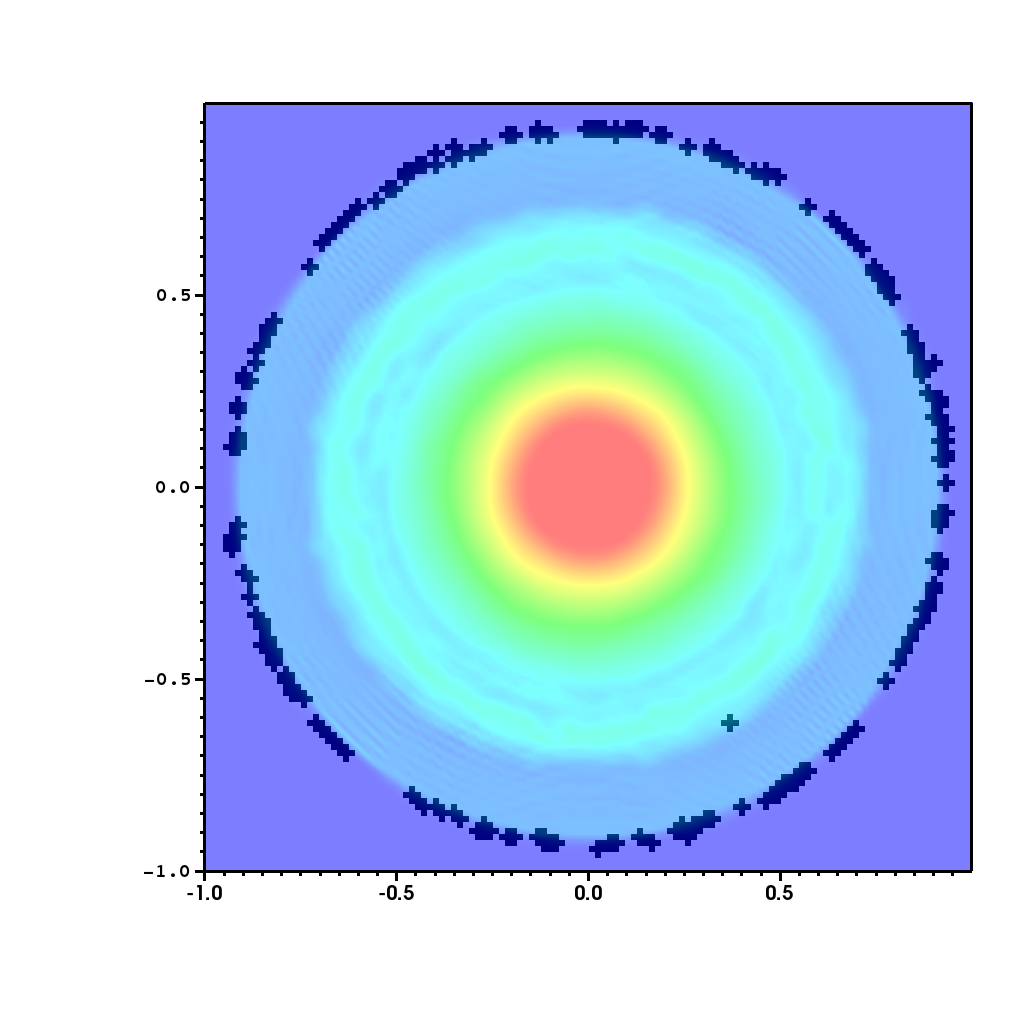}
\caption{Detection comparison between the neural network limiter without retraining (top row) and with retraining (bottom row) on a structured mesh.\label{fig:sod2ddetectioncomparisonquad}}}
\end{figure}

\begin{figure}
\centering{
\includegraphics[width=0.5\textwidth]{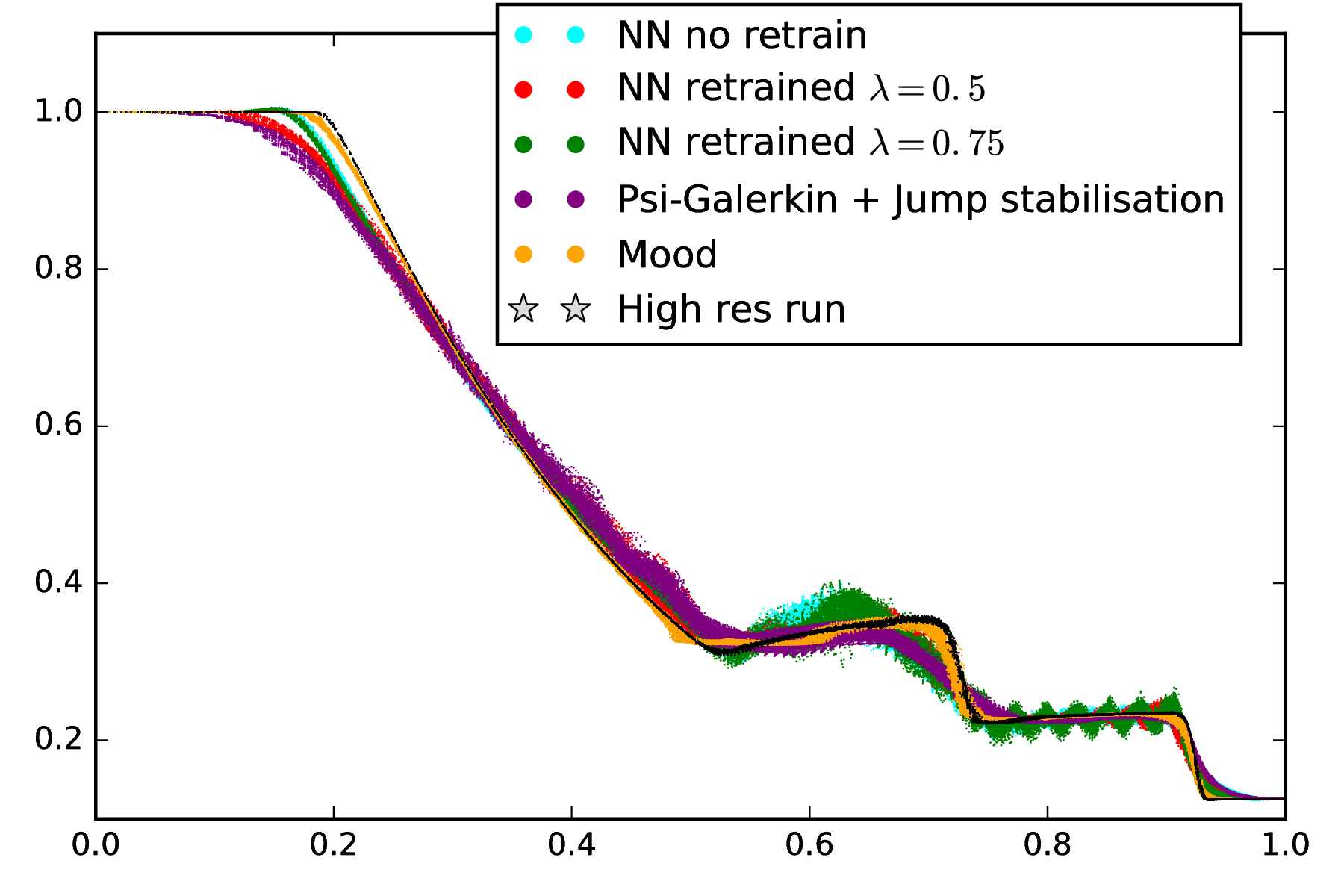}
\caption{Comparison between different limiters for the RD scheme on a structured mesh\label{fig:sod2comparisonquad}}}
\end{figure}

\begin{figure}
\centering{
\includegraphics[width=0.5\textwidth]{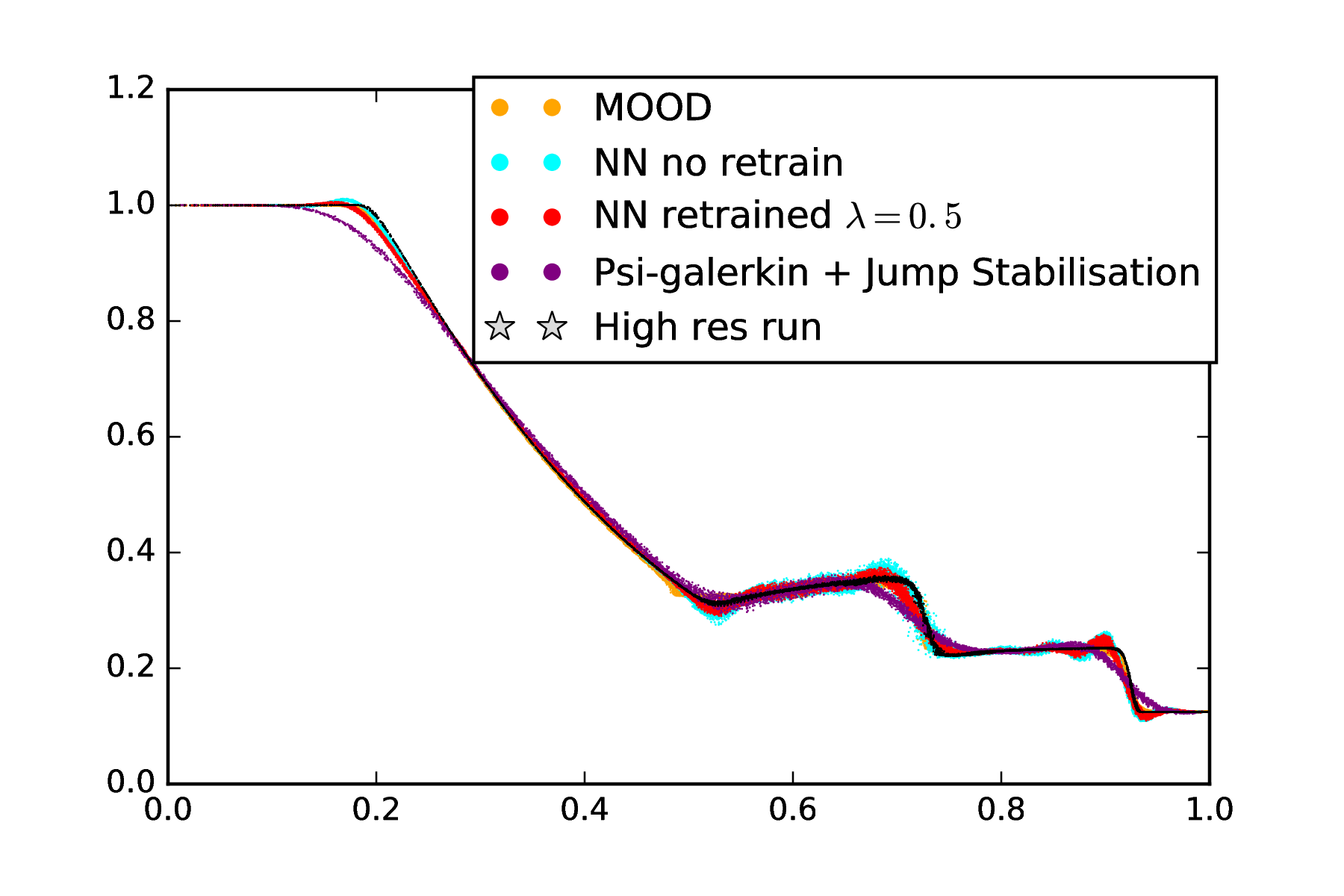}
\caption{Comparison between different limiters for the RD scheme on an unstructured mesh\label{fig:sod2comparisontria}}}
\end{figure}

In Figures \ref{fig:sod2ddetectioncomparisonquad} and \ref{fig:sod2ddetectioncomparisontrian}, we compare the performance of the transferred neural network without and with retraining on residual distribution data, on Cartesian grids and unstructured (triangular) grids, respectively. We can see that the adapted limiter detects troubled cells in the shock fronts, for both grid types, and that additionally, retraining the limiter on a small dataset using examples for the residual distribution scheme appears to make a difference, in particular for the unstructured mesh case. This is better observed in figures \ref{fig:sod2comparisonquad} and \ref{fig:sod2comparisontria}, where the performance of the neural network limiter (with varying retraining strategies) is compared to other limiters in the residual distribution context.

\section{Discussion and conclusion}
\label{sec:conclusion}
The purpose of this work was primarily to demonstrate the potential of using learning algorithms in CFD codes to automate away some parameter tuning which is typical when using stabilisation methods and to explore the idea of transferring knowledge across numerical schemes. In particular, we detail different stages necessary to train a black box limiter that can be integrated in different codes. To this end, we described how to construct a dataset for the training phase, how to train a multi-layer Perceptron that detects shocks and how to integrate the trained model with existing CFD codes. We performed the usual validation of the limiter in the context of scalar and systems of equations. Furthermore, we use the model trained on the data generated with the discontinuous Galerkin code and integrate this with a Residual distribution code, exploring different ways to perform domain adaptation.

We showed empirically that in principle it is possible to learn a \textit{parameter free} limiter (after the training phase) from data. While the performance for the seen functions is relatively good (as observed in \ref{sub:detection}), the on-the-fly performance could be improved with the careful design of the loss function during the training phase (for example, to include information on the maximum preserving principle). For the smooth case, we observed that this limiter was by far less diffuse than the $minmod$ limiter. Certain oscillations were corrected but there were other oscillations which were not stabilized enough. When applied to a system of equations, we observed that this method still produced sensible results, slightly better than the unlimited version.

It is to note that the training set used was extremely reduced - we only used data from solving the linear advection equation with different initial conditions and varying mesh sizes. Even so, the direct application of this model to unseen initial conditions, to the Euler system and to differing grid sizes was somewhat successful. In a realistic setting, one would use a much more complete dataset.

While we found that for the exact task of shock detection, a neural based limiter appears to have more several major drawbacks in comparison to some limiters which are both quite agnostic to the underlying numerical method and require minimal parameter tuning (e.g. MOOD \cite{mood}), it is our belief that these ideas can be applied to other problems which depend on certain local properties of the numerical solution, ultimately contributing towards Computational Fluid Dynamics (CFD) codes which are robust to different initial conditions and that require less parameter tuning to produce readily usable results.

The main issues we found were the following:
\begin{enumerate}
    \item{\textit{Lack of theoretical guarantees and quality of the numerical solution:}\\
    Throughout the time we worked on this problem, it became apparent that  certain properties of the shock capturing function do not arise without a certain amount of careful considerations. For example, it was particularly concerning the fact that the detection was not symmetric with respect to a defined stencil (see figure \ref{fig:invariance}). In order to overcome this problem, the literature typically suggests either a feature transformation which renders a certain property invariant or data augmentation by generating training examples that cover such invariance \cite{data_augmentation}.

\begin{figure}
\centering
%\begin{subfigure}{0.45\textwidth}
\includegraphics[width=0.4\linewidth]{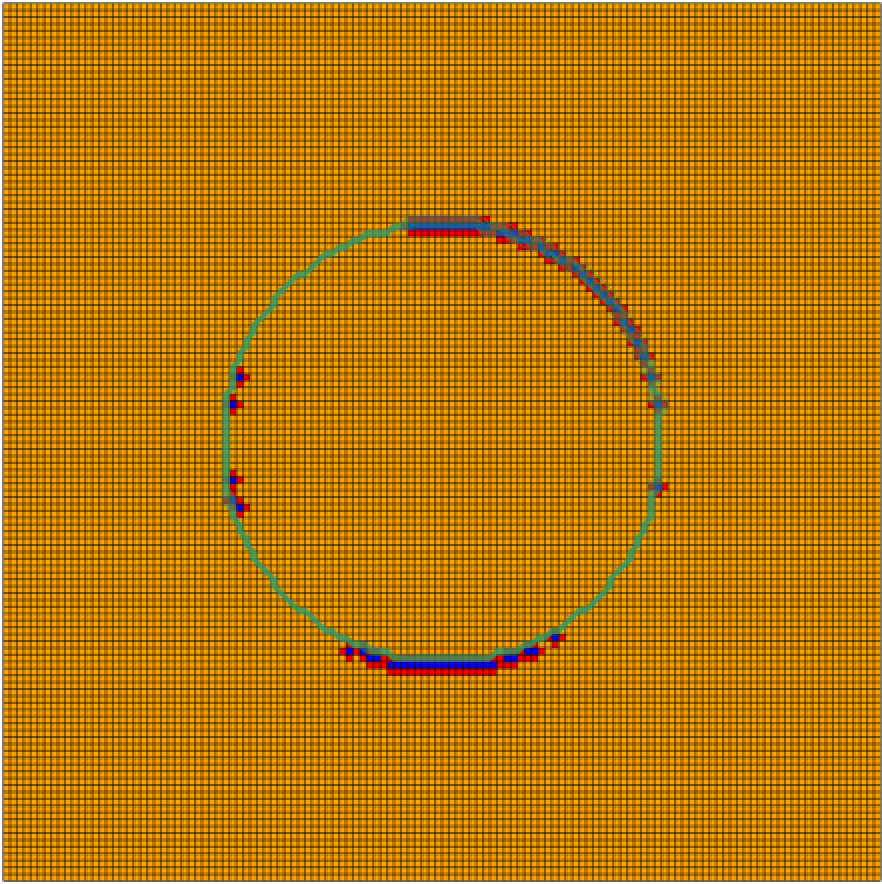}
%\caption{First subfigure} \label{fig:1a}
%\end{subfigure}
%\hspace*{\fill} % separation between the subfigures
%\begin{subfigure}{0.45\textwidth}
\includegraphics[width=0.4\linewidth]{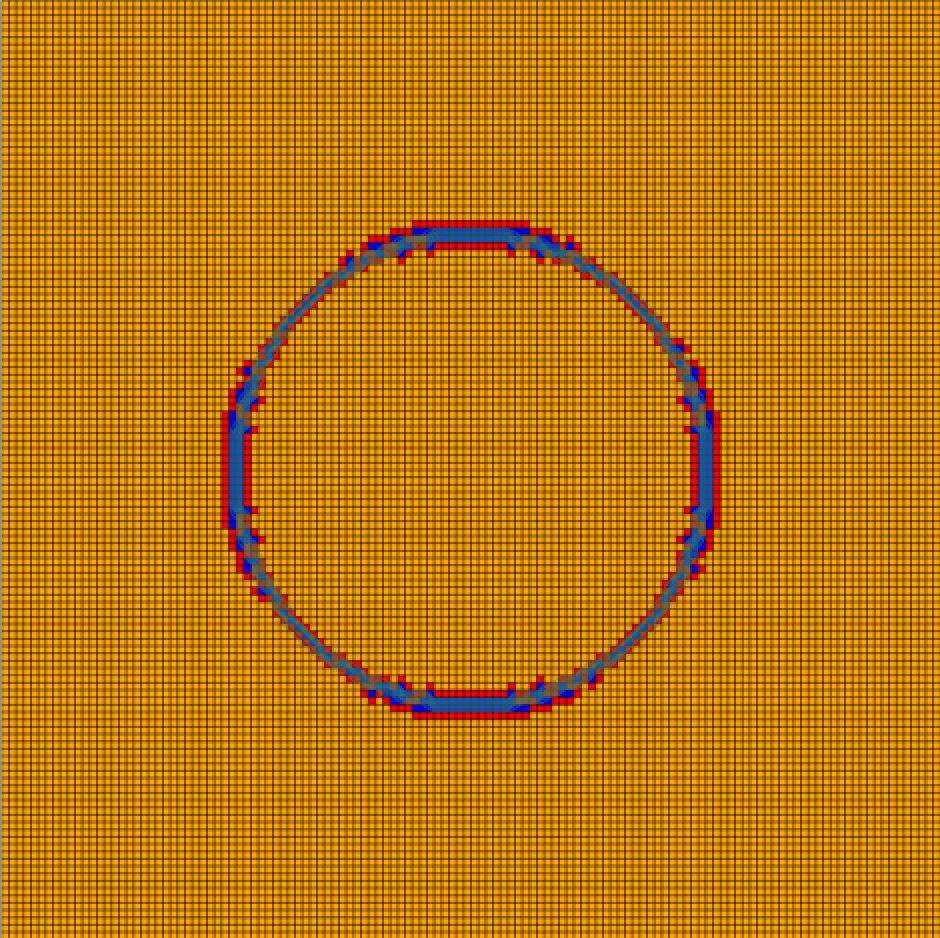}
%\caption{Second subfigure} \label{fig:1b}
%\end{subfigure}
\caption{No symmetry invariance considerations versus explicit symmetry invariance.} \label{fig:invariance}
\end{figure}

We empirically observed that merely performing data augmentation did not perform very well and it lead to longer training times (due to the enlargement of the dataset). We then used an ensemble classification where the stencil is permuted and the prediction is averaged (as detailed in section \ref{sec:directionalinv}), such that, for example, a prediction on a particular stencil and the same stencil mirrored along the x-axis yield the same response. Although \emph{hack} leads to exact invariance with respect to permutations of the neighbours, it adds a some cost to the method during the prediction phase.	
}
    \item{\textit{Computational performance degradation:}\\
    It has been claimed that using a neural network limiter can lead to reduced computational cost because the limiter is activated less times \cite{deepray2018}. However, this claim might have neglected the prediction step, which is necessary for each cell and thus, scales as $\mathcal{O}(N)$, $N$ the number of elements for a scalar problem. We verified that, when integrating this method with existing codes, the prediction step is expensive as it entails, at best, a series of matrix multiplications, and at worst, a feature vector computation (this can be stored, at the expense of larger memory requirements). Namely, we are required to compute $k$ matrix multiplications (where $k$ is the number of layers of the neural network) for each individual cell. This is not to say it is always going to be too computationally expensive - there have been some successful examples of learning a reduced network from a larger neural network, and with this, a large computational gain has been observed \cite{googlevoice}, however, this has not been in the scope of this project, and as such, we verified that the prediction step was overwhelmingly more expensive than applying a limiter, in particular, if the intention is to integrate this methodology with an existing code.    
    }
    \item{\textit{Lack of supporting theory on how to optimize, tune and generate deep neural networks:}\\
    To this day, the optimization of hyperparameters (such as architecture, shape, size, learning rates related to the neural network) is approached mostly through trial and error. While several works \cite{optimalbounds} estimate lower bounds and necessary complexity of the network to capture a given complexity of the function to be learned, these results remain rather far from concrete applications. However, in partice, there exist empirical ways to optimise the hyperparameter space \cite{hyperopt}.
    }
\end{enumerate}

The attractive property of this type of shock detector is that, once trained, it can operate as a black-box, parameter free shock detector. We have verified that both in one and two dimensional problems we were able to attain better results when comparing to limiters which have not been properly hand tuned. However, in the case of residual distribution, the transferred neural network limiter did not perform as well as the  MOOD limiter \cite{loubere,vilar,mood}.

\section{Appendix}
\label{appendix}

\subsection{MLP architectures}
\label{ap:arch}

We detail the architectures tested in this work in table \ref{tab:arch}. For this work we found the lower bound ($2^{28}$) for the number of weights and fixed the activation functions to be ReLU (as detailed in section \ref{sec:arch}). What we observed was that this quantity was not producing very good models. Thus, we used ($2^{32}$) non-zero weights, which presupposes a constant $C$ of order 10. It is left to specify the distribution of the weights across the different layers. The number of weights and neurons is related by multiplying $d$ through the number of neurons per layers to get the total number of weights.

\subsection{Adam algorithm}
\label{nn-ap:adam}
The Adam algorithm \cite{adam} is used for stochastic optimisation. The algorithm updates exponential moving averages of the gradient ($m_t$) and the squared gradients ($v_t$) where the hyper-parameters $\beta_1, \beta_2 \in [0,1)$ control the exponential decay of these moving averages. These moving averages are estimates of the first and second moment (uncentered variance) of the gradient. Because the quantities $m_0, v_0$ are initialised at 0, they are biased towards zero (in particular if the decaying rates $\beta_1$ and $\beta_2$ are small). Hence, there is an additional correction leading to the unbiased quantities $\hat{m}_t$ and $\hat{v}_t$. The algorithm is detailed in \ref{algo:adam}.

\begin{table}[h!]
\centering
%\begin{minipage}{\textwidth}
\begin{tabular}{ c | c c c }
      \hline
      \textbf{Model} & \textbf{Layers}& \textbf{Neurons per layer} & \textbf{Description} \\ \hline\hline
      Model 1 & 2 & 16384:16384\footnote{memory requirements were too high when keeping the same amount of weights (because in this case, weights = neurons), so this model was reduced.} & --- \\ \hline
      Model 2 & 3 & 2048:1024:512 & --- \\ \hline
      Model 3 & 4 & 512:256:256:128 & --- \\ \hline
      Model 4 & 5 & 256:128:64:64:32 & --- \\ \hline
      Model 5 & 5 & 256:128:64:64:32 & with weighted loss $\omega = 5$. \\ \hline
    \end{tabular}
    \caption{Architectures \label{tab:arch}}
    %\end{minipage}
\end{table}

\begin{algorithm}[h!]
 \KwData{$\alpha$: stepsize, $\beta_1$, $\beta_2 \in [0,1]$: exponential decay rates for the moment estimates, $f(\theta)$: stochastic objective function with parameters $\theta$, $\theta_0:$ initial parameter vector}
 \KwResult{ $\theta_t$ (resulting parameters) }
 $m_0$ = 0 (initialise 1st moment vector)\;
 $v_0$ = 0 (initialise 2nd moment vector)\;
 $t$ = 0 (initialise timestep)\;
 \While{$\theta_t$ not converged}{
  $t = t +1$\;
  $g_t = \nabla_\theta f_t (\theta_{t-1})$ (Get gradients w.r.t stochastic objective at stepsize $t$)\;
  $m_t = \beta_1 \cdot m_{t-1} + (1-\beta_1)\cdot g_t$ (Update biased first moment estimate)\;
  $v_t = \beta_2 \cdot v_{t-1} + (1-\beta_2)\cdot g_t^2$ (Update biased second moment estimate)\;
  $\hat{m}_t = m_t/(1-\beta_1^t)$ (Biased corrected first moment estimate)\;
  $\hat{v}_t = v_t/(1-\beta_2^t)$ (Biased corrected second raw moment estimate)\;
  $\theta_t = \theta_{t-1} - \alpha\cdot\hat{m}_t/\sqrt[]{\hat{v}_t +\epsilon}$ (Update parameters)\;
  }
 \caption{Adam algorithm \label{algo:adam}}
\end{algorithm}

\begin{acknowledgements}
MHV has been funded by the UZH Candoc Forschungskredit grant.
\end{acknowledgements}

\section*{Conflict of interest}
On behalf of all authors, the corresponding author states that there is no conflict of interest.

\bibliographystyle{spmpsci}      % mathematics and physical sciences
\bibliography{nnpaper}   % name your BibTeX data base

\begin{thebibliography}{10}
\providecommand{\url}[1]{{#1}}
\providecommand{\urlprefix}{URL }
\expandafter\ifx\csname urlstyle\endcsname\relax
  \providecommand{\doi}[1]{DOI~\discretionary{}{}{}#1}\else
  \providecommand{\doi}{DOI~\discretionary{}{}{}\begingroup
  \urlstyle{rm}\Url}\fi

\bibitem{repository}
Github repository.
\newblock \url{https://github.com/hanveiga/1d-dg-nn}

\bibitem{abg}
Abgrall, R.: Residual distribution schemes: Current status and future trends.
\newblock Computers \& Fluids \textbf{35}(7), 641 -- 669 (2006).
\newblock \doi{https://doi.org/10.1016/j.compfluid.2005.01.007}.
\newblock
  \urlprefix\url{http://www.sciencedirect.com/science/article/pii/S0045793005001738}.
\newblock Special Issue Dedicated to Professor Stanley G. Rubin on the Occasion
  of his 65th Birthday

\bibitem{abgrall2019}
Abgrall, R., Bacigaluppi, P., Tokareva, S.: High-order residual distribution
  scheme for the time-dependent euler equations of fluid dynamics.
\newblock Computers \& Mathematics with Applications \textbf{78}(2), 274 -- 297
  (2019).
\newblock \doi{https://doi.org/10.1016/j.camwa.2018.05.009}.
\newblock
  \urlprefix\url{http://www.sciencedirect.com/science/article/pii/S0898122118302712}.
\newblock Proceedings of the Eight International Conference on Numerical
  Methods for Multi-Material Fluid Flows (MULTIMAT 2017)

\bibitem{tvdlimiter}
Arora, M., Roe, P.L.: A well-behaved tvd limiter for high-resolution
  calculations of unsteady flow.
\newblock Journal of Computational Physics \textbf{132}(1), 3 -- 11 (1997).
\newblock \doi{https://doi.org/10.1006/jcph.1996.5514}.
\newblock
  \urlprefix\url{http://www.sciencedirect.com/science/article/pii/S002199919695514X}

\bibitem{mood}
{Bacigaluppi}, P., {Abgrall}, R., {Tokareva}, S.: {``A Posteriori'' Limited
  High Order and Robust Residual Distribution Schemes for Transient Simulations
  of Fluid Flows in Gas Dynamics}.
\newblock arXiv e-prints arXiv:1902.07773 (2019)

\bibitem{hyperopt}
Bergstra, J., Yamins, D., Cox, D.: Making a science of model search:
  Hyperparameter optimization in hundreds of dimensions for vision
  architectures.
\newblock In: S.~Dasgupta, D.~McAllester (eds.) Proceedings of the 30th
  International Conference on Machine Learning, \emph{Proceedings of Machine
  Learning Research}, vol.~28, pp. 115--123. PMLR, Atlanta, Georgia, USA
  (2013).
\newblock \urlprefix\url{http://proceedings.mlr.press/v28/bergstra13.html}

\bibitem{biswas1994}
Biswas, R., Devine, K.D., Flaherty, J.E.: Parallel, adaptive finite element
  methods for conservation laws.
\newblock Applied Numerical Mathematics \textbf{14}(1), 255 -- 283 (1994).
\newblock \doi{https://doi.org/10.1016/0168-9274(94)90029-9}.
\newblock
  \urlprefix\url{http://www.sciencedirect.com/science/article/pii/0168927494900299}

\bibitem{loubere}
Clain, S., Diot, S., Loub{\`e}re, R.: Multi-dimensional optimal order detection
  (mood) --- a very high-order finite volume scheme for conservation laws on
  unstructured meshes.
\newblock In: J.~Fo{\v{r}}t, J.~F{\"u}rst, J.~Halama, R.~Herbin, F.~Hubert
  (eds.) Finite Volumes for Complex Applications VI Problems {\&} Perspectives,
  pp. 263--271. Springer Berlin Heidelberg, Berlin, Heidelberg (2011)

\bibitem{backprop}
Feldman, J., Rojas, R.: Neural Networks: A Systematic Introduction.
\newblock Springer Berlin Heidelberg (1996).
\newblock \urlprefix\url{https://books.google.ch/books?id=txsjjYzFJS4C}

\bibitem{deeplearning}
Goodfellow, I., Bengio, Y., Courville, A.: Deep Learning.
\newblock MIT Press (2016).
\newblock \url{http://www.deeplearningbook.org}

\bibitem{googlevoice}
{He}, Y., {Sainath}, T.N., {Prabhavalkar}, R., {McGraw}, I., {Alvarez}, R.,
  {Zhao}, D., {Rybach}, D., {Kannan}, A., {Wu}, Y., {Pang}, R., {Liang}, Q.,
  {Bhatia}, D., {Shangguan}, Y., {Li}, B., {Pundak}, G., {Sim}, K.C., {Bagby},
  T., {Chang}, S., {Rao}, K., {Gruenstein}, A.: Streaming end-to-end speech
  recognition for mobile devices.
\newblock In: ICASSP 2019 - 2019 IEEE International Conference on Acoustics,
  Speech and Signal Processing (ICASSP), pp. 6381--6385 (2019).
\newblock \doi{10.1109/ICASSP.2019.8682336}

\bibitem{imbalance}
Hoens, T.R., Chawla, N.V.: Imbalanced Datasets: From Sampling to Classifiers.
\newblock John Wiley \& Sons, Inc. (2013).
\newblock \doi{10.1002/9781118646106.ch3}

\bibitem{adam}
Kingma, D.P., Ba, J.: Adam: {A} method for stochastic optimization.
\newblock CoRR \textbf{abs/1412.6980} (2014).
\newblock \urlprefix\url{http://arxiv.org/abs/1412.6980}

\bibitem{lydia}
{Krivodonova}, L.: {Limiters for high-order discontinuous Galerkin methods}.
\newblock Journal of Computational Physics \textbf{226}, 879--896 (2007).
\newblock \doi{10.1016/j.jcp.2007.05.011}

\bibitem{vision}
Krizhevsky, A., Sutskever, I., Hinton, G.E.: Imagenet classification with deep
  convolutional neural networks.
\newblock In: Proceedings of the 25th International Conference on Neural
  Information Processing Systems - Volume 1, NIPS'12, pp. 1097--1105. Curran
  Associates Inc., USA (2012).
\newblock \urlprefix\url{http://dl.acm.org/citation.cfm?id=2999134.2999257}

\bibitem{kurganov}
Kurganov, A., Tadmor, E.: Solution of two-dimensional riemann problems for gas
  dynamics without riemann problem solvers.
\newblock Numerical Methods for Partial Differential Equations \textbf{18}(5),
  584--608 (2002).
\newblock \doi{10.1002/num.10025}.
\newblock
  \urlprefix\url{https://onlinelibrary.wiley.com/doi/abs/10.1002/num.10025}

\bibitem{data_augmentation}
Mikolajczyk, A., Grochowski, M.: Data augmentation for improving deep learning
  in image classification problem.
\newblock 2018 International Interdisciplinary PhD Workshop (IIPhDW) pp.
  117--122 (2018)

\bibitem{optimalbounds}
{Petersen}, P., {Voigtlaender}, F.: {Optimal approximation of piecewise smooth
  functions using deep ReLU neural networks}.
\newblock ArXiv e-prints  (2017)

\bibitem{earlystopping}
Prechelt, L.: Early Stopping --- But When?, pp. 53--67.
\newblock Springer Berlin Heidelberg, Berlin, Heidelberg (2012).
\newblock \doi{10.1007/978-3-642-35289-8_5}.
\newblock \urlprefix\url{https://doi.org/10.1007/978-3-642-35289-8_5}

\bibitem{deepray2018}
Ray, D., Hesthaven, J.S.: An artificial neural network as a troubled-cell
  indicator.
\newblock Journal of Computational Physics \textbf{367}, 166 -- 191 (2018).
\newblock \doi{https://doi.org/10.1016/j.jcp.2018.04.029}.
\newblock
  \urlprefix\url{http://www.sciencedirect.com/science/article/pii/S0021999118302547}

\bibitem{ricchiuto2010}
Ricchiuto, M., Abgrall, R.: Explicit runge–kutta residual distribution
  schemes for time dependent problems: Second order case.
\newblock Journal of Computational Physics \textbf{229}(16), 5653 -- 5691
  (2010).
\newblock \doi{https://doi.org/10.1016/j.jcp.2010.04.002}.
\newblock
  \urlprefix\url{http://www.sciencedirect.com/science/article/pii/S0021999110001786}

\bibitem{SWjcp}
Ricchiuto, M., Abgrall, R., Deconinck, H.: Application of conservative residual
  distribution schemes to the solution of the shallow water equations on
  unstructured meshes.
\newblock J. Comput. Phys. \textbf{222}(1), 287--331 (2007).
\newblock \doi{10.1016/j.jcp.2006.06.024}.
\newblock \urlprefix\url{http://dx.doi.org/10.1016/j.jcp.2006.06.024}

\bibitem{rojas2}
Rojas, R.: Networks of width one are universal classifiers.
\newblock In: Proceedings of the International Joint Conference on Neural
  Networks, 2003., vol.~4, pp. 3124--3127 vol.4 (2003).
\newblock \doi{10.1109/IJCNN.2003.1224071}

\bibitem{rojas}
{Rojas}, R.: {Deepest Neural Networks}.
\newblock ArXiv e-prints  (2017)

\bibitem{gradientdescent}
Snyman, J.: Practical Mathematical Optimization: An Introduction to Basic
  Optimization Theory and Classical and New Gradient-Based Algorithms.
\newblock Applied Optimization. Springer (2005).
\newblock \urlprefix\url{https://books.google.ch/books?id=0tFmf\_UKl7oC}

\bibitem{ilya}
Sutskever, I., Hinton, G.E.: Deep, narrow sigmoid belief networks are universal
  approximators.
\newblock Neural Comput. \textbf{20}(11), 2629--2636 (2008).
\newblock \doi{10.1162/neco.2008.12-07-661}.
\newblock \urlprefix\url{http://dx.doi.org/10.1162/neco.2008.12-07-661}

\bibitem{s2s}
Sutskever, I., Vinyals, O., Le, Q.V.: Sequence to sequence learning with neural
  networks.
\newblock In: Proceedings of the 27th International Conference on Neural
  Information Processing Systems - Volume 2, NIPS'14, pp. 3104--3112. MIT
  Press, Cambridge, MA, USA (2014).
\newblock \urlprefix\url{http://dl.acm.org/citation.cfm?id=2969033.2969173}

\bibitem{mhvabgrall}
Veiga, M.H., Abgrall, R.: Towards a general stabilisation method for
  conservation laws using a multilayer perceptron neural network: 1d scalar and
  system of equations.
\newblock In: European Conference on Computational Mechanics and VII European
  Conference on Computational Fluid Dynamics, 1, pp. 2525--2550. ECCM (2018).
\newblock \urlprefix\url{https://doi.org/10.5167/uzh-168538}

\bibitem{vilar}
Vilar, F.: A posteriori correction of high-order discontinuous galerkin scheme
  through subcell finite volume formulation and flux reconstruction.
\newblock Journal of Computational Physics \textbf{387}, 245 -- 279 (2019).
\newblock \doi{https://doi.org/10.1016/j.jcp.2018.10.050}.
\newblock
  \urlprefix\url{http://www.sciencedirect.com/science/article/pii/S0021999118307174}

\bibitem{transfer}
Weiss, K.R., Khoshgoftaar, T.M., Wang, D.: A survey of transfer learning.
\newblock Journal of Big Data \textbf{3}, 1--40 (2016)

\bibitem{relus}
Xu, B., Wang, N., Chen, T., Li, M.: Empirical evaluation of rectified
  activations in convolutional network.
\newblock CoRR \textbf{abs/1505.00853} (2015).
\newblock \urlprefix\url{http://arxiv.org/abs/1505.00853}

\end{thebibliography}

\end{document}